\documentclass[reqno,a4paper,11pt]{amsart}
\usepackage{amsmath,amsthm,amssymb,mathrsfs}
\usepackage[
left = 2.25cm,
right = 2.25cm,
top = 2.2 cm,
bottom = 2.2 cm
]{geometry}
\usepackage{graphicx,xcolor,tikz,mathtools}
\usepackage{bm,stmaryrd}
\usepackage{amsaddr}
\usepackage{enumitem}
\usepackage[toc,page]{appendix}
\usepackage[
pdfencoding=auto,
colorlinks = true,
citecolor = blue,
linkcolor = blue,
anchorcolor = blue,
urlcolor = blue, 
filecolor=blue,
pdfkeywords={}
]{hyperref}
\usepackage{esint} 

\usepackage{cite}
\usepackage{todonotes}
\makeatletter
\define@key{todonotes}{Helmut}[]{%
	\setkeys{todonotes}{author=\textbf{Helmut},inline,color=red!10}}%
\define@key{todonotes}{Andrea}[]{%
	\setkeys{todonotes}{author=\textbf{Andrea},inline,color=green!20}}%
\define@key{todonotes}{Yadong}[]{%
	\setkeys{todonotes}{author=\textbf{Yadong},inline,color=cyan!20}}%
\makeatother

\theoremstyle{plain}
\newtheorem{theorem}{Theorem}[section]

\newtheorem{definition}[theorem]{Definition}
\newtheorem{lemma}[theorem]{Lemma}
\newtheorem{proposition}[theorem]{Proposition}
\newtheorem{remark}[theorem]{Remark}
\theoremstyle{definition}

\numberwithin{equation}{section}

\newcommand{\bu}{\mathbf{u}}
\newcommand{\bv}{\mathbf{v}}
\newcommand{\bw}{\mathbf{w}}
\newcommand{\bn}{\mathbf{n}}

\newcommand{\bA}{\mathbf{A}}
\newcommand{\bH}{\mathbf{H}}

\newcommand{\bJ}{\mathbf{J}}
\newcommand{\bP}{\mathbf{P}}

\newcommand{\bW}{\mathbf{W}}

\newcommand{\cA}{\mathcal{A}}

\newcommand{\cE}{\mathcal{E}}

\newcommand{\cH}{\mathcal{H}}

\newcommand{\cL}{\mathcal{L}}

\renewcommand{\d}{\mathrm{d}}
\newcommand{\dx}{\,\d x}
\newcommand{\dt}{\,\d t}
\newcommand{\dxdt}{\,\d x \d t}

\newcommand{\ddt}{\frac{\d}{\d t}}

\newcommand{\ptial}[1]{ \partial_{#1} }
\newcommand{\pt}{\ptial{t}}
\newcommand{\onehalf}{\frac{1}{2}}

\newcommand{\norm}[1]{\left\Vert #1 \right \Vert}
\newcommand{\normm}[1]{\Vert #1 \Vert}
\newcommand{\normmm}[1]{\left\vert #1 \right\vert}

\def\Div{\mathrm{div}\,}
\DeclareMathOperator*{\esssup}{\mathrm{esssup}}

\DeclarePairedDelimiter\ceil{\lceil}{\rceil}
\DeclarePairedDelimiter\floor{\lfloor}{\rfloor}
\def\xhi{\chi}
\def\tTh{_t^{T_h(t)}}
\def\Th#1{{T_h(#1)}}
\def\Ths{\Th{s}}
\def\Tht{\Th{t}}
\def\Mm{\mathcal M_{m_0}}
\def\Hmz{H^{-1}_{(0)}(\Omega)}
\def\Hz{H^{1}_{(0)}(\Omega)}
\def\normh#1{\norm{#1}_{\Hmz}}

\def\Hus{\mathbf H^1_{0,\sigma}(\Omega)}
\def\Hs{\mathbf H^1_{\sigma}(\Omega)}

\def\bxi{\boldsymbol\xi}
\def\bL{\mathbf L}

\def\Ts{{T^*}}
\def\tF{\widetilde{F}}
\def\wchi{\widehat\chi}
\def\tmup{\normmm{\mu_t^{\partial\Omega}}_{\mathbb{S}^{d-1}}}
\def\tmu{\normmm{\mu_t^{\Omega}}_{\mathbb{S}^{d-1}}}
\def\om{\Omega}
\def\pom{{\partial\Omega}}
\def\ov{\overline}
\def\wtilde{\widetilde}
\def\R{\mathbb R}
\def\Ls{\bL^2_\sigma(\Omega)}
\def\bpsi{\boldsymbol\psi}
\def\N{\mathbb N}
\def\R{\mathbb R}
\def\Hds{\mathfrak{D}(\bA)}

\def\Grad{\nabla}
\def\cA{\mathcal A}
\def\pcA{\partial\mathcal A}
\def\pAs{\partial^*\cA}
\def\ovu{\overline{u}}

\def\tw{\widetilde{w}}
\def\tlambda{\widetilde{\lambda}}
\def\uuu{\mathbf v}
\def\vphi{\varphi}
\def\JJJ{\widetilde{\mathbf J}}

\allowdisplaybreaks[4]

\begin{document}
	
	\title[Weak Solutions to a Two-Phase Flow of Incompressible Viscous Fluids]{Weak Solutions to a Sharp Interface Model for a Two-Phase Flow of Incompressible Viscous Fluids with Different Densities}
	
	\author[H. Abels and A. Poiatti]{
		\small
		Helmut Abels$^\ast$  and 
		Andrea Poiatti$^\dagger$
	}
	
	\address{
		$^\ast$Fakult\"at f\"ur Mathematik,
		Universit\"at Regensburg,
		93053 Regensburg, Germany
	}
	\email{helmut.abels@ur.de}
	
	\address{
		$^\dagger$Faculty of Mathematics, University of Vienna, 
        1090 Vienna, Austria
	}
	\email{andrea.poiatti@univie.ac.at}

	
	%
	
		
	
	\maketitle
	
	\begin{abstract}
	In this paper we consider the flow of two incompressible, viscous and immiscible fluids in a bounded domain, with different densities and viscosities. This model consists of a coupled system of Navier-Stokes and Mullins-Sekerka type parts, and can be obtained from the sharp interface limit of the diffuse interface model proposed by the first author, Garcke, and Gr\"{u}n (Math. Models Methods Appl. Sci. 22, 2012). We introduce a new notion of weak solutions and prove its global in time existence, together with a consistency result of smooth weak solutions with the classical Navier-Stokes-Mullins-Sekerka system. Our new notion of solution allows to include the case of different densities of the two fluids, a sharp energy dissipation principle à la De Giorgi, together with a weak formulation of the constant contact angle condition at the boundary, which were left open in the previous notion of solution proposed by the first author and R\"{o}ger (Ann. Inst. H. Poincaré Anal. Non Linéaire 26, 2009). 
	\end{abstract}
	\section{Introduction} 

In this paper we study the flow of two incompressible, viscous and immiscible fluids in a bounded domain $\Omega$, $d = 2, 3$, with smooth boundary $\pom$. We assume the fluids to fill domains $\Omega^+(t)$ and $\Omega^-(t)$, $t>0$, respectively, with a shared interface $\Gamma(t)$ between them. Here the surface $\Gamma(t)$ can also touch the boundary $\pom$ with a fixed constant contact angle $\gamma>0$. The flow is described by means of the velocity $\bv : \om \times(0,\infty)\to \R^d$, and the pressure $\pi : \om \times (0,\infty)\to\R$ in both fluids (in Eulerian coordinates). We typically assume the fluids to be Newtonian, so that the stress tensors are of the form $\nu_{i}D\bv - \pi Id$, $i=1,2$, with constant viscosities $\nu_1,\nu_2>0$ corresponding to the viscosities of the fluids filling $\om^+(t)$ and $\om^-(t)$, respectively. Here $D\bv=\tfrac{\nabla\bv+(\nabla\bv)^T}{2}$ represents the symmetric strain rate tensor. We consider the case in which the surface tension at the interface is present, and we also assume that the two fluids have different (constant) densities $\rho_1>0$ and $\rho_2>0$, respectively.  For the evolution of the two phases we account for diffusional effects and consider a contribution to the flux proportional to the negative gradient of the chemical potential $\overline u$. Due to the different densities of the two fluids, the model accounts for the presence of an extra diffusional flux $\wtilde{\mathbf J}$, proportional to the gradient of $\overline u$, in the velocity equation. We now introduce some notation and then present the complete model. First consider a finite time horizon $\Ts\in(0,\infty)$ and let $\mathcal A =
	(\mathcal A (t))_t\in[0,\Ts)$ be a time-dependent family of smoothly evolving open subsets $\mathcal A (t) \subset \om$ with $\partial\mathcal A(t) = \partial^*\mathcal A(t)$, $t \in [0, \Ts)$,
	where $\partial^*\mathcal A$ is the reduced boundary of $\cA$. We set $\cA_0:=\cA(0)$ and denote by $\xhi_0$ the characteristic function corresponding to the initial smooth set $\cA_0$. In the fluid setting introduced above, we have $\om^+(t)=\cA(t)$, $\om^-(t)=\om\setminus \overline{\cA(t)}$, and $\Gamma(t)={\partial \cA(t)\cap \Omega}$, so that $\om\setminus \pcA(t)=\om^+(t)\cup \om^-(t)$.

    Then define, for every $t \in [0, \Ts)$,
 $V_{\pcA} (t)$ and $\bH_{\pcA (t)}=H_{\pcA(t)}\bn_{\pcA(t)}$ as the associated normal velocity and mean curvature vector with respect to $\bn_{\pcA(t)}$, respectively, whereas $H_{\pcA(t)}$ is the scalar mean curvature. Note that $\bn_{\pcA(t)}$ is assumed to be the
 unit normal vector field along $\pcA(t)$ pointing inside $\om^+(t)=\cA (t)$, and analogously
 $\bn_{\pom}$ is the interior unit normal of the domain boundary $\pom$.  Then denote by $\xhi(t)$ the characteristic function associated to $\cA(t)$, i.e., to the domain $\om^+(t)$ filled by the fluid 1. We say that the family $\cA$ and the velocity $\bv$  evolve by the Navier-Stokes-Mullins-Sekerka system if for each
	$t \in (0, \Ts)$ there exists a chemical potential $\ovu(·, t)$ so that
	\begin{alignat}{2}
		\label{A1}\Delta \ovu(\cdot,t)&=0,&&\text{ in }\om\setminus\pcA(t),\\\label{A2}
		V_{\pcA(t)}&=(\bv\cdot\bn_{\pcA(t)})\bn_{\pcA(t)}-(\bn_{\pcA}\cdot[[\nabla\ovu(\cdot,t)]])\bn_{\pcA(t)}\quad\quad\quad&&\text{ on }\pcA(t)\cap \om,\\\label{A3}
		\ovu(\cdot,t)&=c_0 H_{\pcA(t)},&&\text{ on }\pcA(t)\cap \om,\\\label{A4}
		(\bn_{\pcA(t)}\cdot\nabla )\ovu(\cdot,t)&=0,&& \text{ on }\pom\setminus \overline{\pcA(t)\cap \om},\\\label{A5}
		\bn_{\pom}\cdot\bn_{\pcA(t)}&=\cos\gamma,&&\text{ on }\pom\cap \overline{\pcA(t)\cap \om}.
	\end{alignat}
	and the velocity $\bv$ satisfies
	\begin{alignat}{2}
\pt(\rho(\xhi)\bv)+\Div(\rho(\chi)\bv\otimes \bv+\bv\otimes \widetilde{\mathbf J})-\Div(\nu(\chi)D\bv)+\nabla \pi&=\mathbf 0,&& \text{ in }\om\setminus \pcA(t),\label{A6}\\\label{A7}
		\Div \bv&=0,&&\text{ in }\om\setminus \pcA(t),\\\label{A8}
		-[[\nu(\xhi)D\bv]]\bn_{\pcA(t)}+[[p]]\bn_{\pcA(t)}&=c_0\bH_{\pcA(t)},\qquad&&\text{ on }\pcA(t)\cap 	\om,\\\label{A9}
		\bv&=\mathbf 0&& \text{ on }\pom,\\\label{A10}
		\bv_{|t=0}&=\bv_0,&&\text{ in }\om,
	\end{alignat}
	for any $t\in[0,\Ts)$, where $\gamma\in (0,\pi/2]$ is the (constant) contact angle, and $c_0>0$ is the fixed surface tension constant. Moreover, the jump $[[\cdot]]$ across the interface $\Gamma(t)=\pcA (t)\cap \om$ is taken for a function $f$ as
	\begin{align*}
		[[f]](x):=\lim_{h	\to0^+}(f (x +h\bn_{\pcA(t)})-f (x -h\bn_{\pcA(t)})),\quad \forall x\in \pcA(t).
	\end{align*}
	Then, we have set  $\rho(\xhi)=\rho_1\xhi+\rho_2(1-\xhi)\in \{\rho_1,\rho_2\}$ and $\nu(\xhi)=\nu_1\xhi+\nu_2(1-\xhi)\in\{\nu_1,\nu_2\}$. As already introduced, $\widetilde{\mathbf J}:=-(\rho_1-\rho_2)\nabla\ovu$ is an additional diffusive flux due to the different densities between the two fluids, which is specific for this model. 

Equations \eqref{A6}-\eqref{A7} model the conservation of linear momentum and mass in the two fluids and \eqref{A8} is the balance of normal stresses on the interface $\Gamma(t)$. Then the equations for the velocity $\bv$ are completed by the no-slip condition \eqref{A9} at the boundary $\partial\om$.
On the other hand, conditions \eqref{A1}-\eqref{A5} describe a continuity equation for the mass of the phases, and
connect the chemical potential $\overline u$ to the $L^2$-gradient of the surface area, given by the mean curvature of the separating interface $\Gamma(t)$. They also account for a constant contact angle $\gamma$ between the closure of $\Gamma(t)$ and $\pom$.
Observe that equations \eqref{A1}-\eqref{A5} with $\bv=\mathbf 0$ define the well-known Mullins-Sekerka flow. The main property of this evolution, which will be heavily exploited in the sequel, is that it describes the gradient flow for the surface area functional with
respect to the $\Hmz$ inner product. 

Model \eqref{A1}-\eqref{A10} can be obtained as the sharp interface limit ($\varepsilon\to0$) of the following
diffuse interface model, with constant mobility set equal to 1, first introduced in \cite{AGG} as an extension to the unmatched densities case of the celebrated model H (see \cite{Gurtin, ModelH}):
\begin{equation}
\begin{cases}
\partial _{t}(\rho (\vphi )\uuu )+\mathrm{div}\,\left( \uuu \otimes
\left( \rho (\vphi )\uuu +\JJJ\right) \right) -\mathrm{div}\,(\nu
(\vphi )D\uuu )+\nabla \pi =-\varepsilon\ \mathrm{div}\left( \nabla \vphi \otimes
\nabla \vphi \right) , \\
\mathrm{div}\,\uuu =0, \\
\partial _{t}\vphi +\uuu \cdot \nabla \vphi =\Delta\overline u, \\
\overline u =-\varepsilon\Delta \vphi +\frac1\varepsilon\Psi ^{\prime }(\vphi ),%
\end{cases}
\label{AGGlocal2d}
\end{equation}%
in $\Omega \times (0,+\infty)$, where   
\begin{equation}
\nu(\vphi )=\nu _{1}{\vphi}+\nu _{2}({1-\vphi }),\quad \rho (\vphi )=\rho _{1}\vphi+\rho _{2}({1-\vphi }),\quad
\JJJ=-({\rho _{1}-\rho _{2}})\nabla \overline u .
\label{RHO-Jtilde}
\end{equation}%
Here $\varphi:\Omega\times(0,\infty)\to [0,1]$ is the concentration of one of the fluids, where we observe that a partial mixing of both fluids is considered in the
model, and $\Psi$ is a suitable double-well potential, e.g., in the polynomial case, $\Psi (\varphi) = \vphi^2(1-\varphi)^2$. Furthermore, $\varepsilon>0$ is the small parameter related
to the interface thickness, and $\overline{u}$ is the chemical potential. The fluid mixture is driven by the capillary forces (Korteweg force) $-\varepsilon\mathrm{div}\,(\nabla \vphi \otimes \nabla \vphi )$, accounting for surface tension effects, which corresponds, as $\varepsilon\to 0$, to the effect of the mean curvature vector $\bH_{\pcA}$ in the interface condition \eqref{A8}. Additionally, differently from
the single fluid flow, the density $\rho (\vphi )$ in %
\eqref{AGGlocal2d} does not satisfy the continuity equation with respect to the flux associated with the velocity $\uuu $. Indeed, the density
satisfies the continuity equation with a flux given by the sum of the
transport term $\rho (\vphi )\uuu $ and the term $\JJJ$, due to the diffusion of the concentration in the unmatched densities case. To be precise, $\rho =\rho \left( \vphi \right)$ satisfies the continuity equation $$%
\partial _{t}\rho +\mathrm{div}\,\left( \rho \uuu +\JJJ\right) =0.$$
This is reflected in the corresponding sharp interface system \eqref{A1}-\eqref{A10}, where the continuity equation for $\rho(\xhi)$, simply obtained from the evolution equation of the characteristic function $\xhi(t)$ associated to $\cA(t)$, reads
\begin{align}
		\int_\Omega \pt\rho(\xhi)\eta\dx =(\rho_1-\rho_2)\int_{\pcA(t)\cap \om}\bv\cdot\bn_{\pcA(t)}\eta\d \mathcal H^{d-1}+\int_\om \widetilde{\bJ}\cdot\nabla \eta\dx,\label{densityA}
	\end{align}
	for any $\eta:\overline\om\to \R$ sufficiently smooth. 
    
    System \eqref{AGGlocal2d} is usually complemented with the boundary and
initial conditions
\begin{equation}
\begin{cases}
\uuu=0,\quad \partial _{\bn}\vphi =\partial _{\bn}\mu
=0,\quad  & \text{on }\partial \Omega \times (0,+\infty), \\
\uuu (0)=\uuu _{0},\quad \vphi (0)=\vphi _{0},\quad  & \text{%
in }{\Omega },%
\end{cases}
\label{locAGG-bc}
\end{equation}%
so that system \eqref{AGGlocal2d}, endowed with \eqref{locAGG-bc}, corresponds to the sharp interface problem \eqref{A1}-\eqref{A10} with constant contact angle $\gamma= \tfrac\pi2$. System \eqref{AGGlocal2d} has attracted the attention of many researchers, and its mathematical analysis can be found, for instance, in \cite{ADG, AGGio, ADGdegenerate}. Other variants of this model have been proposed in recent years, ranging from nonlocal versions \cite{Frigeri, GGGP} to multi-component variants \cite{AGGP}, to models in the framework of evolving surfaces \cite{AGGPsurf1,AGGPsurf2}. The sharp interface limit rigorously connecting suitable notions of weak solutions to system \eqref{AGGlocal2d} with \eqref{A1}-\eqref{A10} has been first investigated, in the case of matched densities $\rho_1=\rho_2$, in \cite{AbelsRoger}, and later extended to the general case in \cite{AbelsLengeler}. These works adapted the seminal techniques proposed for the sharp interface limit of the Cahn-Hilliard equation in \cite{Chen} to the case of couplings with hydrodynamic models. 

On the other hand, concerning the existence results of solutions to \eqref{A1}-\eqref{A10}, in the case when the mobility is zero, i.e., when the term related to $\nabla \overline u$ is not present in \eqref{A2} and the interface is transported by the flow of the
fluids, existence of strong solutions locally in time was first proven in \cite{DenisovaSolonnikov}, and we refer to \cite{PruessSimonett} for further references. Existence of generalized solutions globally in time was instead shown in \cite{Plotnikov, AbelsA,Abels2}. Moreover, the validity of a weak-strong uniqueness result has been recently shown in \cite{FischerHensel}. For the case with positive constant mobility, i.e., the Navier-Stokes-Mullins-Sekerka system \eqref{A1}-\eqref{A10}, the existence of global in time weak solutions and general initial data was first proven in \cite{AbelsRoger}, again in the simpler case of same densities $\rho_1=\rho_2$, whereas the existence of strong solutions locally in time, together with the stability of spherical droplets, was demonstrated in \cite{AbelsWilke}. 

In this contribution we primarily aim at extending the existence of weak solutions to \eqref{A1}-\eqref{A10} when the two fluid densities are different, i.e., $\rho_1\not=\rho_2$, which is a highly nontrivial task due to the presence of the additional flux $\wtilde{\bJ}$. Concerning this specific case, the only available result on the existence of weak solutions is the aforementioned \cite{AbelsLengeler}, where the authors show the existence of a varifold solution by means of a sharp interface limit as $\varepsilon\to0$ for the phase-field approximation \eqref{AGGlocal2d}. Here we propose a stronger notion of weak solution in the spirit of the recent \cite{SH}, i.e., including a weak formulation for the constant contact angle condition and a sharp energy dissipation
principle. Namely, we can allow the boundary condition for the interface not only for a constant contact angle $\gamma=\tfrac\pi2$, but also for general $\gamma\in (0,\tfrac\pi2]$. For the formulation of the energy dissipation inequality, we exploit a
gradient flow-like point of view, encoded in terms of a De Giorgi-type inequality. This structure can be seen from the formal energy identity associated to the system. To highlight this, let us introduce the following energy functional
	\begin{align}
		E[\xhi]:=\int_{\om}c_0 \d \normmm{\nabla\xhi}+\int_{ \pom}c_0\xhi\cos\gamma\d \mathcal H^{d-1},
		\label{E1}
	\end{align}
	so that, in the smooth case, it holds, by means of \eqref{A5},
		\begin{align*}
	E[\xhi]&=\int_{\pcA(t)\cap \om}c_0\d \mathcal H^{d-1}+\int_{\pcA(t)\cap \pom}c_0\cos\gamma\d \mathcal H^{d-1}\\&=c_0\cH^{d-1}(\pcA(t)\cap \Omega)+c_0 \cos\gamma\cH^{d-1}(\pcA(t)\cap \pom).
	\end{align*} 
    As it will be shown in greater detail in Appendix \ref{App:A}, the energy identity related to the Mullins-Sekerka part of the problem, i.e., to \eqref{A1}-\eqref{A5}, reads
\begin{align}
		&\ddt E[\xhi(t)]
		=-\int_{\pcA(t)\cap \om }(\bv(t)\cdot\bn_{\pcA(t)})\ovu(\cdot,t)\d \mathcal H^{d-1}-\int_{\om}\normmm{\nabla\ovu(\cdot,t)}^2\dx,\label{enid2A0}
	\end{align}
	whereas from the velocity $\bv$ equation we infer 
	\begin{align}
		&\frac12\ddt \int_\om \rho(\xhi(t))\normmm{\bv(t)}^2\dx +\int_\om \nu(\xhi(t))\normmm{D\bv(t)}^2\dx=\int_{\pcA(t)\cap \om }(\bv(t)\cdot\bn_{\pcA(t)})\ovu(\cdot,t)\d \mathcal H^{d-1}.\label{enid1A0}
	\end{align}
Summing the two, we finally deduce the total energy identity
	\begin{align}
		& \ddt \cE(t)=-\int_\om \nu(\xhi(t))\normmm{D\bv(t)}^2\dx-\int_{\om}\normmm{\nabla\ovu(\cdot,t)}^2\dx,\label{enid1A01}
	\end{align}
    where we have introduced the total energy of the system:
	\begin{align}
		\mathcal{E}(t):=E[\xhi(t)]+\frac12 \int_\om \rho(\xhi(t))\normmm{\bv(t)}^2\dx.
		\label{total1}
	\end{align}
Observe also that, thanks to the divergence-free condition of $\bv$, the total mass is conserved in the flow, i.e.,
	\begin{align}
	&	\ddt \int_\Omega \xhi(t)\dx=\ddt \mathcal L^d(\cA(t))=0.\label{massflow}
	\end{align}
Therefore, in view of \eqref{enid2A0}, the Mullins-Sekerka part of the flow can be equivalently represented with a gradient flow-like structure with respect to interfacial surface
energy, up to the presence of the advective velocity $\bv$, if we think of the time derivative of $\xhi$ as to be the \textit{material} time derivative, advected by $\bv$. Indeed, this can be guessed if,  recalling \eqref{massflow}, we rewrite \eqref{enid2A0} as 
\begin{align}
		&\ddt F(t)
		=-\norm{\partial^\bullet_t \xhi(t)}_{\Hmz}^2,\label{enid2A01}
	\end{align}
    where 
    \begin{align}
    F(t):=E[\xhi(t)]+\int_0^t\int_{\pcA(s)\cap \om }(\bv(s)\cdot\bn_{\pcA(s)})\ovu(\cdot,s)\d \mathcal H^{d-1}\d s,\label{discretization}
    \end{align}
    and where $\partial^\bullet_t \xhi:=\partial_t\xhi+\bv\cdot \nabla \xhi$ in $\Hmz$ represents the material derivative of $\chi$. This shows that the quantity $F(t)$ in \textit{nonincreasing} in time, in analogy with the gradient-flow structure.
In the case of $\bv=\mathbf 0$, say for the classical smooth Mullins-Sekerka flow, it is well known (see \cite{Garcke} and
references therein) that from \eqref{enid2A01} the gradient-flow structure can be realized in terms of a suitable $\Hmz$-type metric on a manifold of smooth surfaces. Additionally, inspired by De Giorgi’s methods for
curves of maximal slope in metric spaces and the approach for Gamma-convergence
of evolutionary equations developed in \cite{Serfaty, Le}, in \cite{SH} the authors are able to interpret even a weak solution to the Mullins-Sekerka flow as a Hilbert space gradient flow with respect to the interfacial energy $E$. 

Here we propose the nontrivial extension of these arguments to deal with a Mullins-Sekerka flow advected by $\bv$, i.e., a divergence-free vector field, tangential at the boundary $\partial \om$, satisfying the Navier-Stokes equations in $\Omega^+(t)\cup\Omega^-(t)$. Although it is well known that the Navier-Stokes equations (and their variant with unmatched densities and viscosities) do not enjoy a gradient flow structure, as we see from \eqref{enid2A01} the advective Mullins-Sekerka flow still conserves many properties in common with an $\Hmz$ gradient flow. We then exploit these properties to give a notion of solution which closely resembles the De Giorgi-type varifold solution of \cite{SH}, also taking into account the coupling with the Navier-Stokes equations. 

Our main result is the proof that such a new weak solution exists globally in time, showing for the first time that the minimizing movement scheme with De Giorgi interpolants is still viable to approximate the advective Mullins-Sekerka equation, as long as we resort to the notion of material derivative. Also, we show that smooth solutions are consistent with the classical Navier-Stokes-Mullins-Sekerka system. We point out that our new solution concept, in accordance to the one for the single Mullins-Sekerka flow of \cite{SH}, presents many novelties with respect to the previous notions of solutions to the Navier-Stokes-Mullins-Sekerka system. First, this notion allows to
include both a sharp energy dissipation principle and a weak formulation of
the contact angle at the intersection of the interface and $\pom$, which were left open by the previous notion proposed in \cite{AbelsRoger}, even in the case of same densities $\rho_1=\rho_2$. This is possible thanks to a single sharp energy
dissipation inequality à la De Giorgi, together with a weak formulation for any fixed contact angle by means of a distributional representation of the first variation of the corresponding capillary energy. Secondly, regarding the more general case with different densities, this notion is considerably stronger than the one so far available, which is a varifold solution obtained as a sharp interface limit in \cite{AbelsLengeler}. Namely, our new notion allows to further show the existence of an integrable generalized mean curvature vector, heavily exploiting the results of \cite{RogerA}. 

In conclusion, the advantage of this formulation is that we expect to be able, in a work in progress, to exploit the recently developed relative entropy approach for curvature
driven interface evolution to show a weak-strong uniqueness result, which is lacking with the previous notions of solution. This would show that weak solutions with sharp energy dissipation rate
are unique in a class of sufficiently regular strong solutions, as long as a strong solution exists, say until it develops, e.g., a singularity in the interface. Results in this direction have been obtained in \cite{FischerHensel} for the case when the mobility is zero, thus when the quantity related to $\nabla \overline u$ in \eqref{A2} is not present. Additionally, in the recent \cite{WSuniqueness} the authors have proven that the solutions constructed in \cite{SH} for the sole Mullins-Sekerka flow enjoy weak-strong uniqueness. Since our new notion of solution shares many properties with the one in \cite{SH}, and it enjoys a sharp energy dissipation inequality as well, we expect that a weak-strong uniqueness principle should be possible also for the complete Navier-Stokes-Mullins-Sekerka system \eqref{A1}-\eqref{A10}.

We now give some details of our existence proof. Following the gradient flow-type structure of the advective Mullins-Sekerka flow, we construct a new minimizing movement scheme with De Giorgi interpolants, i.e., a time discretization approach, for the Mullins-Sekerka flow advected by means of a suitably regularized velocity. At each time step we then update the advective velocity by solving a regularized version of the bulk Navier-Stokes equations with different densities and viscosities. The discretization of the advective Mullins-Sekerka flow takes some inspiration from the fundamental work \cite{Pironneau} and, more recently, from the new variational approach to the Navier-Stokes equations of \cite{Gigli}, in the sense that also here we treat the advective term in a semi-implicit fashion, i.e., the regularized advecting velocity $\bv$ is treated explicitely, whereas the advected quantity $\xhi$ is still implicit. This is achieved by means of a discretization of the \textit{material} derivative,  after introducing the composition between the advected variable $\xhi$ and a suitable flow map. Notice that this argument is considerably different from the one adopted in \cite{AbelsRoger}, where a fully explicit treatment of the advective term is performed. We make use of the well-known De Giorgi interpolation (see, e.g., \cite{Gradientflow}), crucially adapted to the presence of the flow map, in order to retrieve the optimal energy inequality. Due to the gradient flow-like structure of the advective Mullins-Sekerka flow, differently from \cite{Gigli}, here we are able to obtain the energy inequality for the phase indicator, corresponding to the discretized version of \eqref{enid2A01}, directly from the minimizing properties of the discretized solution, analogously to the standard gradient flow approach  \cite{Gradientflow}, without the necessity of resorting to any testing procedure on the Euler-Lagrange equations. We then pass to the limit in the time discretization parameter, which is the only approximating parameter in the scheme. Here we resort to some fine properties of flow maps, together with a delicate approximation of the velocity vector field $\bv$. Note that the use of a single parameter $h>0$ to approximate the problem is necessary, since we aim at applying the well-known results from \cite{Schatzle}, and these results can be apparently applied if the time-discretized solution directly comes from a minimizing movement scheme. To achieve this, we introduce two different approximating scales in the problem, i.e., the time discretization scale $h$ and the regularizing advective velocity scale $k$. By carefully choosing the second scale as a function of $h$, we are able to pass to the limit \textit{at once} as $h\to0$. One might wonder if we can solve first the decoupled advective Mullins-Sekerka flow given a suitable velocity, and then in a second step solve the velocity equation plugging in the result from the advective Mullins-Sekerka flow, exploiting a kind of decoupled approach and constructing a fixed point argument. This does not appear feasible, since the solution to the advective Mullins-Sekerka flow is in general nonunique, and thus we cannot expect the contraction map to be single-valued. As a consequence, at each time step of the minimizing movement scheme we also need to update the velocity $\bv$ solving the regularized Navier-Stokes equations, so that we avoid resorting to a contraction mapping technique.\\

\textbf{Structure of the paper. }The structure of the paper is as follows. The basic notation, some necessary preliminaries, and the main assumptions are summarized in Section \ref{Sec::preliminaries}. Then Section \ref{Sec::mainresults} is devoted to the definition of our new notion of weak solution, together with the statement of the main results. Section \ref{Sec::technical} deals with the proof of some essential technical lemmas, which are needed in the proof of our main existence theorem, which is the subject of Section \ref{proofthm1}. Section \ref{proofthm0} is then devoted to the proof of the consistency of any smooth weak solution with the classical Navier-Stokes-Mullins-Sekerka system. Finally, in Appendix \ref{App:A} we give the justification of some formulas used in the introduction, in Appendix \ref{measurability} we deal with some measurability issues appearing in the mininimizing movement scheme for the proof of our main result, whereas in Appendix \ref{App:C} we conclude with the proof of an embedding lemma in the case of weakly*-measurable vector-valued functions. 
\section{Functional setting and assumptions}
    \label{Sec::preliminaries}
	 Let $\Omega\subset \R^d$, $d=2,3$, be a bounded domain with smooth boundary. We denote by $\mathcal L^d$ the $d$-dimensional Lebesgue measure, whereas $\mathcal H^m$ represents the $m$-dimensional Hausdorff measure. We then denote by $C^k(\om)$ ($C^k(\overline{\om})$, respectively), $k\in\N\cup\{\infty\}$, the space of continuous functions on $\om$ ($\overline{\om}$, respectively) which are $k$-times continuously differentiable, as well as by $C^k_c(\om)$ ($C^k_c(\overline{\om})$,  respectively) the space of continuous functions with compact support in $\om$ ($\overline{\om}$, respectively) which are $k$-times continuously differentiable. The Sobolev spaces are denoted as usual by $W^{k,p}(\Omega )$%
	, where $k\in \mathbb{N}$ and $1\leq p\leq \infty $, with norm $\Vert \cdot
	\Vert _{W^{k,p}(\Omega )}$. The Hilbert space $W^{k,2}(\Omega )$ is classically indicated by $H^{k}(\Omega )$ with norm $\Vert \cdot \Vert _{H^{k}(\Omega )}$. Also, we denote by $BV(\om;\{0,1\})$ ($BV(\pom;\{0,1\})$, respectively) the space of BV functions over $\om$ ($\pom$, respectively) with discrete values in $\{0,1\}$. Recall that, given $f\in BV(\om;\{0,1\})$, we have
	$$
	\norm{\nabla f}_{\mathcal M(\om)}:=\sup_{\varphi\in C_c(\om),\ \norm{\varphi}_{C(\om)}\leq 1}\int_\om f \Div \varphi \dx.
	$$
	Moreover, given a (real) vector space $X$, we denote by $\mathbf{X}$ the generic space of vectors or matrices, with each component belonging to $X$. In this case $\vert \mathbf{v} \vert$ is the Euclidean norm
	of $\mathbf{v}\in  \mathbf{X}$, i.e., $\vert\mathbf{v}\vert^2=\sum_{j}\|v_j\|_X^2$.
	We then denote by $(\cdot,\cdot )_\om$ the inner product in ${L}^{2}(\Omega )$ (or $\bL^2(\om)$) and by $\Vert \cdot \Vert $
	the corresponding induced norm. Also, we denote by $\langle\cdot,\cdot \rangle_{X',X}$ the duality paring between $X'$ and $X$, given $X$ a generic Banach space. Then we set 
	\begin{align*}
		\Hz:=\left\{f\in H^1(\om):\ \int_\om f\dx=0\right\},\quad \Hmz:=(\Hz)'.\end{align*}
It is well-known that the homogeneous Neumann Laplace operator
\begin{equation*}
	-\Delta_N:\Hz\rightarrow H_{(0)}^{-1}(\Omega ) ,
	\quad
	\left\langle -\Delta_N  u,v\right\rangle_{	\Hmz,\Hz} =
	(\Grad u,\Grad v)
	\quad\text{for all $v\in H_{(0)}^{1}(\Omega )$}
\end{equation*}
is a continuous linear isomorphism. For any $f,g\in H_{(0)}^{-1}(\Omega )$, we set
\begin{equation}
(f,g)_{\Hmz}=\int_\om \nabla (-\Delta_N)^{-1}f\cdot \nabla (-\Delta_N)^{-1}g\dx.\label{Hmz1}
\end{equation}
This defines a bilinear form which is an inner product on the Hilbert space $H_{(0)}^{-1}(\Omega )$. Its induced norm $\Vert \cdot\Vert_{\Hmz}$ is equivalent to the standard operator norm on this space. 
Additionally, given $m_0\in(0,\mathcal L^d(\om))$, we let
\begin{align*}
	\mathcal M_{m_0}:=\left\{\xhi\in BV(\Omega;\{0,1\}):\ \int_\Omega \xhi\dx=m_0\right\}.
\end{align*}
We then introduce a class of regular test functions associated to $\xhi\in \mathcal M_{m_0}$, giving rise to infinitesimally volume preserving inner variations, denoted by
\begin{align}
	\mathcal S_{\xhi}:=\left\{B\in C^1(\overline\om;\R^d):\ \int_\om\chi\Div B\dx=0,\quad B\cdot\bn_{\partial\om}\equiv 0\text{ along }\pom \right\}.
	\label{Schi}
\end{align}
Also, we define the space of regular mass preserving normal velocities
generated on the interface associated with $\xhi\in \mathcal M_{m_0}$:
\begin{align}
	\mathcal W_{\xhi}:=\left\{B\cdot \nabla\xhi\in \Hmz:\ B\in \mathcal S_\xhi\right\},
	\label{Wchi}
\end{align}
with the convention that $B\cdot \nabla\xhi\in \Hmz$ is defined as 
\begin{align*}
	\langle B\cdot\nabla \xhi,\varphi\rangle_{\Hmz,	\Hz}:=-\int_\om \xhi \Div (B \varphi)\dx,\quad 	\forall \varphi\in \Hz.
\end{align*}
We also define its closure in $\Hmz$:
\begin{align}
	\mathcal V_{\xhi}:=\overline{\mathcal W_\xhi}^{\Hmz}\subset \Hmz.
	\label{Vchi}
\end{align}
The space $\mathcal V_\chi$ carries a Hilbert space structure directly induced by the natural Hilbert space structure of $\Hmz$, with the inner product induced by \eqref{Hmz1}. In particular, we set
\begin{equation}
	(f,g)_{\mathcal V_\xhi}:=(f,g)_{\Hmz},\quad 	\forall f,g\in \mathcal V_\xhi\label{Hmz2},
\end{equation}
 as well as $\norm{f}_{\mathcal V_\xhi}:=\norm{f}_{\Hmz}$ for any $f\in \mathcal V_\chi$.
In conclusion, we set
\begin{align}
	\mathcal G_{\xhi}:=-(-\Delta_N)^{-1}\mathcal V_	\xhi.
	\label{Gchi}
\end{align}
As shown in \cite[Section 2.3]{SH}, it holds 
\begin{align}
	\mathcal G_\xhi\subset \{f\in \Hz:\ \Delta f=0\text{ in }\om\setminus\text{supp}\normmm{\nabla\xhi}\}.
	\label{pchi}
\end{align}
Furthermore for the velocity field we set
\begin{align*}
	 \mathbf{C}^\infty_{c,\sigma}(\Omega):=\{\bu\in \mathbf C^\infty_c(\om):\ \Div \bu\equiv 0\},
\end{align*}
as well as
\begin{align*}
&\mathbf{L}^2_\sigma(\om):=\overline{ \mathbf{C}^\infty_{c,\sigma}(\Omega)}^{\mathbf{L}^2(\Omega)},\quad \mathbf H^1_\sigma(\Omega):=\mathbf H^1(\Omega)\cap \bL^2_\sigma(\Omega),\\&  	\mathbf{H}^1_{0,\sigma}(\om):=\overline{ \mathbf{C}^\infty_{c,\sigma}(\Omega)}^{\mathbf{H}^1(\Omega)},\quad	\mathbf{H}^k_{0,\sigma}:=\overline{ \mathbf{C}^\infty_{c,\sigma}(\Omega)}^{\mathbf{H}^k(\Omega)},
\end{align*}
for any $k\in \N\setminus\{0,1\}$. 

We set $\mathbf{P}_\sigma$ as the Leray projector from $\mathbf{L}^2(\Omega)$ to $\bL_\sigma(\om)$ and we also introduce $\bA=-\bP_\sigma\Delta: \mathfrak{D}(\bA)=\bH^2(\om)\cap \Hus\subset \bL^2_\sigma(\om)\to \bL^2_\sigma(\om)$ as the standard Stokes operator.

We then recall that Korn's inequality yields
\begin{equation}
	\Vert \mathbf{u}\Vert \leq \sqrt{2}\Vert D\mathbf{u}%
	\Vert\leq \sqrt{2}\Vert \nabla \mathbf{u}\Vert
	\quad \text{ for all } \mathbf{u}\in \Hus.
	\label{korn}
\end{equation}
Hence, $\|\nabla\cdot\|$ is a norm on $\Hus$ that is equivalent to the standard norm.

 Let $X$ be a Banach space. We denote by $L^q (a,b; X )$, $0\leq a\leq b$, $q\in[1,\infty]$, the Bochner space of $X$-valued $q$-integrable (or essentially bounded functions).
Moreover, $f:[0,\infty)\to X$ belongs to $L^q_{loc} ([0, \infty); X )$ if and only if it is strongly measurable with values in $X$ and $f \in L^q (0, T ; X )$ for every $T > 0$. Moreover, given $I=[0,T]$, we set $C_w(I;X)$ to be the topological vector space of weakly continuous functions $f:I\to X$. Given a generic interval $J$, the function space 
$C_c^\infty (J ; X)$ denotes the vector space of all $C^\infty$-functions $f : J\to X$ with compact support in $J$. Furthermore, $W^{1,p} (0, T ; X)$, $1 \leq p < \infty$, is the space of functions $f$ such that $\ddt f\in L^p(0,T;X)$ and $f\in L^p(0,T;X)$, where $\ddt$ denotes the vector-valued distributional derivative of $f$. 

In conclusion, if $Y=X'$ is a dual space, then $L^\infty_{w*}(0,T;Y)$, $T>0$, denotes the space of all functions $\eta: (0,T)\to Y$ which are weakly* measurable and essentially bounded, i.e. the map
$$
t\mapsto \langle \eta(t),f\rangle_{Y,X}
$$ 
is measurable for any $f\in X$, and 
\begin{align}
\norm{\eta}_{L^\infty_{w*}(0,T;Y)}:=\esssup_{t\in(0,T)}\norm{\eta(t)}_Y.
	\label{norm1}
\end{align}
Recalling (see \cite{AmbrosioFuscoPallara}) that $BV(\om)$ is the dual of a separable Banach space $X$, then we deduce $L^\infty_{w*}(0,T;BV(\om))=(L^1(0,T;X))'$ and that uniformly bounded sets in $L^\infty_{w*}(0,T;BV(\om))$ are weakly* compact.

	\subsection{Measures and varifolds}
	For a locally compact and separable metric space $X$, we denote
	by $\mathcal M(X)$ the space of (finite) Radon measures on $X$. An (oriented) varifold $\mu$ on $\overline\om$ is a positive Radon measure $\mu \in \mathcal M(\overline\om\times\mathbb S^{d-1})$. The corresponding varifold mass measure is denoted by $\normmm{\mu}_{\mathbb S^{d-1}}\in \mathcal M(\overline{\om})$, with $\normmm{\mu}_{\mathbb S^{d-1}}(A):=\mu(A\times \mathbb S^{d-1})$, for any Borel-measurable $A\subset \overline\om$.  We say that the varifold $\mu$ is
$(d-1)$-rectifiable if $\normmm{\mu}_{\mathbb S^{d-1}}$ is $(d-1)$-rectifiable. We say that $\mu$ is $(d - 1)$-integer-rectifiable if, additionally, the $(d - 1)$-density $\theta^{d-1}(\normmm{\mu}_{\mathbb S^{d-1}},\cdot)$ of $\normmm{\mu}_{\mathbb S^{d-1}}$ is integer valued.  

Observe that a simple example of an oriented varifold is the measure $\normmm{\nabla\xhi}\otimes \delta_{\frac{\nabla\xhi}{\normmm{\nabla\xhi}}(x)}$, associated to a set of finite perimeter represented via the characteristic function $\xhi\in BV(\om;\{0,1\})$. 

In conclusion, we indicate as tangential variation a vector field $B \in C^1(\overline{\om};\R^d )$, such that $\bn_{\pom}\cdot B\equiv 0$ along $\pom$. Then the tangential first variation of a varifold $\mu$ is defined as the first inner variation of the varifold computed with respect to a tangential variation and is here denoted by $\delta\mu(B)$ (see \cite{Allard}), i.e.,
\begin{align}
	\delta \mu(B)=\int_{\overline{\Omega}\times\mathbb{S}^{d-1}}(Id-s\otimes s):\nabla B(x)\d\mu(x,s),
	\label{allard1}
\end{align}
for any $B \in C^1(\overline{\om};\R^d )$, such that $\bn_{\pom}\cdot B\equiv 0$ along $\pom$. 

\subsection{Main assumptions}
	Concerning the unmatched fluid densities and viscosities, we consider the following assumptions:
	\begin{enumerate}
		\item[(\textbf{M0})] \label{M0}
		The viscosity $\nu:\{0,1\}\to \R$ is such that
		\begin{align}
			\nu(s):=\begin{cases}
				\nu_1,\quad s=1,\\
				\nu_2,\quad s=0,
			\end{cases}
			\label{nu1b}
		\end{align}
		where $\nu_1,\nu_2>0$,  such that 
		\begin{align}
			0<\min\{\nu_1,\nu_2\}=:\nu_*\leq  \nu^*:=\max\{\nu_1,\nu_2\}.
			\label{nu1}
		\end{align}
	\item[(\textbf{M1})]\label{M1a} The density $\rho:\{0,1\}\to \R$ is such that 
		\begin{align}
	\rho(s):=\begin{cases}
		\rho_1,\quad s=1,\\
		\rho_2,\quad s=0,
	\end{cases}
	\label{rho0}
\end{align}
	where $\rho_1,\rho_2>0$, so that 
	\begin{align}
		0<\min\{\rho_1,\rho_2\}=:\rho_*\leq \rho^*:=\max\{\rho_1,\rho_2\}.
		\label{rhocontrol}
	\end{align}
	\end{enumerate}

	\section{Main results}
    \label{Sec::mainresults}
	First we recall here \cite[Definition 3]{SH}, which gives the admissibility conditions for a couple given by evolving phase indicators and varifolds.
	\begin{definition}[Admissible couples of evolving phase indicators and varifolds]
	Let $d=2,3$, and set $\Ts\in(0,\infty)$. Let $\xhi\in L^\infty_{w*}(0,\Ts;\mathcal M_{m_0})\cap C([0,\Ts);\Hmz)$. Consider a family $\mu=(\mu_t)_{t\in(0,\Ts)}$ of oriented varifolds $\mu_t\in \mathcal M(\overline\Omega\times \mathbb S^{d-1})$, for all $t\in(0,\Ts)$. The couple $(\chi,\mu)$ is said to be admissible if the following properties are satisfied: 
	\begin{enumerate}
		\item\textit{(Structure of oriented varifolds).} For almost every $t\in(0,\Ts)$ the oriented varifold $\mu_t$ decomposes as $\mu_t=c_0\mu_t^\Omega+(\cos\gamma)c_0\mu_t^{\partial\Omega}$ for two separate oriented varifolds given in their disintegrated form by 
		\begin{align}
			\label{pA}
			\mu_t^\Omega=: \normmm{\mu_t^\Omega}_{\mathbb S^{d-1}}\otimes (\lambda_{x,t})_{x\in \overline{\Omega}}\in \mathcal M(\overline{\Omega}\times \mathbb S^{d-1})
		\end{align} 
		and 
			\begin{align}
			\label{pB}
			\mu_t^{\partial\Omega}=: \normmm{\mu_t^{\partial\Omega}}_{\mathbb S^{d-1}}\otimes (\delta_{\bn_{\pom}(x)})_{x\in {\partial\Omega}}\in \mathcal M(\overline{\Omega}\times \mathbb S^{d-1}).
		\end{align}
		For almost any $t\in(0,\Ts)$ the measure $\normmm{\mu_t^\Omega}_{\mathbb S^{d-1}}$ is $(d-1)$-integer rectifiable, and the measure $\normmm{\mu_t^{\partial\Omega}}_{\mathbb S^{d-1}}$ is given by $g_t\mathcal H^{d-1}\llcorner \partial \Omega$, where $g_t\in BV(\partial\Omega;\{0,1\})$.
		\item\textit{(Compatibility with the phase indicator).} The oriented varifolds $\mu^\Omega_t$ and $\mu_t^{\partial\Omega}$ contain the interfaces associated   with the phases modeled by $\chi$, in the sense that 
		\begin{align}
			&c_0\normmm{\nabla\xhi}\nonumber
			\\& 
			=c_0\sum_{k=1}^\infty \frac1{2k-1}\normmm{\mu_t^\Omega}_{\mathbb S^{d-1}}\llcorner \left(\Omega \cap \{\theta^{d-1}(\normmm{\mu_t^\Omega}_{\mathbb S^{d-1}}\llcorner \Omega,\cdot)=2k-1\}\right),
			\label{PC}
		\end{align}
		and 
		\begin{align}
			(\cos\gamma)\xhi(\cdot,t)\mathcal H^{d-1}\llcorner \partial\Omega\leq \cos\gamma\left(\normmm{\mu_t^\Omega}_{\mathbb S^{d-1}}\llcorner\partial\Omega+\normmm{\mu_t^{\partial\Omega}}_{\mathbb S^{d-1}}\right),
			\label{PD}
		\end{align}
		for almost every $t\in(0,\Ts)$.  Additionally, for almost every $t\in(0,\Ts)$, for every $\eta\in C^{1}(\overline\Omega;\mathbb R^d)$ such that $\mathbf n_{\partial\Omega} \cdot \eta=0$ on $\partial\Omega$ and every $\xi\in C^1(\overline\Omega;\mathbb R^d)$ with $\bn_{\partial\Omega}\cdot \xi=\cos\gamma$ on $\partial\Omega$, it holds that
		\begin{align}
			&\int_{\overline\Omega\times \mathbb{S}^{d-1}}s\cdot \eta(x)\d\mu_t^\Omega(x,s)=\int_\Omega \frac{\nabla \xhi(\cdot,t)}{\normmm{\nabla\xhi(\cdot,t)}}\cdot \eta(\cdot)\d\normmm{\nabla\chi(\cdot,t)},\label{PE}\\&
					-\int_{\overline\Omega\times \mathbb{S}^{d-1}}s\cdot \xi(x)\d\mu_t^\Omega(x,s)\nonumber\\&=-\int_\Omega \frac{\nabla \xhi(\cdot,t)}{\normmm{\nabla\xhi(\cdot,t)}}\cdot \xi(\cdot)\d\normmm{\nabla\chi(\cdot,t)}
					+\cos\gamma\left(\tmup(\pom)-\int_\Omega \xhi(\cdot,t)\d\mathcal H^{d-1}\right)
					.\label{PF}
		\end{align} 
		
	\item\textit{(Existence of a generalized mean curvature vector).} For almost every $t\in(0,\Ts)$ there exists a map $\mathbf H_{\tmu\llcorner\om}:\text{supp}\left(\tmu\llcorner\om\right)\to \mathbb R^d$ such that 
	\begin{align}
		\mathbf H_{\tmu\llcorner\om}\in L^s(\Omega; \d\tmu),
		\label{PG}
	\end{align} 
	for $s\in[1,4]$ if $d=3$, $s\in[1,\infty)$ if $d=2$, and such that the first variation $\delta\mu_t$ of $\mu_t$ in the direction of a tangential field $B\in C^1(\ov\om; \mathbb R^d)$, $\bn_\pom\cdot B\equiv 0$ along $\pom$ is given by 
	\begin{align}
		\delta\mu_t(B)=-\int_\om c_0\mathbf H_{\tmu\llcorner\om}\cdot B\d\tmu\llcorner \Omega.
		\label{PH}
	\end{align}
Also, $\mu_t$ is of bounded first variation on $\overline\om$ and such that
\begin{align}
	\sup_{B\in C^1(\overline\om;\R^d), \norm{B}_{L^\infty(\Omega)}\leq 1}\normmm{\delta\mu_t(B)}\leq C(\om)\tmu(\pom)+\wtilde C \normmm{\bH_{\tmu\llcorner \Omega}}_{L^1(\Omega,\d \normmm{\mu}_{\mathbb{S}^{d-1}})},
	\label{PI}
\end{align}	
for some universal constant $\wtilde C>0$ and some $C(\om)>0$, this last one depending only on the second fundamental form of $\pom$.
\item \textit{(Structure of the generalized mean curvature).} The generalized mean curvature vector 	
$\bH_{\tmu\llcorner \Omega}$ satisfies 
\begin{align}
	&\bH_{\tmu\llcorner \Omega}=0\quad\text{ in }\Omega\cap \{\theta^{d-1}(\tmu\llcorner\om,\cdot)>1\},
	\label{PJ}\\&
	\label{PK}
	\bH_{\tmu\llcorner \Omega}=H_{\xhi(\cdot,t)}\frac{\nabla \xhi(\cdot,t)}{\normmm{\nabla\xhi(\cdot,t)}},\quad\text{ in }\text{supp}(\normmm{\nabla\xhi(\cdot,t)}),
\end{align}
$H_{\xhi(\cdot,t)}:\text{supp}(\normmm{\nabla\xhi(\cdot,t)})\to \R$ stands for the generalized mean curvature of $\text{supp}(\normmm{\nabla\xhi(\cdot,t)})$ in the sense of \cite[Definition 1.2]{Roger}.
\item\textit{(Measurability in time of the energy).} The total mass measure associated with the oriented varifold $\mu_t$, that is
\begin{align}
	E[\mu_t]:=\normmm{\mu_t}_{\mathbb S^{d-1}}(\overline{\om})=c_0\tmu(\overline\om)+c_0\cos\gamma\tmup(\pom),
	\label{PL}	
\end{align}
is a measurable map from $t\in(0,\Ts)$ to $E[\mu_t]\in[0,\infty)$.
\end{enumerate}
\label{admissibility}
\end{definition}
\begin{remark}
		As it has also been observed in \cite[Remark 4]{SH}, in general we cannot expect the contribution of $\normmm{\mu_t^\Omega}$
		to be zero on $\partial\Omega$, not even in the case $\gamma=\tfrac\pi2$. Indeed,	these measures will be constructed from a local minimization problem, but in the limit
		this will allow for surfaces laying tangentially on $\partial\Omega$. This means that in general the inequality (in the sense of measures) $\normmm{\mu_t^\om}_{\mathbb S^{d-1}}\llcorner \Omega\leq  \normmm{\mu_t^\Omega}_{\mathbb S^{d-1}}$ might be strict.
	\end{remark}

	\begin{remark}
	Observe that, by setting $\cA(t):=\{x\in \Omega:\ \xhi(x,t)=1\}$, for almost any $t\in (0,\Ts)$, we easily obtain that $\xhi(\cdot,t)=\xhi_{\cA(t)}(\cdot)$ and the family $\{\cA(t)\}_{t\in(0,\Ts)}$ is a family of finite perimeter sets contained in $\Omega$, with (inner) unit normal to the reduced boundary $\partial^*\cA(t)$ given by $\bn_{\partial^*A(t)}:=\tfrac{\nabla \xhi(\cdot,t)}{\normmm{\nabla\chi(\cdot,t)}}$ on $\partial^*\cA(t)$. Also, it holds $\text{supp}({\normmm{\nabla\xhi(\cdot,t)}})=\partial^*\cA(t)$ and $\normmm{\nabla\xhi(\cdot,t)}=\mathcal H^{d-1}\llcorner \partial^*\cA(t)$ for almost any $t\in(0,\Ts)$.\label{A(t)}
	\end{remark}	
	We can now give our new notion of solution to the Navier-Stokes-Mullins-Sekerka system with variable densities and viscosities: 
	\begin{definition}
		[Varifold solutions to Navier-Stokes-Mullins-Sekerka system]
		Let $d=2,3$, consider $\Ts\in(0,\infty)$, $\rho,\nu$ as in assumptions (\textbf{M0})-(\textbf{M1}), and let $\om\in\R^d$ be bounded with smooth boundary $\pom$. Fix also $\chi_0\in \mathcal M_{m_0}$, $m_0\in(0,\cL^d(\om))$, and $\bv_0\in \bL^2_\sigma(\Omega)$. 
		
		A measurable map $\xhi:\Omega\times(0,\Ts)\to \{0,1\}$, together with a family $\mu=(\mu)_{t\in(0,\Ts)}$ of oriented varifolds $\mu_t\in \mathcal(\ov\om\times\mathbb{S}^{d-1})$, $t\in(0,\Ts)$, and a vector field $\bv:\Omega\times(0,\Ts)\to \R^d$ is called a \textit{varifold solution }for the Navier-Stokes-Mullins-Sekerka system with variable densities and viscosities, with time horizon $\Ts$ and initial data $(\xhi_0,\bv_0)$, if 
		\begin{enumerate}
			\item\textit{(Admissibility and regularity).} $(\xhi,\mu)$ is an admissible couple in the sense of Definition \ref{admissibility}, and additionally 
			\begin{align}
			\xhi\in	H^1(0,\Ts;\Hmz)\cap C([0,T^*];L^p(\om));\label{regularitychi}
			\end{align}
		for any $p\in[1,\infty)$.
			\item\textit{(Fluid velocity regularity).} The fluid velocity $\bv$ is such that 
				\begin{align}
				&\bv\in C_w([0,\Ts];\Ls)\cap L^2(0,\Ts;\Hus).
				\label{regv}
			\end{align} 
		\item\textit{(Consistency with initial data.)} It holds $\xhi(0)=\xhi_0$ in $\Hmz$ and $\bv(0)=\bv_0$ almost everywhere in $\Omega$.
		\item\textit{(Curvature potential for the mean curvature). }There exists a potential $w\in L^2(0,T;H^1(\Omega))$, with $w_0:=w-\tfrac{\int_\om w\dx}{\mathcal L^d(\om)}$ such that $w_0(t)\in \mathcal G_{\xhi(\cdot,t)}$ for almost any $t\in(0,\Ts)$. 
		Additionally, the potential $w$ satisfies the Gibbs-Thomson law in the sense that 
		\begin{align}
			\int_{\overline\Omega\times \mathbb S^{d-1}}(Id-s\otimes s):\nabla B(x)\d \mu_t(x,s)=\int_{\Omega}\xhi(x,t)\Div(w(x,t)B(x))\dx,
			\label{GThom}
		\end{align}
		for almost every $t\in(0,\Ts)$ and all $B\in C^1(\overline\Omega;\mathbb R^d)$ such that $B\cdot\bn_{\partial\Omega}\equiv0$ on $\partial\Omega$;
		\item\textit{(Kinetic potential).} There exists a potential $u\in L^2(0,\Ts;H^1(\Omega))$ with ${\int_\Omega u(\cdot,t)\dx}=0$, for almost any $t\in(0,\Ts)$, such that 
		\begin{align}
		&\int_{\om}\xhi(\cdot,\Ts)\zeta(\cdot, \Ts)\dx-\int_{\om}\xhi_0\zeta(\cdot,0)\dx\nonumber\\&
        =\int_0^{\Ts}\int_\Omega \xhi(\cdot,t)\pt \zeta(\cdot,t)\dxdt+\int_0^{\Ts}\int_{\om}\xhi\bv\cdot \nabla \zeta\dx\dt-\int_0^{\Ts}\int_\Omega \nabla u\cdot\nabla \zeta\dx\dt,
			\label{kinetic}
		\end{align} 
		for all $\zeta\in C^\infty_c(\overline\Omega\times[0,\Ts])$; 
		\item \textit{(Fluid velocity equation).} It holds
		\begin{align}
			&	\int_0^{\Ts}\int_{\om}\left(-\rho( \xhi)\bv\cdot \pt\bpsi-\rho(\xhi)(\bv\otimes\bv):\nabla  \bpsi+\nu(\xhi)D\bv:D\bpsi\right)\dx\dt\nonumber\\&
			-\int_0^{\Ts}\int_{\om}\left((\bv\otimes \widetilde{\mathbf J}):\nabla\bpsi\right)\dx\dt\nonumber\\&=\int_\om \rho(\xhi_0)\bv_0\cdot\bpsi(0)\dx-\int_{0}^{\Ts}\int_{\om}\xhi\nabla w\cdot \bpsi\dx\dt,
			\label{velocity}
		\end{align}
		for any $\bpsi\in  C^\infty_c([0,\Ts);\Hds)$, with
		\begin{align}
			\widetilde{\mathbf J}:=-(\rho_1-\rho_2)\nabla u.
			\label{tildeJ}
		\end{align}
		\item \textit{(Sharp energy dissipation inequalities).} The following energy inequalities hold: 
		\begin{align}
			\nonumber &E[\mu_{T'}]+\int_s^{T'}\int_\Omega \frac12 \normmm{\nabla u}^2\dx \dt +\int_s^{T'}\int_\om \onehalf\normmm{\nabla w}^2\dx\dt\\&\leq E[\mu_s]+\int_s^{T'}\int_{\Omega}\xhi\bv\cdot \nabla w\dx\dt,
			\label{advectiveMullins}	
		\end{align}
		 for almost any $T'\in(0,\Ts)$ and almost any $s\in[0,T')$, including $s=0$, where
		 $$
		 E[\mu_t]:=\begin{cases}
		 	c_0\tmu(\overline\om)+c_0\cos\gamma\tmup(\pom),\quad t>0,\\
		 	E[\xhi_0],\quad t=0,
		 	\end{cases}
		 $$
		 with $E[\xhi_0]$ defined in \eqref{E1}, as well as 
		 \begin{align}
		 &\nonumber	\onehalf \int_\Omega \rho(\xhi(T'))\normmm{\bv(T')}^2\dx+\int_s^{T'}\int_\Omega\nu(\xhi)\normmm{D\bv}^2\dx\dt \\&\leq \onehalf \int_\Omega \rho(\xhi(s))\normmm{\bv(s)}^2\dx-\int_s^{T'}\int_\Omega \xhi\bv\cdot\nabla w\dx\dt,
		 	\label{NS}
		 \end{align}  
		  for any $T'\in(0,\Ts)$ and almost any $s\in[0,T')$, including $s=0$.
		\end{enumerate}	 
		\label{weaksol}   
	\end{definition}
	\begin{remark}[Total energy dissipation inequality]
	It is trivial to notice that Definition \ref{weaksol}, due to \eqref{advectiveMullins}-\eqref{NS}, entails the validity of the standard (sharp) energy inequality:
	\begin{align}
		\mathcal{E}(T')+\int_s^{T'}\int_\Omega\nu(\xhi)\normmm{D\bv}^2\dx\dt+\int_s^{T'}\int_\Omega \frac12 \normmm{\nabla u}^2\dx \dt +\int_s^{T'}\int_\om \onehalf\normmm{\nabla w}^2\dx\dt
		\leq \mathcal E(s),
		\label{energytot}
	\end{align}
	 for almost any $T'\in(0,\Ts)$ and almost any $s\in[0,T')$, including $s=0$, where 
	 \begin{align*}
	 	\mathcal E(t):=\nonumber	\onehalf \int_\Omega \rho(\xhi(t))\normmm{\bv(t)}^2\dx+E[\mu_t],
	 \end{align*}
	 for almost any $t\in(0,\Ts)$, is the total energy of the system. On the other hand, we point out that the inequality is valid only \textit{for almost any }$0<s<T'$, due to the low regularity of $E[\mu_t]$.\label{toten1}
	\end{remark}
	\begin{remark}
		In view of \eqref{PC}, \eqref{PH}-\eqref{PK}, and \eqref{GThom}, if $(\chi,\mu,\bv)$  satisfies Definition \ref{weaksol}, we immediately deduce that
		\begin{align}
		&	c_0H_{\chi(\cdot,t)}=w(\cdot,t),\quad	\text{on }\text{supp}(\normmm{\nabla\xhi(\cdot,t)}),\label{wkk1}\\&
			w(\cdot,t)=0,\quad\text{ on }\Omega\cap \{\theta^{d-1}(\tmu,\cdot)\in 2\mathbb N+1\}.
		\end{align}
		\label{identif}
		Therefore, the surface tension term in the right-hand side of \eqref{velocity} can be rewritten as
		$$
		-\int_{0}^{\Ts}\int_{\om}\xhi\nabla w\cdot \bpsi\\dx\dt=\int_0^\Ts\int_\om c_0 H_{\xhi(\cdot,t)}{\frac{\nabla\xhi(\cdot,t)}{\normmm{\nabla \xhi(\cdot,t)}}}\cdot \bpsi\d \normmm{\nabla \xhi(\cdot,t)}.
		$$
	\end{remark}
\begin{remark}[Kinetic and curvature potentials]
		One of the main features of the notion of solution given in the definition above is the fact that there are two (possibly) different potentials at the same time, i.e., $u$ and $w$, with different roles. This distinction is reflected into the Navier-Stokes equations \eqref{velocity}, where the diffusive flux $\widetilde{\mathbf J}$ is related to the kinetic potential $u$, since it is defined in \eqref{tildeJ} by means of the evolution equation of the density $\rho(\xhi)$, coming from \eqref{kinetic}. On the other hand, the surface tension term in the right-hand side of \eqref{velocity} is related to the curvature potential $w$, since the surface tension depends on the mean curvature (see also Remark \ref{identif}).
	\end{remark}
In view of the definition of varifold solution given in Definition \ref{weaksol}, we can state our main global existence result, whose proof is shown in Section \ref{proofthm1}:
\begin{theorem}[Existence of global varifold solutions to the Navier-Stokes-Mullins-Sekerka system]
	Let $d \in \{2,3\}$, $\om\subset \R^d$ be a bounded domain with smooth boundary, and $\rho,\nu$ as in assumptions  (\textbf{M0})-(\textbf{M1}). Let $\bv_0\in \bL^2_\sigma(\om)$, $m_0 \in (0,\mathcal L^d(	\om))$, $\xhi_0 \in \mathcal M_{m_0}$, $c_0\in (0,\infty)$, $\gamma \in (0, \pi/2]$.
	Then, there exists a varifold solution for the Navier-Stokes-Mullins-Sekerka system with variable densities and viscosities,  globally defined on $(0,\infty)$, with
	initial data ($\chi_0$,$\bv_0$), in the sense of Definition \ref{weaksol} for any $\Ts>0$.
	Additionally, for any $\Ts>0$ there exists $C_p=C_p(d,\om,\Ts,\xhi_0,\bv_0)>0$ such that the curvature potential satisfies
	\begin{align}
		&\norm{w(\cdot,t)}_{H^1(\om)}\leq C_p(1+\norm{\nabla w(\cdot,t)}),\label{C1}
	\end{align}
	for almost every $t\in(0,\Ts)$. 
	\label{thm1}
\end{theorem}
In conclusion, we show that the notion of varifold solutions on $(0,\Ts)$ introduced in Definition \ref{weaksol} is consistent, extending \cite[Lemma 4]{SH}. The proof of this Theorem is proposed in Section \ref{proofthm0}.
\begin{theorem}[Consistency]
	\label{consistency}
	Let the assumptions and notation of Theorem \ref{thm1} hold.
	\begin{enumerate}
	\item Any classical (i.e., smooth) solution $(\mathcal A, \bv)$ of the Navier-Stokes-Mullins-Sekerka system \eqref{A1}-\eqref{A10} on $[0,\Ts)$, $\Ts>0$, gives
	rise to a (smooth) varifold solution $(\xhi,\mu,\bv)$, with $\xhi=\xhi_\cA$, in the sense of Definition \ref{weaksol}.
	\item On the other hand, let $(\xhi,\mu,\bv)$ be a varifold solution for the Navier-Stokes-Mullins-Sekerka system in the
	sense of Definition \ref{weaksol} which is smooth, i.e., $\xhi(x, t) = \xhi_{\cA(t)}(x)$
	for a smoothly evolving family $\cA = (\cA (t))_t\in[0,\Ts)$, $\Ts>0$. Also, $\mu$ and $\bv$ are smooth in space and time. Moreover, assume
	that \eqref{PH} holds with  $\delta E[\xhi(\cdot, t)]$ in place of $\delta\mu_t$. Then, $(\mathcal A,\bv)$ is a classical solution for the Navier-Stokes-Mullins-Sekerka system satisfying \eqref{A1}-\eqref{A10}.
	Furthermore, it holds 
	\begin{align}
		&\tmu\llcorner \om =\normmm{\nabla\xhi(\cdot,t)},\quad \tmu\llcorner \pom=0,
		\label{H1},\\&
		(\cos\gamma)\tmup=(\cos\gamma)\xhi(\cdot,t)\cH^{d-1}\llcorner \pom,
		\label{H2}
	\end{align}
	for almost any $t\in(0,\Ts)$.
	
	\end{enumerate}
\end{theorem}
\begin{remark}
	As observed in \cite[Remark 8]{SH}, there are cases in which the assumption that \eqref{PH} holds with  $\delta E[\xhi(\cdot, t)]$ in place of $\delta\mu_t$ is already satisfied. For instance, one could construct slightly more regular (i.e., BV) solutions in place of the varifold solution of Definition \ref{weaksol}, for which this condition is satisfied. Nevertheless we point out that the existence of such solutions is so far only conditional even for the Mullins-Sekerka flow (see \cite[Theorem 1]{SH}).
\end{remark}
		\section{Technical Lemmas}
        \label{Sec::technical}
Before starting the proofs of Theorems \ref{thm1} and \ref{consistency}, we aim at presenting some technical lemmas, which are essential in the sequel. 
\subsection{Preliminary results}
We recall here a proposition whose proof is given in \cite[Proposition 5]{SH} and \cite[Corollary 6, eq. (37)]{SH}. This will be of some use in the sequel.  
\begin{proposition}[First variation estimate up to the boundary for tangential variations]\label{prop5}
	Let $d=2,3$ and $\Omega\subset \R^d$ be a bounded domain with orientable $C^2$-boundary $\pom$. Let also $w\in H^1(\Omega)$, $	\xhi\in BV(\Omega;\{0,1\})$ and assume $\mu\in \mathcal M(\overline\Omega\times \mathbb S^{d-1})$ is an oriented varifold such that $c_0\normmm{\nabla \xhi}\leq \normmm{\mu}_{\mathbb S^{d-1}}\llcorner\om$ in the sense of measures for some constant $c_0>0$. Assume moreover that it holds
	\begin{align*}
		\int_{\overline\Omega\times \mathbb S^{d-1}}(Id-s\otimes s):\nabla B\d\mu=\int_\Omega \xhi\Div(w B)\dx,
	\end{align*}  
	for all $B\in C^1(\overline{\Omega},\R^d)$, with $\bn_{\pom}\cdot B\equiv 0$ along $\pom$. Then, for all $s\in[2,4]$ if $d=3$ or $s\in[2,\infty)$ otherwise, there exists a constant $C=C(\pom,s,d)$ such that $w$ is integrable with respect to $\normmm{\mu}_{\mathbb S^{d-1}}$ and satisfies 
	\begin{align*}
	\left(\int_{\overline{\om}}\normmm{w}^s \d\normmm{\mu}_{\mathbb S^{d-1}}\right)^{\frac 1s} \leq C(1+\normmm{\mu}_{\mathbb S^{d-1}}(\overline\Omega)+\norm{w}^d_{H^1(\om)})^{1+\frac1s}.
	\end{align*} 
	Additionally, the first variation  of $\mu$ with respect to tangential variations is given by a generalized mean curvature vector $\mathbf H^\om=\rho^\om \frac{w}{c_0}\frac{\nabla \xhi}{\normmm{\nabla \xhi}}$, where $\rho^\om:=\frac{c_0\normmm{\nabla \xhi}}{\normmm{\mu}_{\mathbb S^{d-1}}\llcorner \om}\in[0,1]$, in the sense of \eqref{PH}.
 Also, $\mu$ is of bounded variation on $\overline \om$ and its first variation $\delta \mu$ satisfies 
 \begin{align}
 	\sup_{B\in C^1(\overline\om),\norm{B}_{L^	\infty(\om)}\leq1}\normmm{\delta\mu(B)}\leq C\normmm{\mu}_{\mathbb S^{d-1}}(\overline\om)+\norm{\mathbf H^\om}_{L^1(\om,\d\normmm{\mu}_{\mathbb S^{d-1}})}.\label{totalvar}
 \end{align} 
\end{proposition}
We also recall here the following lemma, which is fundamental to pass from a Gibbs-Thomson relation holding
for infinitesimally mass-preserving velocities to arbitrary
variations tangential at the boundary. Namely, it gives a control on the Lagrange multiplier arising from the mass constraint. A similar lemma can be found in \cite{AbelsRoger}, whereas the proof of the statement below is shown in \cite[Proof of Lemma 9]{SH}.
\begin{lemma}
	\label{lemma9}
	Let $\xhi\in \mathcal M_{m_0}$, $w\in \Hz$, and $\mu\in \mathcal M(\overline\om\times\mathbb S^{d-1})$ be an oriented varifold such that 
	\begin{align*}
		\delta\mu(B)=\int_\Omega \xhi\Div(w B)\dx,\quad\forall B\in \mathcal S_{\xhi}.
	\end{align*}
	Then there is some $\lambda \in \R$ such that \begin{align*}
		\delta\mu(B)=\int_\Omega \xhi\Div((w+\lambda)B)\dx,\quad \forall B\in C^1(\overline\om;\R^d),\ \text{with }B\cdot\bn_\pom\equiv 0\text{ along }\pom,
	\end{align*}
	and there exists $C=C(\om,d,m_0)$ such that 
	\begin{align*}
		\norm{w+\lambda}_{H^1(\om)}\leq C(1+\normmm{\nabla\xhi}(\om))(\normmm{\mu}_{\mathbb S^{d-1}}(\overline \om)+\norm{\nabla w}).
	\end{align*}
\end{lemma}

\subsection{Flow map and related properties}
We consider the following flow map: given a $\mathbf C^1(\overline\Omega)$  divegence-free vector field $\bv$, find $X_t^\bv:\overline\om\to \overline\om$ satisfying, for any $t\in \mathbb R$,
	\begin{align}
		\begin{cases}
			\ddt X_t^{\bv}=\bv(X_t^\bv),\\
			X_0^\bv=Id.
		\end{cases}\label{flowmap}
	\end{align}
	This problem admits a unique solution $X^\bv\in C^1((-\infty,+\infty),\mathbf C^1(\overline\Omega))$. Also, the flow map $X_t^\bv$ is a $C^1(\overline{\Omega})$-diffeomorphism for any $t\in \mathbb R$. 
We now show some results concerning the composition between BV functions and flow maps. First, we have the following lemma, whose proof can be adapted from \cite[Lemma 10]{SH}. 
	\begin{lemma}\label{vC1}
		Let $\bv\in \mathbf C^1(\overline{\Omega})\cap \bH^1_\sigma(\om)$ and consider the corresponding flow map $X_t^\bv$ as in \eqref{flowmap}. Then, for any $\xhi\in \mathcal M_{m_0}$ it holds 
		\begin{align}
			\frac{\xhi-\xhi\circ X_{-s}^\bv}{s}\to \bv\cdot \nabla \xhi,\quad\text{strongly in }\Hmz,
			\label{convHm1}
		\end{align}
		as $s\to0^+$. 
		\label{basic1}
	\end{lemma}
	\begin{remark}
		The above lemma also entails, 
		\begin{align}
			\lim_{s\to 0}\frac{\normh{\xhi-\xhi\circ X_{-s}^\bv}^2}{s}=0,
			\label{convHm2}
		\end{align} 
		since, by \eqref{convHm1},
		\begin{align*}
			\lim_{s\to 0}\frac{\normh{\xhi-\xhi\circ X_{-s}^\bv}^2}{s^2}=\normh{\bv\cdot \nabla \chi}^2.
		\end{align*}
	\end{remark}
	\begin{proof}[Proof of Lemma \ref{vC1}]
		We aim at applying \cite[Lemma 10]{SH}, and thus we need to identify the family of diffeomorphisms $\Psi_s$, as well as the field $B$ with suitable properties. Now, in our case we have by assumption $B:=\bv\in S_{\xhi}$, where $S_\xhi$ is defined in \eqref{Schi}. Moreover, we can define $\Psi_s:=X_s^{\bv}:\overline{\Omega}\to \overline{\Omega}$, so that $\Psi_s^{-1}=X_{-s}^\bv$, for $s\in(-1,1)$. Notice that, from the properties of the flow map, we immediately infer that $\{\Psi_{s}\}_{s\in(-1,1)}$ is a family of $C^1$-diffeomorphism which depends differentiably on $s\in(-1,1)$. Also, again by \eqref{flowmap} it holds $\Psi_0(0)=X_0^\bv(x)=x$ for any $x\in\overline{\Omega}$, as well as $\partial_s\Psi_s|_{s=0}=\partial_s X_s^\bv|_{s=0}=\bv(X_s^\bv)|_{s=0}=\bv(x)=B(x)$. Additionally, we have  
		\begin{align*}
			\int_\Omega \xhi\circ \Psi^{-1}_s\dx=\int_\Omega \xhi\circ X_{-s}^\bv\dx=\int_\Omega \xhi\dx=m_0,
		\end{align*}
		since $\xhi\in \mathcal M_{m_0}$, and ${X_s^\bv}_\sharp\mathcal L^d=\mathcal L^d$ for any $s\in \mathbb R$. Therefore, with this choice of $\Psi_s$ and $B$, all the assumptions of \cite[Lemma 10]{SH} are satisfied and we thus infer \eqref{convHm1}, concluding the proof.
	\end{proof}
Additionally, we also have the estimates stated in the following lemma, which is a subtle refinement of the results in \cite{Gigli}.
\begin{lemma}
   	Let $\bv\in \mathbf C^1(\overline{\Omega})\cap \Hs$. Then it holds 
		\begin{align}
			\normh{\xhi-\xhi\circ X_{-s}^\bv}\leq s\norm{\bv}_{\mathbf C(\overline{\Omega})},\quad \forall \xhi\in \mathcal M_{m_0},\quad \forall s\geq 0,
			\label{H1c}
		\end{align}
		as well as
		\begin{align}
			&\nonumber\norm{\xhi-\xhi\circ X_{-s}^\bv}_{L^1(\Omega)} =\norm{\xhi\circ X_s^\bv-\xhi}_{L^1(\Omega)}\\&\leq \sqrt d\left(\int_0^s\ e^{\tau\norm{\bv}_{\mathbf C^1(\overline{\Omega})}}\d \tau\right)\normmm{\nabla\xhi}(\Omega)\norm{\bv}_{\mathbf C(\overline{\Omega})},\quad \forall \xhi\in \mathcal M_{m_0},\quad \forall s\geq 0.
			\label{H2c}
		\end{align} 
\end{lemma}
\begin{proof}
    First recall that ${X_s^\bv}_\sharp\mathcal L^d=\mathcal L^d$ for any $s\in \mathbb R$, since $\Div \bv=0$ and $\bv\cdot \bn_{\pom}=0$ on $\partial\om$, so that, given $\xhi\in \Mm$ and $\phi\in \Hz$, with $\norm{\phi}_{\Hz}\leq1$,
		\begin{align*}
			\normmm{\int_\Omega (\xhi-\xhi\circ X_{-s}^\bv)\phi\dx}=\normmm{ 	\int_\Omega \xhi(\phi-\phi\circ X_{s}^\bv)\dx}\leq \norm{\xhi}_{L^\infty(\Omega)}\norm{\phi-\phi\circ X_{s}^\bv}_{L^1(\Omega)}\leq \norm{\phi-\phi\circ X_{s}^\bv}_{L^1(\Omega)}.
		\end{align*} 
		Moreover, we can show that
		\begin{align}
			\norm{\phi\circ X_s^\bv-\phi}_{L^1(\om)}\leq s\norm{\bv}_{\mathbf C(\overline{\Omega})}\norm{\phi}_{\Hz},\quad \forall \phi\in \Hz,
			\label{H1cv}
		\end{align}
		and thus conclude the proof of \eqref{H1c}. 
		To prove \eqref{H1cv}, observe that, as in \cite{Gigli}, by the Fundamental Theorem of Calculus and Fubini's Theorem,
		\begin{align*}
			\norm{\phi\circ X_s^\bv-\phi}^2&\leq \int_\Omega \normmm{\int_0^s{\partial_\tau (\phi\circ X_\tau^\bv)}\d \tau}^2\dx
			= \int_\Omega \normmm{\int_0^s{(\nabla \phi)\circ X_\tau^\bv\cdot \bv\circ X_\tau^\bv}\d \tau}^2\dx\\&
			\leq s\int_0^s\int_\Omega\normmm{(\nabla \phi)\circ X_\tau^\bv\cdot \bv\circ X_\tau^\bv}^2\dx\d \tau=
			s\int_0^s\int_\Omega\normmm{\nabla \phi\cdot \bv}^2\dx\d \tau\\&
			\leq s^2\norm{\bv}_{{\mathbf C(\overline{\Omega})}}^2\norm{\nabla \phi}^2,
		\end{align*}
		which entails \eqref{H1cv}. 
		
		To prove \eqref{H2c}, let us first consider $\phi\in C_c^1(\Omega)$ with $\norm{\phi}_{C^0(\overline{\Omega})}\leq 1$. Then we have, for $\xhi\in \Mm$, again by the Fundamental Theorem of Calculus and Fubini's Theorem,
		\begin{align}
			\nonumber	&\int_\Omega (\xhi\circ X_{s}^\bv-\xhi)\phi\dx \nonumber
			= 	\int_\Omega \xhi(\phi\circ X_{-s}^\bv-\phi)\dx\\&
			=\int_\Omega \chi\left(\int_{0}^{-s} \partial_\tau \phi\circ X_{\tau}^\bv \d \tau\right) \dx\nonumber
			\nonumber
			=-\int_0^s\int_\Omega \xhi\circ X_{\tau}^\bv(\text{div} (\phi\bv))\dx \d \tau \nonumber
			\\&\leq \norm{\phi\bv}_{\mathbf C(\overline{\Omega})}\int_0^s\normmm{\nabla (\xhi\circ  X_{\tau}^\bv)}(\Omega)\d \tau
			\leq \norm{\bv}_{\mathbf C(\overline{\Omega})}\int_0^s\normmm{\nabla (\xhi\circ  X_{\tau}^\bv)}(\Omega)\d\tau.\label{Bvconvergence}
		\end{align} 
		In order to conclude we need to find a uniform bound on the total variation of $\nabla (\xhi\circ  X_{\tau}^\bv)$. Let us first observe that, from \eqref{flowmap}, it holds, for any $s\geq0$, since $\Omega\subset \mathbb R^d$,
		\begin{align*}
			\normmm{\nabla X_s^\bv}&\leq \int_0^s\normmm{\nabla \partial_\tau X_\tau^\bv}\d\tau+\normmm{\nabla Id}
			\\&
			=\int_0^s\normmm{\nabla \bv(X_\tau^\bv)}\d\tau+\sqrt d
			\leq \int_0^s\norm{\bv}_{\mathbf C^1(\overline{\Omega})}\normmm{\nabla X_\tau^\bv}\d\tau+\sqrt d,
		\end{align*}
		entailing, by Gronwall's Lemma,
		\begin{align}
			\norm{\nabla X_s^\bv}_{\mathbf C(\overline{\Omega})}\leq \sqrt d\ e^{\int_0^s\norm{\bv}_{\mathbf C^1(\overline{\Omega})}\d\tau}=\sqrt d\ e^{s\norm{\bv}_{\mathbf C^1(\overline{\Omega})}},\quad \forall s\geq0.
			\label{Gronw}
		\end{align}
		Now, since $\chi\in \mathcal M_{m_0}$ by a well-known result for BV functions there exists a sequence of functions $\{\xhi_m\}_{m\in\mathbb N}\subset C^\infty(\Omega)\cap W^{1,1}(\om)$, such that $\xhi_m\to \xhi$ in $L^1(\Omega)$ and 
		\begin{align}
			\int_\Omega \vert \nabla \xhi_m\vert \dx\to \normmm{\nabla \chi}(\Omega),
			\label{cv}\end{align}
		as $m\to \infty$. Now, since ${X_s^\bv}_\sharp\mathcal L^d=\mathcal L^d$ for any $s\in \mathbb R$, it is immediate to infer
		\begin{align}
			\norm{\xhi_m\circ X_\tau^\bv-\xhi\circ X_\tau^\bv}_{L^1(\Omega)}=\norm{\xhi_m-\xhi}_{L^1(\Omega)}\to 0,\quad\text{ as }m\to \infty.\label{L1conv}
		\end{align}
		Moreover, it holds, from \eqref{Gronw},
		\begin{align}
			\nonumber\int_\Omega \vert \nabla (\xhi_m\circ X_\tau^\bv)\vert \dx&\leq \norm{\nabla X_\tau^\bv}_{\mathbf C(\overline{\Omega})}\int_\Omega \normmm{\nabla \xhi_m}\dx\\&
			\leq \sqrt d\ e^{\tau\norm{\bv}_{\mathbf C^1(\overline{\Omega})}}\int_\Omega \normmm{\nabla \xhi_m}\dx\leq C,\quad \forall \tau\geq0,\label{control1}
		\end{align}
		where $C>0$ is uniform in $m$ thanks to \eqref{cv}.
		This entails that the sequence $\xhi_m\circ X_\tau^\bv$ converges in $L^1(\Omega)$ to $\xhi\circ X_\tau^\bv$, and it has uniform-in-$m$ bounded total variation of $\nabla(\xhi_m\circ X_\tau^\bv)$. By \cite[Proposition 3.13]{AmbrosioFuscoPallara} this entails that $\xhi_m\circ X_\tau^\bv\to \xhi\circ X_\tau^\bv$ weakly* in BV, and thus, by lower semicontinuity of the total variation with respect to the convergence in measure, together with \eqref{cv} and \eqref{control1}, it holds
		\begin{align}
			\normmm{\nabla(\xhi\circ X_\tau^\bv)}(\Omega)\leq \liminf_{m\to \infty}\int_\Omega \vert \nabla (\xhi_m\circ X_\tau^\bv)\vert \dx\leq \sqrt d\ e^{\tau\norm{\bv}_{\mathbf C^1(\overline{\Omega})}}\normmm{\nabla \xhi}(\Omega),\quad \forall \tau\geq 0.
			\label{semic}
		\end{align}
		Plugging this estimate in \eqref{Bvconvergence} and taking the supremum over $\phi\in C_c^1(\Omega)$ with $\norm{\phi}_{C^0(\overline{\Omega})}\leq 1$, we infer 
		\begin{align*}
			\sup_{\phi\in C_c^1(\Omega),\   \norm{\phi}_{C^0(\overline{\Omega})}\leq 1}	\int_\Omega({\xhi\circ X_s^\bv-\xhi})\phi\leq \sqrt d\left(\int_0^s\ e^{\tau\norm{\bv}_{\mathbf C^1(\overline{\Omega})}}\d \tau\right)\normmm{\nabla\xhi}(\Omega)\norm{\bv}_{\mathbf C(\overline{\Omega})},
		\end{align*}
		entailing \eqref{H2c} by density of $C^1_c(\Omega)$ in $C_c^0(\Omega)$ together with Riesz Representation Theorem.
\end{proof}
Furthermore, in the following lemma we prove some convergence results, as $t\to0$, of the composition between a function in $\mathcal M_{m_0}$ and the flow map $X_t^\bv$. 
	\begin{lemma}
		Let $T>0$, $\{h_n\}_{n\in \mathbb{N}}\subset(0,1)$, with $h_n\to0$ as $n\to\infty$, and assume that a sequence $\{\bv_n\}_n \subset C([0,T];\mathbf C^1(\overline\Omega)\cap \Hs)$ is such that, for some $\gamma\in [0,1)$ and some constant $C>0$, 
			\begin{align}
				\norm{\bv_n}_{C([0,T];\mathbf C^1(\overline\Omega))}\leq Ch_n^{-\gamma},\quad \forall n\in \mathbb N.
					\label{unifC1}
			\end{align}
		Consider the map $X_{s}^{\bv_n(t)}$, $t\in(0,T)$, $s\in[-h_n,h_n]$, which is a solution to 
		\begin{align}
			\begin{cases}\frac{\d}{\d s} X_{s}^{\bv_n(t)}=\bv_n(X_{s}^{\bv_n(t)}, t),\\
				X_0=Id.
			\end{cases}
			\label{map1}
		\end{align}
	
		Then for any bounded sequence $\{w_n\}_n\subset L^2(0,T;\Hz)$, such that $w_n\rightharpoonup w$ weakly in $L^2(0,T;\Hz)$ as $n\to\infty$, it holds 
		\begin{align}
			w_n\circ X_{\cdot-\floor{\frac \cdot{h_n}}}^{\bv_n(\cdot)}\rightharpoonup  w,\quad\text{ weakly in }L^2(0,T;\Hz)\text{ as }n\to\infty,
			\label{w}
		\end{align}
		and 
		\begin{align}
			\nabla	(w_n\circ X_{\cdot-\floor{\frac \cdot{h_n}}}^{\bv_n(\cdot)})\rightharpoonup \nabla w,\quad\text{ weakly in }L^2(0,T;\mathbf L^2(\Omega))\text{ as }n\to \infty.
			\label{w1}
		\end{align}
		Furthermore, if, additionally to \eqref{unifC1}, it holds
		\begin{align}
		\bv_n\to\bv\quad\text{ strongly in }L^2(0,T;\mathbf L^2_\sigma(\om))\text{ as }n\to\infty,\qquad 	\norm{\bv}_{L^3(\om\times(0,T))}\leq C,
			\label{b1}
		\end{align}
		for some $C>0$, then for any sequence $\{\xhi_n\}_n\subset L^\infty(0,T;\mathcal M_{m_0})$, such that there exists $C>0$ so that
		\begin{align}
			\esssup_{t\in(0,T)}\normmm{\nabla \xhi_n(t)}(\om)\leq C,\quad\forall n\in\N,
			\label{b2}
		\end{align}
		as well as it holds
		\begin{align}
			\xhi_n\to \xhi\quad\text{ strongly in }L^1(\om\times(0,T)),\label{strongL1}
		\end{align}
 as $n\to\infty$, then  
		\begin{align}
			\frac{\xhi_n-\xhi_n\circ X_{-h_n}^{\bv_n(\cdot)}}{h_n}\to \bv\cdot \nabla \xhi\quad\text{ strongly
				in }L^2(0,T;\Hmz)\text{ as }n\to\infty,
			\label{ww1}
		\end{align}
		where $\bv\cdot \nabla\xhi\in H^{-1}_{(0)}(\Omega)$ means
		\begin{align*}
			\langle\bv\cdot \nabla\xhi, \varphi\rangle_{\Hmz,\Hz}=-(\bv\cdot \nabla \varphi, \xhi)_\om,
		\end{align*}
		for any $\varphi\in \Hz$.
		\label{convh1}
	\end{lemma}
	\begin{proof}
		To prove \eqref{w}, let us consider a bounded sequence $\{w_n\}_n\subset L^2(0,T;\Hz)$ such that $w_n\rightharpoonup w$ weakly in $L^2(0,T;\Hz)$ as $n\to\infty$, which clearly entails $w_n\rightharpoonup w$ weakly in $L^2(\Omega\times(0,T))$. Now, repeating exactly the arguments leading to \eqref{H1cv}, applied to $X_{t-\floor{\frac t{h_n}}}^{\bv_n(t)}$, we get, for any $\xi\in L^2(0,T;\Hz)$, for almost any $t\in(0,T)$,
		\begin{align*}
			\norm{\xi( X_{t-\floor{\frac t{h_n}}}^{\bv_n(t)},t)-\xi(t) }^2&
			=\norm{\int_0^{t-\floor{\frac t{h_n}}}\frac{\d}{\d s}\xi(X_{s}^{\bv_n(t)},t)\d s}^2\\&
			=\norm{\int_0^{t-\floor{\frac t{h_n}}}\nabla \xi(X_{s}^{\bv_n(t)},t)\cdot \bv_n(X_{s}^{\bv_n(t)},t)\d s}^2
			\\&\leq (t-\floor{\frac t{h_n}})^2\norm{\bv_n}^2_{C([0,T];\mathbf C(\overline\Omega))}\norm{\nabla\xi(t)}^2\\&\leq h_n^2\norm{\bv_n}_{C([0,T];\mathbf C(\overline\Omega))}^2\norm{\nabla\xi(t)}^2,
		\end{align*}
		so that it holds, for any $\phi\in L^2(\Omega\times(0,T))$,
		\begin{align*}
			&\normmm{\int_0^T({w_n\circ X_{t-	\floor{\frac{t}{h_n}}}^{\bv_n(t)}-w},\phi)_\om\dt}= \normmm{\int_0^T(w_n\circ X_{t-	\floor{\frac{t}{h_n}}}^{\bv_n(t)}-w_n,\phi)_\om\dt+\int_0^T({w_n-w},\phi)_\om\dt}\\&
			\leq h_n\norm{\bv_n}_{C([0,T];\mathbf C(\overline{\Omega}))}\norm{w_n}_{{L^2(0,T;\Hz)}}\norm{\phi}_{L^2(0,T;L^2(\Omega))}+\normmm{\int_0^T({w_n-w},\phi)\dt}\\&
			\leq Ch_n^{1-\gamma}+\normmm{\int_0^T({w_n-w},\phi)\dt} \to 0,\quad \text{as }n\to\infty,
		\end{align*}
		since $h_n\to0$, $\gamma\in[0,1)$, and, due to \eqref{unifC1}, it holds $$\norm{\bv_n}_{C(0,T;\mathbf C(\overline{\Omega}))}\sup_{n}\norm{w_n}_{_{L^2(0,T;\Hz)}}\leq Ch_n^{-\gamma}.$$
		Therefore, 
		\begin{align}
			w_n\circ  X_{\cdot-	\floor{\frac{	\cdot}{h_n}}}^{\bv_n(\cdot)}\rightharpoonup w\text{ weakly in }L^2(\Omega\times(0,T)),\label{weak}
		\end{align}
		as $n\to\infty$. Repeating the same steps leading to \eqref{Gronw}, this time with the velocity $\bv(t)$, we infer 
		\begin{align}
			\nonumber\sup_{s\in[0,t-\floor{\frac{t}{h_n}}]}\norm{\nabla X_s^ {\bv_n(t)}}_{\mathbf C(\overline{\Omega})}&\leq \sqrt d\ e^{\int_0^{t-\floor{\frac{t}{h_n}}}\norm{\bv_n(t)}_{\mathbf C^1(\overline{\Omega})}\d s}\\&\nonumber
			\leq \sqrt d\ e^{h_n\norm{\bv_n}_{C([0,T];\mathbf C^1(\overline{\Omega}))}}\\&
			\leq \sqrt d\ e^{Ch_n^{1-\gamma}}, \quad \forall t\in[0,T].
			\label{Gronw2}
		\end{align}
		Therefore, it holds (this can be shown for instance by approximating $w_n$ with smooth functions),
		\begin{align}
			&\nonumber	\int_0^T\normm{\nabla (w_n\circ X_{t-\floor{\frac{t}{h_n}}}^{\bv_n(t)})}^2\dt\leq \int_0^T\norm{\nabla X_{t-\floor{\frac{t}{h_n}}}^{\bv_n(t)}}_{\mathbf C(\overline{\Omega})}^2\norm{\nabla w_n}^2\dt\\&\leq {d}\ e^{2Ch_n^{1-\gamma}}\sup_{n}\int_0^T\norm{w_n}_{\Hz}^2\dt\leq C<+\infty,\quad \text{ for any }n\in \mathbb N,\label{dr}
		\end{align}
		for $\gamma\in[0,1)$, due to \eqref{unifC1}.
		As a consequence, we immediately infer that, since $w_n\circ X_{\cdot-\floor{\frac{\cdot}{h_n}}}^{\bv_n(\cdot)}$ is bounded in $L^2(0,T;\Hz)$, it also converges weakly, up to a subsequence, in the same space, and thus we deduce from \eqref{weak} that $w_n\circ X_{\cdot-\floor{\frac{\cdot}{h_n}}}^{\bv_n(\cdot)}\rightharpoonup w$ weakly in $L^2(0,T;\Hz)$. By the uniqueness of the weak limit we also deduce that the entire sequence $\{w_n\circ X_{\cdot-\floor{\frac{\cdot}{h_n}}}^{\bv_n(\cdot)}\}_n$ weakly converges to the same limit, entailing \eqref{w} and, by a similar argument, \eqref{w1}.
		
			In order to prove \eqref{ww1}, let us first observe that, for almost any $t\in(0,T)$,
		\begin{align*}
			\xhi_n(t)\circ X_{(\cdot)}^{\bv_n(t)}\in C([-1,1];L^p(\om)),\quad \forall p\in [1,\infty).
		\end{align*}
		Indeed, it holds, exploiting \eqref{map1} and by the Fundamental Theorem of Calculus and Fubini's Theorem (the computations are formal, but they can be easily obtained after regularizing $\xhi_n(t)$, as to obtain \eqref{control1}), for $t\in(0,T)$, $r,s\in[-1,1]$, $r\leq s$,
		\begin{align}
			\nonumber	&\int_\Omega\normmm{\xhi_n(t)\circ X_{s}^{\bv_n(t)}-\xhi_n(t)\circ X_{r}^{\bv_n(t)}}\dx \\&\nonumber
			\leq \int_r^s\int_\om\normmm{\nabla \xhi_n\circ X_{\tau}^{\bv_n(t)}\cdot \bv_n(t)\circ X_{\tau}^{\bv_n(t)}}\dx\d\tau
			\\&\label{hn}
			=\int_r^s\int_\om\normmm{\nabla \xhi_n\cdot \bv_n(t)}\dx\d\tau
			\leq (s-r)\normmm{\nabla \xhi_n(t)}(\Omega)\norm{\bv_n(t)}_{\mathbf C(\overline{\Omega})}\leq Ch_n^{-\gamma}\normmm{s-r},
		\end{align}
		where we used \eqref{unifC1} and \eqref{b2}. This gives $\xhi_n(t)\circ X_{(\cdot)}^{\bv_n(t)}\in C([-1,1];L^1(\om))$. Continuity in $L^p(\om)$, $p\in[1,\infty)$, is then obtained by interpolation, since $\xhi_n(t)\in\{0,1\}$ almost everywhere in $\om$. Now, due to continuity, we immediately have, by the Mean Value Theorem, that
		\begin{align}
		\int_{0}^{h_n}	\frac{ \xhi_n(t)\circ X_{-\tau}^{\bv_n(t)}}{h_n}\d\tau=\xhi_n(t)\circ X_{-q_n(t)}^{\bv_n(t)},
			\label{strconv}
		\end{align} 
		for some $0\leq q_n(t)\leq {h_n}$. Now, exploiting again \eqref{hn}, we have that 
		\begin{align}
			\nonumber	&\int_\Omega\normmm{\xhi_n(t)\circ X_{-q_n(t)}^{\bv_n(t)}-\xhi_n(t)}\dx 
			\leq Ch_n^{-\gamma}q_n(t)\leq Ch_n^{1-\gamma}\to 0,
		\end{align}
		for almost any $t\in(0,T)$. As a consequence, recalling \eqref{strconv}, by Lebesgue's Dominated Convergvence Theorem we easily infer that 
		\begin{align}
			\label{DOM0}	&\int_0^T\int_\Omega\normmm{\int^{h_n}_0\frac{ \xhi_n(t)\circ X_{-\tau}^{\bv_n(t)}}{h_n}\d\tau-\xhi_n(t)}\dx\dt \to 0,
		\end{align}
		entailing, again by interpolation,
		\begin{align}
			\label{DOM}	&\int^{h_n}_0\frac{ \xhi_n(\cdot)\circ X_{-\tau}^{\bv_n(\cdot)}}{h_n}\d\tau-\xhi_n(\cdot) \to 0,
		\end{align}
		strongly in $L^p(\om\times(0,T))$, for any $p\in[1,\infty)$, as $n\to\infty$.
		
		Observe now that it holds, again by the Fundamental Theorem of Calculus and Fubini's Theorem, 
		\begin{align}
		\nonumber	&	\int_0^T\int_\om \frac{\xhi_n(t)-\xhi_n(t)\circ X_{-h_n}^{\bv_n(t)}}{h_n}\varphi\dx\dt
			= \int_0^T\int_\om \xhi_n(t)\frac{\varphi-\varphi\circ X_{h_n}^{\bv_n(t)}}{h_n}\dx\dt\\&\nonumber
			=-\frac1{h_n}\int_0^T\int_\om \xhi_n(t)\int_0^{h_n}\nabla \varphi\circ X_{\tau}^{\bv_n(t)}\cdot \bv_n(t)\circ X_{\tau}^{\bv_n(t)}\d\tau\dx\dt\\&\nonumber
			=-\int_0^T\int_0^{h_n}\int_\om \frac{\xhi_n(t)\circ X_{-\tau}^{\bv_n(t)}}{h_n}\nabla \varphi\cdot \bv_n(t)\dx\d\tau\dt\\&
			=-\int_0^T\int_\om \nabla \varphi\cdot \bv_n(t)\int_0^{h_n} \frac{\xhi_n(t)\circ X_{-\tau}^{\bv_n(t)}}{h_n}\d\tau\dx\dt,\label{meanv}
		\end{align}	
		for any $\varphi\in \Hz$.
		We can thus compare this integral with the one we expect from the limit, obtaining, for $\varphi\in L^2(0,T;\Hz)$, 
		\begin{align*}
			&\normmm{\int_0^T\int_\om \left(\frac{\xhi_n(t)-\xhi_n(t)\circ X_{-h_n}^{\bv_n(t)}}{h_n}\varphi(t)+\xhi(t)\bv(t)\cdot\nabla \varphi(t)\right)\dx\dt}\\&
			=\normmm{\int_0^T\int_\om\left( -\nabla \varphi(t)\cdot \bv_n(t)\int_0^{h_n} \frac{\xhi_n(t)\circ X_{-\tau}^{\bv_n(t)}}{h_n}\d\tau+\xhi(t)\bv(t)\cdot\nabla \varphi(t)\right)\dx\dt}\\&
			\leq \norm{ \varphi}_{L^2(0,T;\Hz)}\norm{\bv}_{L^3(\om\times(0,T))}\norm{\int_0^{h_n} \frac{\xhi_n(\cdot)\circ X_{-\tau}^{\bv_n(\cdot)}}{h_n}\d\tau-\chi_n(\cdot)}_{L^6(\om\times(0,T))}\\&
			\quad +\norm{ \varphi}_{L^2(0,T;\Hz)}\norm{\bv_n-\bv}_{L^2(0,T;\bL^2(\om))}\norm{\int_0^{h_n} \frac{\xhi_n(\cdot)\circ X_{-\tau}^{\bv_n(\cdot)}}{h_n}\d\tau}_{L^\infty(\om\times(0,T))}\\&\quad +\norm{ \varphi}_{L^2(0,T;\Hz)}\norm{\bv}_{L^3(\om\times(0,T))}\norm{\xhi-\xhi_n}_{L^6(\om\times(0,T))},
		\end{align*}
		and thus, by using $\langle \frac{\xhi_n(t)-\xhi_n(t)\circ X_{-h_n}^{\bv_n(t)}}{h_n},\varphi\rangle_{\Hmz,\Hz}=\int_\om \frac{\xhi_n(t)-\xhi_n(t)\circ X_{-h_n}^{\bv_n(t)}}{h_n}\varphi\dx$, for $\varphi\in \Hz$, we infer 
		\begin{align*}
			&\sup_{\varphi \in L^2(0,T;\Hz)}\normmm{\int_0^T\langle \frac{\xhi_n(t)-\xhi_n(t)\circ X_{-h_n}^{\bv_n(t)}}{h_n}-\nabla\xhi(t)\cdot\bv(t),\varphi(t)\rangle_{\Hmz,\Hz}\dt}\\&\leq C\left(\norm{\int_0^{h_n} \frac{\xhi_n(\cdot)\circ X_{-\tau}^{\bv_n(\cdot)}}{h_n}\d\tau-\chi_n(\cdot)}_{L^6(\om\times(0,T))}\right.\\&\quad\left.+\norm{\bv_n-\bv}_{L^2(0,T;\bL^2(\om))}+\norm{\xhi-\xhi_n}_{L^6(\om\times(0,T))}\right)\to 0,
		\end{align*}
		as $n\to\infty$, where we have used \eqref{b1}, \eqref{strongL1} (which holds in $L^6(\om\times(0,T))$ by interpolation), as well as \eqref{DOM} with $p=6$. Note that we have also exploited the fact that, since $\xhi_n\in\{0,1\}$ almost everywhere, 
		\begin{align*}
		\norm{\int_0^{h_n} \frac{\xhi_n(\cdot)\circ X_{-\tau}^{\bv_n(\cdot)}}{h_n}\d\tau}_{L^\infty(\om\times(0,T))}\leq 1.
		\end{align*}
		
		 This convergence exactly corresponds to \eqref{ww1}, which is thus proven.
	\end{proof}
	In conclusion, we also need a specific chain rule for the composition $\xhi\circ X_t^\bv$ for $\xhi\in \mathcal M_{m_0}$. In particular, we have
	\begin{lemma}
		Let $\xhi_A,\xhi\in L^\infty(\Omega)$, with $\xhi_A,\xhi\in \{0,1\}$ almost eveywhere in $\om$, $\int_\Omega \chi_A\dx=\int_\Omega \chi\dx=m_0$, and $\bv\in \mathbf C^1(\overline{\Omega})\cap \Hs$. Then, for almost any $t\in \mathbb R$ it holds
		\begin{align}
			\frac12\ddt\norm{\xhi_A-\xhi\circ X_{-t}^\bv}_{\Hmz}^2=-\left(\xhi\bv, \nabla ((-\Delta_N)^{-1}(\xhi_A-\xhi\circ X_{-t}^\bv)\circ X_t^\bv)\right)_\om\label{cv1}
		\end{align}
		\label{convH1}
	\end{lemma}
	\begin{proof}
		Since $\xhi\in L^\infty(\Omega)$, there exists a sequence of functions $\{\xhi_m\}_{m\in\mathbb N}\subset H^2(\Omega)$, with $0\leq \xhi_m\leq 1$ and $\int_\Omega \xhi_m\dx=m_0$, such that $\xhi_m\to \xhi$ in $L^2(\Omega)$ as $m\to\infty$. Indeed, it is enough to choose:
		\begin{align*}
			\xhi_m:=(I-\frac1m\Delta_N)^{-1}\xhi \in H^2(\Omega),
		\end{align*}
		since, by maximum principle, $0\leq \xhi_m\leq 1$, and $\xhi_m$ satisfies, pointwise almost everywhere,
		\begin{align*}
			-\frac1m\Delta \xhi_m+\xhi_m=\xhi,
		\end{align*}
		so that, integrating over $\Omega$, we get $\int_\Omega \xhi_m\dx=\int_\Omega\xhi\dx=m_0$. Note now that $\int_\Omega(\xhi_A-\xhi_m\circ X_{-t}^\bv)\dx=\int_\Omega(\xhi_A-\xhi_m)\dx=0$, and thus we can compute 
		\begin{align*}
			\onehalf\ddt\normh{\chi_A-\xhi_m\circ X_{-t}^\bv}^2&=-(\pt(\xhi_m\circ X_{-t}^\bv), (-\Delta_N)^{-1}(\chi_A-\xhi_m\circ X_{-t}^\bv))_\om\\&
			=(\partial_{s}(\xhi_m\circ X_{s}^\bv)_{|s=-t}, (-\Delta_N)^{-1}(\chi_A-\xhi_m\circ X_{-t}^\bv))_\om
			\\&
			=((\nabla \xhi_m) \circ X_{-t}^\bv\cdot \bv\circ X_{-t}^\bv, (-\Delta_N)^{-1}(\chi_A-\xhi_m\circ X_{-t}^\bv))_\om\\&
			= (\nabla \xhi_m\cdot \bv, ((-\Delta_N)^{-1}(\chi_A-\xhi_m\circ X_{-t}^\bv)) \circ X_{t}^\bv)_\om\\&
			=-(\xhi_m \bv, \nabla ((-\Delta_N)^{-1}(\chi_A-\xhi_m\circ X_{-t}^\bv) \circ X_{t}^\bv))_\om,
		\end{align*}
		so that, integrating over $(s,t)$, $0\leq s<t$, we get  
		\begin{align}
			&	\nonumber\onehalf\normh{\chi_A-\xhi_m\circ X_{-t}^\bv}^2-\onehalf	\normh{\chi_A-\xhi_m\circ X_{-s}^\bv}^2\\&
			=-\int_s^t(\xhi_m \bv, \nabla ((-\Delta_N)^{-1}(\chi_A-\xhi_m\circ X_{-\tau}^\bv) \circ X_{\tau}^\bv))_\om\d \tau.\label{final1}
		\end{align}
		Now, the left-hand side easily passes to the limit since, for any $t\in \mathbb R$ and any $\phi\in \Hz$, 
		\begin{align*}
			(\chi_A-\xhi_m\circ X_{-t}^\bv-(\xhi_A-\chi\circ X_{-t}^\bv),\phi)_\om=(\xhi_m-\xhi, \phi\circ X^\bv_{t})_\om\leq \norm{\xhi_m-\xhi}\norm{\phi}_{\Hz},
		\end{align*}
		and $\xhi_m\to \xhi$ in $L^2(\Omega)$ by construction. Concerning the right-hand side, we first notice that
		$$
		\norm{\xhi_m\circ X_{-t}^\bv-\xhi\circ X_{-t}^\bv}=\norm{\xhi-\xhi_m}\to 0,\text{ as } m\to \infty.
		$$
		Now, arguing as for estimate \eqref{dr}, since $\bv\in \mathbf C^1(\overline\om)$, it is immediate to infer 
		\begin{align*}
		&\norm{\nabla((-\Delta_N)^{-1}((\xhi_A-\xhi_m\circ X_{-t}^\bv)-(\xhi_A-\xhi\circ X_{-t}^\bv))\circ X_t^\bv)}\\&
		=\norm{\nabla((-\Delta_N)^{-1}(\xhi_m\circ X_{-t}^\bv-\xhi\circ X_{-t}^\bv)\circ X_t^\bv)}\\&
		\leq C(t)\norm{\nabla (-\Delta_N)^{-1}(\xhi_m\circ X_{-t}^\bv-\xhi\circ X_{-t}^\bv)}\\&\leq C(t)\norm{\xhi_m\circ X_{-t}^\bv-\xhi\circ X_{-t}^\bv}\to 0,\text{ as }m\to \infty,
		\end{align*}
		so that 
		$$
		\nabla ((-\Delta_N)^{-1}(\chi_A-\xhi_m\circ X_{-t}^\bv) \circ X_{t}^\bv)\to \nabla ((-\Delta_N)^{-1}(\chi_A-\xhi\circ X_{-t}^\bv) \circ X_{t}^\bv),\quad \text{as }m\to \infty,\text{ strongly in }\mathbf L^2(\Omega),
		$$
		for any $t\in\R$.
		Since also $\xhi_m\to \xhi$ strongly in $L^2(\Omega)$ and $\bv\in \mathbf C(\overline{\Omega})$, we have 
		$$
		(\xhi_m \bv, \nabla ((-\Delta_N)^{-1}(\chi_A-\xhi_m\circ X_{-\tau}^\bv) \circ X_{\tau}^\bv))_\om\to (\xhi \bv, \nabla ((-\Delta_N)^{-1}(\chi_A-\xhi\circ X_{-\tau}^\bv) \circ X_{\tau}^\bv))_\om
		$$
		as $m\to \infty$ for any $\tau\in(s,t)$. Therefore, by the Dominated Convergence Theorem we immeditaley infer that
		$$
		\int_s^t(\xhi_m \bv, \nabla ((-\Delta_N)^{-1}(\chi_A-\xhi_m\circ X_{-\tau}^\bv) \circ X_{\tau}^\bv))_\om\d\tau\to \int_s^t(\xhi \bv, \nabla ((-\Delta_N)^{-1}(\chi_A-\xhi\circ X_{-\tau}^\bv) \circ X_{\tau}^\bv))_\om\d \tau.
		$$
		We can thus pass to the limit in \eqref{final1}, obtaining
		\begin{align*}
			&	\nonumber\onehalf\normh{\chi_A-\xhi\circ X_{-t}^\bv}^2-\onehalf	\normh{\chi_A-\xhi\circ X_{-s}^\bv}^2\\&
			=-\int_s^t(\xhi \bv, \nabla ((-\Delta_N)^{-1}(\chi_A-\xhi\circ X_{-\tau}^\bv) \circ X_{\tau}^\bv))_\om\d \tau,\quad \forall 0\leq s<t,
		\end{align*}
		from which \eqref{cv1} follows.
		
	\end{proof}

\section{Proof of Theorem \ref{thm1}}
\label{proofthm1}
The proof is quite technical and it is subdivided in many subsections for the reader's convenience.
\subsection{Approximating velocity}
Let us fix $h>0$. For $0<b$, given $\bv_h\in L^2(0,b;\bL^2_\sigma(\om))$, we define $\bu^k_h:\Omega\times(0,b)\to \R^d$, $k\in \mathbb N$, as the solution to:
\begin{align}
	\frac1k\pt\bA^2 \bu^k_h+\frac1k\bA^2\bu^k_h+\bu^k_h=\bP_\sigma M_k(\bv_h),\quad \bu^k_h(0)=\mathbf 0\in \bH^4(\om)\cap \Hus,
	\label{bvk}
\end{align}
where, for any component $i=1,\ldots,d$, 
\begin{align}
	M_k(\bv_h)_i:=\begin{cases}
		k,\quad &\text{ if }\bv_{h,i}\geq k,\\
		\bv_{h,i},\quad&\text{ if }-k<\bv_{h,i}< k,\\
		-k,\quad &\text{ if }\bv_{h,i}\leq -k.
	\end{cases}
	\label{cutoff0}
\end{align}
The existence of a solution to \eqref{bvk} can be easily shown by means of a Galerkin scheme, whereas uniqueness directly comes from linearity together with standard energy estimates. We only point out that, by multiplying \eqref{bvk} by $\pt\bA^2\bu^k_h$ we obtain 
\begin{align}
	\nonumber &\frac12\ddt\left(\norm{\bA\bu^k_h}^2+\frac1k\norm{\bA^2\bu^k_h}^2\right)	+\frac1k\norm{\pt\bA^2\bu^k_h}^2\\&
	\leq \norm{\bP_\sigma M_k(\bv_h)}\norm{\pt\bA^2\bu^k_h}\nonumber
	\leq \frac{1}{2k}\norm{\pt\bA^2\bu^k_h}^2+\frac{k}{2}\norm{M_k(\bv_h)}^2\\&
	\leq \frac{1}{2k}\norm{\pt\bA^2\bu^k_h}^2+C{k^3},\label{ctrl0}
\end{align}
entailing, by Gronwall's Lemma,
\begin{align}
	\norm{\bu^k_h}_{C^\frac12([0,b];\mathbf C^1(\overline{\om}))}\leq C\norm{\bu^k_h}_{ H^1(0,b;\bH^4(\om)\cap \Hus)}\leq  C(b)k^2.
	\label{basic0}
\end{align}
Notice that here the essential property of this approximating scheme is that the higher-order norm above is controlled \textit{independently} of $\bv_h$, and only depends on $k$. This is possible since the lower-order norm which is needed to be under control, uniformly in $k$, by means of $\bv_h$, is only the $L^2$-norm (see \eqref{controls} below). All the other norms can be approximated without any necessity of a relation with $\bv_h$, and thus they can be set to depend \textit{only} on $k$.
From now on we crucially set $$k:=k(h)=h^{-\frac1{8}},$$ so that from \eqref{basic0} we infer
\begin{align}
	\norm{\bu^k_h}_{C^\frac12([0,b];\mathbf C^1(\overline{\om}))}\leq C\norm{\bu^k_h}_{ H^1(0,b;\bH^4(\om)\cap \Hus)}\leq  C(b)h^{-\frac14}.
	\label{basic}
\end{align}
For the sake of clarity, to distinguish between the time discretization scale and the $k$ scale, we will still distinguish between $h$ and $k$, having in mind that $k$ depends on $h$.

We then set $\bu_n^k:=\bu^k_h(nh)$, for $n\in \N$, $nh\leq b$, together with $\bu_0^k:=\mathbf 0$, so that we have a sequence $\{\bu_n^k\}_n\subset\bH^4(\om)\cap \Hus\hookrightarrow \mathbf C^1(\overline\Omega)$, for any fixed $h>0$.\label{approxv}
	\subsection{First discretized step}\label{firstminmovement}
	Set, for $s\geq0$, 
	$$T_h(s):=s-(n-1)h,\quad s\in((n-1)h,nh],$$
	where $n=\ceil{\tfrac sh}$, so that $T_h((mh)^-):=\lim_{s\to (mh)^-}T_h(s)=h$, for any $m\in\N$, and $T_h(t)=t$ for $t\in(0,h]$.
	\subsubsection{Minimizing movement step}
	 We introduce the following minimizing problem,
	\begin{align}
		&\chi^h_{0}=\chi_0,\\
		&\chi_{1}^{h}\in\arg\min_{\chi\in \mathcal M_{m_0}}\left(E[\chi]+\frac 1{2h}\normh{\chi-\chi_{0}^{h}\circ X_{-h}^{\bu_0^k}}^2\right),\label{M1}
	\end{align}
	with, as in \eqref{E1}, 
	$$
	E[\xhi]=\int_\Omega c_0d\normmm{\nabla \xhi}+\int_{\partial\Omega} c_0\cos\gamma\chi \d \mathcal H^{d-1},
	$$
	$c_0>0$, $\gamma\in(0,\frac\pi2]$, and where $\bu_0^k=\mathbf 0\in \mathbf H^4_\sigma(\om)$ is defined as in Section \ref{approxv}. Here we keep this quantity even if it is zero, since the same argument will be repeated for any time step $n$, then involving the generic nontrivial $\bu_n^k$.
	Then, we need a further interpolation, as follows:
	\begin{align}
		& \xhi^{T_h(0)}_{0}:=\xhi_{0}^h,\\&
		\xhi_t^{T_h(t)}\in\arg\min_{\chi\in \mathcal M_{m_0}}\left(E[\chi]+\frac 1{2t}\normh{\chi- \xhi^{h}_{0}\circ X_{-T_h(t)}^{\bu_0^k}}^2\right),
		\text{ for }t\in (0,h].\label{M2}
	\end{align}
	For any fixed $t$, the problem clearly admits at least one solution by the Direct Method of the Calculus of Variations. For instance, see \cite[ Proposition 1.2]{Modica} for the lower-semicontinuity of the capillary energy.

Let us now define the function 
 \begin{align}
 	f\tTh:=E[\xhi_t^{T_h(t)}]+\frac{1}{2t}\normh{\chi\tTh-\xhi_0\circ X_{-T_h(t)}^{\bu_0^k}}^2,
 	\label{ftt}
 \end{align}
 where $\xhi_t^\Th{t}$ is defined in \eqref{M2}.
 Now, we aim at showing that, for almost any $t\in(0,h]$,
 \begin{align}
 	\nonumber \ddt f\tTh&=-\frac{1}{2t^2}\normh{\chi\tTh-\chi_0\circ X_{-T_h(t)}^{\bu_0^k}}^2\\&\quad -\frac{1}{t}\left( \xhi_{0}\cdot \bu_0^k,\nabla \left((-\Delta_N)^{-1}(\chi\tTh-\xhi_0\circ X_{-T_h(t)}^{\bu_0^k})\circ X_{T_h(t)}^{\bu_0^k}\right)\right)_\om.\label{derivative}
 \end{align}
 \begin{proof}[Proof of \eqref{derivative}]
 	The main steps of the proof are the following:
 	\begin{itemize}
 		\item Show that $f\tTh$ is locally Lipschitz in time. 
 		
 		Indeed, letting $t,s\in (0,h]$, by the optimality properties \eqref{M2}, we have
 		\begin{align}
 			& \nonumber f\tTh-f_s^\Th{s}\\&\nonumber\leq E[\xhi_s^\Ths]+\frac{1}{2t}\normh{\chi_s^\Ths-\xhi_0\circ X_{-T_h(t)}^{\bu_0^k}}^2 \\&\quad \nonumber-E[\xhi_s^\Ths]-\frac{1}{2s}\normh{\chi_s^\Ths-\xhi_0\circ X_{-T_h(s)}^{\bu_0^k}}^2\\&\nonumber
 			=\left(\frac1{2t}-\frac1{2s}\right)\normh{\chi_s^\Ths-\xhi_0\circ X_{-T_h(t)}^{\bu_0^k}}^2\nonumber\\&\quad+\frac1{2s}\left(\normh{\chi_s^\Ths-\xhi_0\circ X_{-T_h(t)}^{\bu_0^k}}^2-\normh{\chi_s^\Ths-\xhi_0\circ X_{-T_h(s)}^{\bu_0^k}}^2\right).\label{est}
 		\end{align}
 Now first note that, since $\bu_0^k\in \mathbf C^1(\overline{\Omega})$, by \eqref{Gronw} it holds 
 		\begin{align}
 			\nonumber	&\norm{\nabla\left( (-\Delta_N)^{-1}\left(\chi_s^\Ths-\xhi_0\circ X_{-T_h(r)}^{\bu_0^k}\right)\circ X_{T_h(r)}^{\bu_0^k}\right)}\\\nonumber&\leq C(h,\norm{\bu_0^k}_{\mathbf C^1(\overline\Omega)})\norm{\left(\nabla (-\Delta_N)^{-1}\left(\chi_s^\Ths-\xhi_0\circ X_{-T_h(r)}^{\bu_0^k}\right)\right)\circ X_{-T_h(t)}^{\bu_0^k}}\\&\nonumber
 			=C(h,\norm{\bu_0^k}_{\mathbf C^1(\overline\Omega)})\norm{\nabla (-\Delta_N)^{-1}\left(\chi_s^\Ths-\xhi_0\circ X_{-T_h(r)}^{\bu_0^k}\right)}\\&
 			\leq C(h,\norm{\bu_0^k}_{\mathbf C^1(\overline\Omega)})\normh{\chi_s^\Ths-\xhi_0\circ X_{-T_h(t)}^{\bu_0^k}}.\label{Nm1}
 		\end{align}
 		Moreover, we have
 		\begin{align*}
 			\left(\chi_s^\Ths-\xhi_0\circ X_{-T_h(t)}^{\bu_0^k},\varphi\right)_\om\leq C\norm{\varphi}_{\Hz},
 		\end{align*}
 		so that 
 		$$
 		\normh{\chi_s^\Ths-\xhi_0\circ X_{-T_h(t)}^{\bu_0^k}}\leq C,
 		$$
 		and thus
 		\begin{align}
 			&\norm{\nabla\left( (-\Delta_N)^{-1}\left(\chi_s^\Ths-\xhi_0\circ X_{-T_h(r)}^{\bu_0^k}\right)\circ X_{T_h(r)}^{\bu_0^k}\right)}\leq C(h).\label{estim}
 		\end{align}
 		Note now that, by Lemma \ref{convH1} and recalling that $T_h(r)=r$ since $r\in(0,h]$, for almost any $r\in(0,h]$ it holds,
 		\begin{align}
 			&\nonumber\onehalf\partial_r\normh{\chi_s^\Ths-\xhi_0\circ X_{-T_h(r)}^{\bu_0^k}}^2\\&
 			=-\left( \chi_0\cdot \bu_0^k, \nabla\left( (-\Delta_N)^{-1}\left(\chi_s^\Ths-\xhi_0\circ X_{-T_h(r)}^{\bu_0^k}\right)\circ X_{T_h(r)}^{\bu_0^k}\right)\right)_\om,
 			\label{derivativeA}
 		\end{align}
 		so that, recalling \eqref{estim},
 		\begin{align*}
 			&\onehalf\normmm{\partial_r\normh{\chi_s^\Ths-\xhi_0\circ X_{-T_h(r)}^{\bu_0^k}}^2}\\&\leq \norm{\chi_0}_{L^\infty(\Omega)}\norm{\bu_0^k}\norm{\nabla\left( (-\Delta_N)^{-1}\left(\chi_s^\Ths-\xhi_0\circ X_{-T_h(r)}^{\bu_0^k}\right)\circ X_{T_h(r)}^{\bu_0^k}\right)}
 			\leq C.
 		\end{align*}
 		Therefore, we obtain 
 		\begin{align*}
 			& \frac1{2s}\left(\normh{\chi_s^\Ths-\xhi_0\circ X_{-T_h(t)}^{\bu_0^k}}^2-\normh{\chi_s^\Ths-\xhi_0\circ X_{-T_h(s)}^{\bu_0^k}}^2\right)\\&\leq 
 			\frac1{2s}\int_s^t\normmm{\partial_r\normh{\chi_s^\Ths-\xhi_0\circ X_{-T_h(r)}^{\bu_0^k}}^2}dr\\&
 			\leq \frac{C}{s}\normmm{t-s},
 		\end{align*}
 		and thus, plugging into \eqref{est}, we easily obtain
 		\begin{align*}
 			f\tTh-f_s^\Th{s}\leq \frac{C\normmm{t-s}}{2ts}+\frac{C\normmm{t-s}}{s}.
 		\end{align*}
 		
 		Arguing in the same way, but inverting the role of $s$ and $t$, we also obtain
 		\begin{align*}
 			f^\Ths_s-f_t^\Th{t}\leq \frac{C\normmm{t-s}}{2ts}+\frac{C\normmm{t-s}}{s},
 		\end{align*}
 		so that one can write
 		\begin{align}
 			\normmm{f\tTh-f_s^\Th{s}}\leq \left(\frac{C}{2ts}+\frac{C}{s}\right)\normmm{t-s},\quad \forall s,t\in(0,h],
 			\label{Lipschitz}
 		\end{align}
 		entailing that $f\tTh$ is locally Lipschitz on $(0,h]$. Therefore by the well known Rademacher's Theorem, we deduce that $f\tTh$ is almost everywhere differentiable on $(0,h]$.
 		\item Show that some left/right derivative in time of $f\tTh$ is controlled by the desired term.
 		
 		Indeed, from \eqref{est}, by assuming $0<s<t\leq h$, and dividing by $t-s$, we obtain 
 		\begin{align}
 			& \nonumber \frac{f\tTh-f_s^\Th{s}}{t-s}\\&\nonumber
 			\leq -\frac{1}{2ts}\normh{\chi_s^\Ths-\xhi_0\circ X_{-T_h(t)}^{\bu_0^k}}^2\nonumber\\&\quad+\frac1{2s}\frac{\normh{\chi_s^\Ths-\xhi_0\circ X_{-T_h(t)}^{\bu_0^k}}^2-\normh{\chi_s^\Ths-\xhi_0\circ X_{-T_h(s)}^{\bu_0^k}}^2}{t-s}.\label{est1}
 		\end{align}
 		Since, recalling \eqref{derivativeA}, it holds
 		\begin{align*}
 			&\lim_{t\to s^+}\frac12\frac{\normh{\chi_s^\Ths-\xhi_0\circ X_{-T_h(t)}^{\bu_0^k}}^2-\normh{\chi_s^\Ths-\xhi_0\circ X_{-T_h(s)}^{\bu_0^k}}^2}{t-s}\\&
 			=\frac12 \left(\partial_r\normh{\chi_s^\Ths-\xhi_0\circ X_{-T_h(r)}^{\bu_0^k}}^2\right)_{r=s}\\&
 			=-\left( \chi_0 \bu_0^k, \nabla \left((-\Delta_N)^{-1}\left(\chi_s^\Ths-\xhi_0\circ X_{-T_h(s)}^{\bu_0^k}\right)\circ X_{T_h(s)}^{\bu_0^k}\right)\right)_\om,
 		\end{align*}
 		and since $f\tTh$ is differentiable almost everywhere on $(0,h]$, by letting $t\to s^+$ we immediately infer that, for almost any $s\in(0,h]$, 
 		\begin{align}
 			\nonumber & \partial_s f_s^\Th{s}\leq  -\frac{1}{2s^2}\normh{\chi_s^\Ths-\xhi_0\circ X_{-T_h(s)}^{\bu_0^k}}^2\\&-\frac1s\left( \chi_0 \bu_0^k, \nabla \left((-\Delta_N)^{-1}\left(\chi_s^\Ths-\xhi_0\circ X_{-T_h(s)}^{\bu_0^k}\right)\circ X_{T_h(s)}^{\bu_0^k}\right)\right)_\om.\label{upper}
 		\end{align}
 		
 		To obtain the reverse inequality, note that, by applying the optimizing properties \eqref{M2} to $f_s^\Ths$, we obtain
 		\begin{align}
 			& \nonumber f\tTh-f_s^\Th{s}\\&\nonumber
            \geq E[\xhi_t^\Tht]+\frac{1}{2t}\normh{\chi_t^\Tht-\xhi_0\circ X_{-T_h(t)}^{\bu_0^k}}^2\\&\nonumber\quad -E[\xhi_t^\Tht]-\frac{1}{2s}\normh{\chi_t^\Tht-\xhi_0\circ X_{-T_h(s)}^{\bu_0^k}}^2\\&\nonumber
 			=\left(\frac1{2t}-\frac1{2s}\right)\normh{\chi_t^\Tht-\xhi_0\circ X_{-T_h(s)}^{\bu_0^k}}^2\nonumber\\&\quad+\frac1{2t}\left(\normh{\chi_t^\Tht-\xhi_0\circ X_{-T_h(t)}^{\bu_0^k}}^2-\normh{\chi_t^\Tht-\xhi_0\circ X_{-T_h(s)}^{\bu_0^k}}^2\right),\label{estA1}
 		\end{align}
 		entailing, assuming again $0<s<t\leq h$,
 		\begin{align*}
 			& \nonumber \frac{f\tTh-f_s^\Th{s}}{t-s}\\&\geq -\frac{1}{2ts}\normh{\chi_t^\Tht-\xhi_0\circ X_{-T_h(s)}^{\bu_0^k}}^2\nonumber\\&\quad+\frac{1}{2t}\frac{\normh{\chi_t^\Tht-\xhi_0\circ X_{-T_h(t)}^{\bu_0^k}}^2-\normh{\chi_t^\Tht-\xhi_0\circ X_{-T_h(s)}^{\bu_0^k}}^2}{t-s},
 		\end{align*}
 		so that, again, by letting this time $s\to t^-$, we infer, for almost any $t\in(0,h]$, 
 		\begin{align*}
 			& \partial_t f_t^\Th{t}\\&\geq  -\frac{1}{2t^2}\normh{\chi_t^\Tht-\xhi_0\circ X_{T_h(t)}^{\bu_0^k}}^2\nonumber
            \\&\quad -\frac1t\left( \chi_0 \bu_0^k, \nabla \left((-\Delta_N)^{-1}\left(\chi_t^\Tht-\xhi_0\circ X_{-T_h(t)}^{\bu_0^k}\right)\circ X_{T_h(t)}^{\bu_0^k}\right)\right)_\om,
 		\end{align*}
 		which together with \eqref{upper} entails \eqref{derivative}.    
 	\end{itemize}
 \end{proof}
 \subsubsection{Energy estimates for the advective Mullins-Sekerka flow
 }
 First note that, to obtain \eqref{derivative}, \textit{any }selection of the function $\xhi_t^{T_h(t)}$ to the minimizing problem \eqref{M2} is fine. Then from \eqref{derivative} we infer that the right-hand side of the equality is Bochner-measurable in time since the left-hand side is the pointwise a.e. derivative of a locally Lipschitz function on $(0,h]$. In order to integtate \eqref{derivative} and decompose the integral of the right-hand side in the sum of two integrals, we need to ensure that both the two summands are measurable. To this aim, in Appendix \ref{measurability} it is proven that we can choose $\xhi_t^{T_h(t)}$ to be (strongly) Bochner-measurable on $(0,h]$ with values in $L^p(\Omega)$ for any $p\geq1$.
 By definition of Bochner measurability, it is then easy to deduce, thanks to the regularity of the flow map $X_t^{\bu_0^k}$, that  $\chi\tTh-\chi_0\circ X_{-T_h(t)}^{\bu_0^k}$ is Bochner-measurable with values in $H^{-1}_0(\Omega)$, as well as $\xhi_{0} \bu_0^k\cdot \nabla \left((-\Delta_N)^{-1}(\chi\tTh-\xhi_0\circ X_{-T_h(t)}^{\bu_0^k})\circ X_{T_h(t)}^{\bu_0^k}\right)$ is Bochner-measurable with values in $L^1(\Omega)$. This is enough to ensure that the two summands of the right-hand side of \eqref{derivative} are measurable and thus it is rigorous to integrate each of them separately.
 We can now safely integrate \eqref{derivative} over $(s,t)$, $0<s<t<h$, to obtain
 \begin{align}
 	\nonumber & f\tTh+\int_s^t\frac{1}{2\tau^2}\normh{\chi_\tau^{T_h(\tau)}-\chi_0\circ X_{-T_h(\tau)}^{\bu_0^k}}^2\d \tau\\& =f_s^\Th{s}-\int_s^t\frac{1}{\tau}\left(\xhi_{0} \bu_0^k\cdot \nabla \left((-\Delta_N)^{-1}(\chi\tTh-\xhi_0\circ X_{-T_h(t)}^{\bu_0^k})\circ X_{T_h(t)}^{\bu_0^k}\right)\right)_\om\d \tau.\label{derivative2}
 \end{align}
 Now we observe that, by the minimizing properties of \eqref{M2}, 
 \begin{align*}
 	f_s^\Th{s}\leq E[\xhi_0]+\frac{s}{2}\normh{\frac{\xhi_0-\xhi_0\circ X_{-\Th{s}}^{\bu_0^k}}{s}}^2,
 \end{align*}
  Plugging this into \eqref{derivative2}, we end up with 
 \begin{align}
 	& \nonumber f\tTh+\int_s^t\frac{1}{2\tau^2}\normh{\chi_\tau^{T_h(\tau)}-\chi_0\circ X_{-T_h(\tau)}^{\bu_0^k}}^2\d \tau\\& \leq E[\xhi_0]+\frac{s}{2}\normh{\frac{\xhi_0-\xhi_0\circ X_{-\Th{s}}^{\bu_0^k}}{s}}^2\nonumber\\&\quad -\int_s^t\frac{1}{\tau}\left( \xhi_{0} \bu_0^k\cdot \nabla \left((-\Delta_N)^{-1}(\chi^\Th{\tau}_\tau-\xhi_0\circ X_{-T_h(\tau)}^{\bu_0^k})\circ X_{T_h(\tau)}^{\bu_0^k}\right)\right)_\om\d \tau.\label{derivative3}
 \end{align}
 
Observe that, crucially, we have by Lemma \ref{basic1} 
 \begin{align*}
 	\frac{ \xhi_0-\xhi_0\circ X_{-\Th{s}}^{\bu_0^k}}{s}\to \bu_0^k\cdot \nabla \xhi_0,\quad \text{ in }\Hmz,\text{ as }s\to 0,
 \end{align*}
 so that 
 \begin{align}
 	\frac{s}{2}\normh{\frac{\xhi_0-\xhi_0\circ X^{\bu_0^k}_{-\Th{s}}}{s}}^2\to 0,\quad\text{ as }s\to 0.\label{limitA}
 \end{align}

 We need to pass from $f\tTh$ to an evaluation in $t=h$ of the energy, but since $f\tTh$ is Lipschitz continuous on $(0,h]$, we can evaluate $f_t^\Th{t}$ at $t=h$, and also 
 \begin{align}
 	\lim_{t\to h^-} f\tTh=f^{\Th{h}}_h.\label{limitB}
 \end{align}
 At this point it is enough let $t\to h^-$ and $s\to 0^+$, recalling \eqref{limitA} and \eqref{limitB}, to obtain the desired inequality
 \begin{align}
 	\nonumber& E[\xhi_h^{T_h(h)}]+\frac1{2h}\normh{\xhi_h^{T_h(h)}-\xhi_0\circ X_{-\Th{h}}^{\bu_0^k}}^2+\int_0^h\frac{1}{2\tau^2}\normh{\chi_\tau^{T_h(\tau)}-\chi_0\circ X_{-T_h(\tau)}^{\bu_0^k}}^2\d \tau\\& \leq E[\xhi_0]-\int_0^h\frac{1}{\tau}\left( \xhi_{0} \bu_0^k\cdot \nabla \left((-\Delta_N)^{-1}(\chi^\Th{\tau}_\tau-\xhi_0\circ X_{-T_h(\tau)}^{\bu_0^k})\circ X_{T_h(\tau)}^{\bu_0^k}\right)\right)_\om\d \tau.\label{derivative3b}
 \end{align}
 Note that, since $\bu_0^k\in \bL^2_\sigma(\Omega)$, we have, recalling \eqref{Nm1},
 \begin{align*}
 	&\normmm{\int_0^h\frac{1}{\tau}\left( \xhi_{0} \bu_0^k\cdot \nabla \left((-\Delta_N)^{-1}(\chi^\Th{\tau}_\tau-\xhi_0\circ X_{-T_h(\tau)}^{\bu_0^k})\circ X_{T_h(\tau)}^{\bu_0^k}\right)\right)_\om\d \tau}\\&
 	\leq C(h,\norm{\bu_0^k}_{\mathbf C^1(\overline\Omega)})\norm{\xhi_0}_{L^\infty(\Omega)}\norm{\bu_0^k}\int_0^h\frac1\tau\normh{\chi_\tau^{T_h(\tau)}-\chi_0\circ X_{-T_h(\tau)}^{\bu_0^k}}\d\tau
 	\\&\leq C(h,\norm{\bu_0^k}_{\mathbf C^1(\overline\Omega)})h+\int_0^h\frac1{4\tau^2}\normh{\chi_\tau^{T_h(\tau)}-\chi_0\circ X_{-T_h(\tau)}^{\bu_0^k}}^2\d \tau.
 \end{align*}
 
 Thus, recalling the definition of $\bu_n^k$ and \eqref{basic}, 
 \begin{align}
 	\nonumber& E[\xhi_h^{T_h(h)}]+\frac1{2h}\normh{\xhi_h^{T_h(h)}-\xhi_0\circ X_{-\Th{h}}^{\bu_0^k}}^2+\int_0^h\frac{1}{4\tau^2}\normh{\chi_\tau^{T_h(\tau)}-\chi_0\circ X_{-T_h(\tau)}^{\bu_0^k}}^2\d \tau\\& \leq E[\xhi_0]+C(h,\norm{\bu_0^k}_{\mathbf C^1(\overline\Omega)})h\nonumber\\&
 	\leq  E[\xhi_0]+C(h).\label{step0}
 \end{align}
 Observe also that, thanks to \eqref{H1c}, it holds, recalling again \eqref{basic},
 \begin{align*}
 	&\normh{\xhi_h^{T_h(h)}-\xhi_0}^2\\&	\leq \left(\normh{\xhi_h^{T_h(h)}-\xhi_0\circ X_{-\Th{h}}^{\bu_0^k}}+\normh{\xhi_0-\xhi_0\circ X_{-\Th{h}}^{\bu_0^k}}\right)^2\\&
 	\leq \left(\normh{\xhi_h^{T_h(h)}-\xhi_0\circ X_{-\Th{h}}^{\bu_0^k}}+h\norm{\bu_0^k}_{\mathbf C(\overline\Omega)}\right)^2,
 \end{align*}
 entailing 
 \begin{align}
 	\frac{h}{2}\normh{\frac{\xhi_h^{T_h(h)}-\xhi_0}{h}}^2\leq \frac1{2h}\left(\normh{\xhi_h^{T_h(h)}-\xhi_0\circ X_{-\Th{h}}^{\bu_0^k}}+C(h)\right)^2.
 	\label{convergenza}	
 \end{align}\label{s1A}
 \subsubsection{Approximating scheme for the fluid equation: interval $[0,h]$}
We construct the following two potentials:
 \begin{align*}
 	&u_h^1:=-(-\Delta_N)^{-1}\frac{\xhi_1^{h}-\xhi_0\circ X_{-h}^{\bu_0^k}}{h},\\&
 	w_h(t):=-(-\Delta_N)^{-1}\frac{\chi_t^{T_h(t)}-\chi_0\circ X_{-T_h(t)}^{\bu_0^k}}{t},\quad t\in(0,h),
 \end{align*}
 We also denote by $\rho_h^1:=\rho_1\xhi^{h}_1+\rho_2(1-\xhi_1^{h})$, $\rho_h^0=\rho_1\xhi_{0}+\rho_2(1-\xhi_{0})$, so that, from the definition of $u_h^1$, it holds (in a suitable weak sense) 
 \begin{align}
 	{\partial}^{\bullet,h}_t\rho^1_h :=\frac{\rho_h^1-\rho_0\circ X_{-h}^{\bu_0^k}}{h} ={(\rho_1-\rho_2)}\Delta u_h^1.\label{rho}
 \end{align}
 
 Now, to write the equation for $\bv_h$, we need the discretized and approximated version of the (formal) equation
 \begin{align*}
 	\partial^\bullet_t \rho \bv +\partial_t^\bullet \bv \rho-\Div(\bv\otimes (\rho_1-\rho_2)\nabla u)-\Div(\nu(\xhi)D\bv)=-\nabla w \xhi,
 \end{align*}
 where $\partial^\bullet_t=\partial_t +\bv\cdot\nabla $ is the total material derivative (in Eulerian coordinates). For the sake of brevity, here and in the following we denote by $\rho$ the composition $\rho(\xhi)$. Note that the discrete operator $\partial_t^{\bullet,h}$ is a discretized version of the material derivative $\partial^{\bullet}$, since
 $$
 \partial_t^{\bullet,h}v_1=\frac{v_1-v_0\circ X_{-h}^{\bu_0^k}}{h}=\frac{v_1\circ X_{h}^{\bu_0^k}-v_0}{h}\circ  X_{-h}^{\bu_0^k}.
 $$
 Also, we define the notion of approximated time derivative as
 $$
 \partial_t^{h}\rho^1_h=\frac{\rho^1_h-\rho_0}{h},
 $$
 with the corresponding affine interpolant
 \begin{align*}
 	\rho_h(t):=\frac{t}{h}\rho_h^1+\frac{h-t}{h}\rho_0,\quad t\in[0,h],
 \end{align*}
 so that $\pt^h\rho_h^1=\pt \rho_h(t)$, $t\in[0,h]$.
 Therefore, we look for $\bv_h$ on $[0,h]$ satisfying, in a suitable weak formulation,  \begin{align}
 	&\nonumber\beta \bA^2\bv_h+\frac12\bP_\sigma\left(\pt \rho_{h}\bv_h\right)+\onehalf\bP_\sigma\left(\bv_h\partial_t^{\bullet,h}\rho^1_h\right)+\onehalf\bP_\sigma (\bv_h\cdot \nabla \rho_h)\bv_h\\&+\bP_\sigma({\rho}_h\pt \bv_h+{\rho}_h(\bv_h\cdot \nabla )\bv_h)-\bP_\sigma\Div(\bv_h\otimes ({\rho_1-\rho_2})\nabla u_h^1) -\bP_\sigma\Div(\nu(\xhi_0)D \bv_h)
 	\nonumber\\&=-\bP_\sigma(\nabla w_h(t)\xhi_0),\\&
 	\Div\bv_h=0,\\&
 	\bv_h(0)=\bv_0,
 	\label{NS1}
 \end{align}
 for $\beta:=\beta(k)=\frac1k>0$. This means that, for almost any $t\in(0,h)$, we have 
 \begin{align}
 	&\nonumber \langle \pt \bP_\sigma(\rho_h\bv_h),\bxi\rangle_{\Hds',\Hds}+\int_\Omega \beta\bA\bv_h\cdot \bA\bxi\dx-\frac12\langle\pt \rho_h,\bv_h\cdot 	\bxi\rangle_{H^1(\om)',H^1(\om)}\\&+\frac12\langle \partial_t^{\bullet,h}\rho^1_h,\bv_h\cdot 	\bxi\rangle_{H^1(\om)',H^1(\om)}\nonumber-\onehalf\int_\Omega \rho_h(\bxi\cdot \nabla)\bv_h\cdot\bv_h\dx-\onehalf\int_\Omega \rho_h(\bv_h\cdot \nabla)\bxi\cdot\bv_h\dx\\&+\int_\Omega{\rho}_h(\bv_h\cdot \nabla )\bv_h\cdot \bxi\dx+\int_\Omega(\bv_h\otimes ({\rho_1-\rho_2})\nabla u_h^1):\nabla \bxi\dx  +\int_\Omega\nu(\xhi_{0})D \bv_h:D \bxi\dx 
 	\nonumber\\&=-\int_\Omega\chi_0\nabla w_h(t)\cdot \bxi\dx,\quad\label{weakformulation} \forall \bxi\in \Hds,\\&
 	(\rho_h\bv_h)(0)=\rho_0\bv_0,\quad \text{ in }\Hds'.
 \end{align}
 We point out the implicit treatment of the advective terms $(\bv_h\cdot\nabla)\bv_h$. This means that the problem is nonlinear, but it can be easily solved by means, e.g., of a Galerkin scheme or a time discretization scheme, thanks to the regularization terms. All the estimates below have to be performed in the approximated scheme and then obtained for $\bv_h$ passing to the limit in the approximating parameters (notice that strong convergence of the velocity as the approximating parameters go to zero can be obtained as in Section \ref{compact} below). For the sake of brevity we omit the details and we only derive the formal estimates.
 
 Note that $\bv_h$ is defined on $(0,h)$ at a continuous level (differently from the approach in \cite{AbelsRoger, ADG, AGGP}, where a time discretization scheme is used), since the surface tension effect has the structure as in the energy \eqref{derivative3}, namely $w_h$ is defined for almost any $t\in(0,h)$. Also, we need to use material derivatives since the estimates from \eqref{step0} are on the discrete version of the material derivative of $\xhi_s^{T_h(s)}$.
 
 The energy estimate, by testing with $\bv_h$ and using the definition of $u_h^1$, reads
 
 \begin{align}
 	&\nonumber\onehalf\ddt \int_\Omega {{\rho}}_h \normmm{\bv_h(t)}^2\dx+\beta \norm{\bA\bv_h}^2+\frac12\int_\Omega \partial_{t}^{\bullet,h}\rho_h^1\vert \bv_h\vert^2\dx
 	\nonumber+\norm{\sqrt{\nu(\xhi_0)}D\bv_h}^2\nonumber\\&
 	=\frac12\int_\Omega \normmm{\bv_h}^2(\rho_1-\rho_2)\Delta u_h^1 dx
 	-\int_\Omega \nabla w_h(s)\chi_0\cdot \bv_h(s) dx.
 \end{align}
 Also, by means of \eqref{rho}, we deduce 
 \begin{align*}
 	\frac12\int_\Omega (\partial_{t}^{\bullet,h}\rho_h^1)\vert \bv_h\vert^2\dx=\frac12\int_\Omega \normmm{\bv_h}^2(\rho_1-\rho_2)\Delta u_h^1 dx,
 \end{align*}
 so that the final energy estimate becomes
 \begin{align}
 	&\nonumber\onehalf\ddt \int_\Omega {{\rho}}_h \normmm{\bv_h(t)}^2\dx+\beta \norm{\bA\bv_h}^2
 	\nonumber+\norm{\sqrt{\nu(\xhi_0)}D \bv_h}^2\nonumber\\&
 	=
 	-\int_\Omega \nabla w_h(s)\chi_0\cdot \bv_h(s) dx.\label{energy0}
 \end{align}
 Recalling that $\bv_0\in \bL_\sigma^2(\Omega)$, by standard inequalities, we also have
 \begin{align}
 	\nonumber &\onehalf\ddt \int_\Omega{\rho}_h \normmm{\bv_h(t)}^2\dx+\norm{\sqrt{\nu(\xhi_0)}D \bv_h}^2+\beta \norm{\bA\bv_h}^2\nonumber\\&
 	\leq \norm{\nabla w_h}\norm{\xhi_0}_{ L^\infty}\norm{\bv_h} 
 	\leq C\norm{\bv_h}^2+C\norm{{\nabla w_h}}^2,\label{abu}
 \end{align}
 from which we deduce the existence of $\bv_h$ on $(0,h)$ with the desired properties, thanks to the regularity of the initial datum $\bv_0$ and \eqref{step0}. Clearly the solution is globally defined over $(0,h)$ by the above control, even if the problem is nonlinear. Namely, from \eqref{abu}, recalling \eqref{step0} and the definition of $w_h$, we infer by means of Gronwall's Lemma,
 \begin{align}
 	\norm{\bv_h}_{L^\infty(0,h;\bL^2_\sigma(\om))}+\norm{\bv_h}_{L^2(0,h;\Hus)}+\sqrt\beta\norm{\bv_h}_{L^2(0,h;\Hds)}\leq C(h,k,\norm{\bv_0}),\label{energyh}
 \end{align}
 recalling $\rho\in[\rho_*,\rho^*]$ and the fact we set $k=h^{-\frac18}$.
 We can also study the regularity of the time derivative. Indeed, let us consider $\bxi\in \Hds$. Then we have, from \eqref{weakformulation}, recalling \eqref{energyh} as well as the embeddings $\Hds\hookrightarrow \bL^\infty(\om)$, $\Hds\hookrightarrow\bW^{1,6}(\om)$, and $\Hus\hookrightarrow \bL^3(\om)$, together with Gagliardo-Nirenberg's inequalities (here we are considering $d=3$, the two-dimensional case being similar),
 \begin{align}
 	&\nonumber\langle \pt \bP_\sigma(\rho_h\bv_h), \bxi\rangle_{\Hds',\Hds} \\&
 	\leq \nonumber\beta\norm{\bA\bv_h}\norm{\bxi}_{\bH^2(\om)}\\&\nonumber\quad +C\left(\norm{\pt\rho_h}_{\Hmz}+\norm{\partial_t^{\bullet,h}\rho_h^1}_{\Hmz}\right)(\norm{\bv_h}_{\Hus(\om)}\norm{\bxi}_{\bL^\infty(\om)}+\norm{\bv_h}_{\bL^3(\om)}\norm{\bxi}_{\bW^{1,6}(\om)})\\&\nonumber
\quad
    +C\norm{\bxi}_{\bL^\infty(\om)}\norm{\bv_h}_{\Hus}\norm{\bv_h}+C\norm{\bv_h}_{\bL^3(\om)}\norm{\bxi}_{\bW^{1,6}(\om)}\norm{\bv_h}+C\norm{\bv_h}\norm{\nabla\bv_h}\norm{\bxi}_{\bL^\infty(\om)}\\&\nonumber
 \quad 	+\norm{\bv_h}_{\bL^3(\om)}\norm{\nabla u_h^1}\norm{\bxi}_{\bW^{1,6}(\om)}+C\norm{\bv_h}_{\Hus}\norm{\bxi}_{\Hus}+\norm{\chi_0}_{L^\infty(\om)}\norm{\nabla w_h}\norm{\bxi}\\&\nonumber
 	\leq C\left(\beta\norm{\bA\bv_h}+\left(\norm{\pt\rho_h}_{\Hmz}+\norm{\partial_t^{\bullet,h}\rho_h^1}_{\Hmz}\right)\norm{\bv_h}_{\Hus(\om)}+\norm{\bv_h}_{\Hus}\right.\\&\quad \left.+\norm{\bv_h}_{\Hus}^\onehalf\norm{\nabla u_h^1}+\norm{\nabla w_h}\right)\norm{\bxi}_{\Hds},\label{dtest}
 \end{align}
 where we also exploited, since $\partial^{\bullet,h}_t \rho_h^1, \pt \rho_h\in \Hmz$,
 \begin{align*}
 	&\langle \partial_t\rho_h, \bv_h\cdot \bxi\rangle_{H^1(\om)',H^1(\om)}=\left(\nabla (-\Delta_N)^{-1}\partial_t\rho_h,\nabla(\bv_h\cdot \bxi)\right)_\om,\\&
 	\langle \partial_t^{\bullet,h}\rho^1_h, \bv_h\cdot \bxi\rangle_{H^1(\om)',H^1(\om)}=\left(\nabla (-\Delta_N)^{-1}\partial_t^{\bullet,h}\rho^1_h,\nabla(\bv_h\cdot \bxi)\right)_\om.
 \end{align*}
 This entails from \eqref{step0} and \eqref{energyh}, recalling the definitions of $w_h$ and $u_h^1$, that 
 
 \begin{align}
 	\norm{\pt \bP_\sigma(\rho_h\bv_h)}_{L^\frac43(0,h;\Hds')}\leq C(h),\label{r1}
 \end{align}
 recalling that we chose both $k$ and $\beta$ as functions of $h$.
 Also, from \eqref{energyh} we infer
 \begin{align}
 	\norm{\bP_\sigma (\rho_h\bv_h)}_{L^2(0,h;\bL^2_\sigma(\om))}\leq C(h).\label{r2}
 \end{align}
 
 From \eqref{abu} and \eqref{r1} we infer that there exists $\bv_h:\Omega\times(0,h)\to \R^d$ solution to \eqref{weakformulation} and such that 
 \begin{align}
 	\nonumber&	\norm{\bv_h}_{L^\infty(0,h;\bL^2_\sigma(\Omega))}+\norm{\bv_h}_{L^2(0,h;\Hus)}+\sqrt\beta\norm{\bv}_{L^2(0,h;\Hds)}	\\&+\norm{\bP_\sigma (\rho_h\bv_h)}_{L^2(0,h;\bL_\sigma^2(\om))}+	\norm{\pt \bP_\sigma(\rho_h\bv_h)}_{L^\frac43(0,h;\Hds')}\leq C(h,\norm{\bv_0}).\label{regv1}
 \end{align}
 We can also easily prove that such a velocity is unique, indeed, let us consider $\bv_h^A,\bv_h^B$ two solutions to \eqref{weakformulation} on $(0,h)$ with the same initial datum $\bv_0$, and the regularity \eqref{regv1}. Then taking the difference and testing by $\bxi=\bv_h:=\bv_h^A-\bv_h^B$, we get (to be rigorous, this estimate should be integrated also in time)
 \begin{align}
 	&\nonumber\onehalf\ddt \int_\Omega {{\rho}}_h \normmm{\bv_h(t)}^2\dx+\beta \norm{\bA\bv_h}^2
 	\nonumber+\norm{\sqrt{\nu(\xhi_0)}D \bv_h}^2\nonumber\\&
 	\nonumber-\int_\om \rho_h(\bv_h\cdot \nabla) \bv_h\cdot \bv_h^A\dx+\onehalf \int_\om \rho_h(\bv_h\cdot\nabla)\bv_h^B\cdot \bv_h\dx -\frac12 \int_\om \rho_h(\bv_h^B\cdot\nabla)\bv_h\cdot\bv_h\dx\\&
 	+\int_\om \rho_h(\bv_h\cdot\nabla)\bv_h^A\cdot\bv_h\dx= 0,\label{energy_0}
 \end{align}
 so that, by H\"{o}lder's and Young's inequalities, together with the embedding $\bH^2(\om)\hookrightarrow\bW^{1,6}(\om)$, we have, recalling $\rho_h\in[\rho_*,\rho^*]$,
 \begin{align}
 	&\nonumber\onehalf\ddt \int_\Omega {{\rho}}_h \normmm{\bv_h(t)}^2\dx+\beta \norm{\bA\bv_h}^2
 	\nonumber+\norm{\sqrt{\nu(\xhi_0)}D \bv_h}^2\nonumber\\&
 	\leq \nonumber C\norm{\bv_h}\norm{\nabla\bv_h}_{\bL^3(	\om)}(\norm{\bv_h^A}_{\bL^6(\om)}+\norm{\bv_h^B}_{\bL^6(\om)})\\&\nonumber\quad+C\norm{\bv_h}\norm{\bv_h}_{\bL^3(	\om)}(\norm{\nabla\bv_h^A}_{\bL^6(\om)}+\norm{\nabla\bv_h^B}_{\bL^6(\om)})\\&
 	\leq C(\beta)(\norm{\bv^A_h}_{\bH^2(\om)}^2+\norm{\bv_h^B}_{\bH^2(\om)}^2)\int_\Omega {{\rho}}_h \normmm{\bv_h(t)}^2\dx+\frac\beta2\norm{\bA\bv_h}^2,\label{energy_1}
 \end{align}
 so that, by Gronwall's Lemma, recalling the regularity \eqref{regv1} for $\bv_h^A,\bv_h^B$, we deduce that $\bv_h^A=\bv_h^B$ almost eveywhere in $\om\times(0,h)$. 
 
 In conclusion, observe that the rigorous energy estimate, which can be obtained by passing to the limit in a suitable approximating scheme, reads: 
 \begin{align}
 	&\onehalf\int_\Omega {\rho}_h(t) \normmm{\bv_h(t)}^2\dx+\int_s^t\beta \norm{\bA\bv_h}^2\d \tau+\int_s^t\norm{\sqrt{\nu(\xhi_0)}D \bv_h}^2\d \tau\nonumber\\&
 	\leq 
 	\onehalf\int_\Omega {\rho_h(s)} \normmm{\bv(s)}^2\dx-\int_s^t\int_\Omega \nabla w_h(\tau)\xhi_0\cdot \bv_h(\tau) dx\d \tau,\label{energy}
 \end{align}
 for almost any $0\leq s\leq t<h$, with $s=0$ included.

\subsection{Iterative scheme at step $n>1$}
Let $\bv_h$ and $\xhi_t^{\Th{t}}$ be defined on $[0,(n-1)h]$, $n\in\N$, $n>1$. Let also $\xhi^h_{m-1}$ be given for any $0\leq m\leq n$. Then compute $\bu^k_h$, defined on $[0,(n-1)h]$, from \eqref{bvk}, and consider the sequence $\{\bu_m^k\}_{m\in\N,\ m\leq {n-1}}$, where we recall that $\bu^k_m:=\bu^k_h(mh)$, for $m\in \N$, $m\leq n-1$.
\subsubsection{Minimizing movements scheme}
We set the following minimizing sequence:
\begin{align}
	&\chi_{n}^{h}\in\arg\min_{\chi\in \mathcal M_{m_0}}\left(E[\chi]+\frac 1{2h}\normh{\chi-\chi_{n-1}^{h}\circ X_{-h}^{\bu_{n-1}^k}}^2\right),\label{M1b}
\end{align}

Then, we define the interpolation sequence as:
\begin{align}
	& \xhi^{T_h((n-1)h)}_{(n-1)h}:=\xhi_{n-1}^h,\\&
	\xhi_t^{T_h(t)}\in\arg\min_{\chi\in \mathcal M_{m_0}}\left(E[\chi]+\frac 1{2(t-(n-1)h)}\normh{\chi- \xhi^{h}_{n-1}\circ X_{-T_h(t)}^{\bu_{n-1}^k}}^2\right),\nonumber\\&
	\text{ for }t\in ((n-1)h,nh].\label{M2b}
\end{align}
As already observed, for any fixed $t$, the problem admits at least one solution by the Direct Method of the Calculus of Variations. Observe that all the estimates in Section \ref{s1A} are still valid up to substituting $\xhi_0$ and $\bu_0^k$ with $\xhi_{m}^h$ and $\bu_{m}^k$, respectively, for all $1\leq m\leq n-1$.

%
%
\subsubsection{Approximating scheme for the fluid equation: interval $[0,nh]$}
We define, for any $1\leq m\leq n$, 
\begin{align}
&\label{uh1}u_h^m:=-(-\Delta_N)^{-1}\frac1{h}\left(\xhi_m^{h}-\xhi_{m-1}^h\circ X_{-h}^{\bu_{m-1}^k}\right),\\&
w_h(t):=-(-\Delta_N)^{-1}\frac{1}{\Th{t}}\left(\chi_{t}^{T_h(t)}-\xhi_{m-1}^h\circ X_{-T_h(t)}^{\bu_{m-1}^k}\right),\quad t\in((m-1)h,mh).\label{wh}
\end{align}
Denote by $\rho_h^m:=\rho_1\xhi^{h}_{m}+\rho_2(1-\xhi_{m}^{h})$, so that
\begin{align}
    {\partial}^{\bullet,h}_t\rho^m_h=\frac{\rho_h^m-\rho^{m-1}_h\circ X_{-h}^{\bu_{m-1}^k}}{h}={(\rho_1-\rho_2)}\Delta u_h^m.\label{rho1}
\end{align}
 Also, we extend the affine interpolant $\rho_h$ on the interval $[0,nh]$:
\begin{align*}
	\rho_h(t):=\frac{t-(m-1)h}{h}\rho_h^m+\rho_h^{m-1}\frac{mh-t}{h},\quad t\in[(m-1)h,mh],\quad 1\leq m\leq n.
\end{align*}
We then introduce the piecewise constant interpolants as 
\begin{align*}
	\xhi^h(t):=\xhi_{m-1}^h,\quad  u^h(t):=u_h^{m-1}\quad t\in[(m-1)h,mh),\quad 1\leq m\leq n,
\end{align*}
together with the affine interpolant

\begin{align*}
\widehat\xhi_h(t):=\frac{t-(m-1)h}{h}\xhi_m^h+\chi^h_{m-1}\frac{mh-t}{h},\quad t\in[(m-1)h,m+h], \quad 1\leq m\leq n.
\end{align*}
Also, we define the piecewise constant interpolant of $\bu_n^k$ as
$$
\widehat{\bu}_h^k(t):=\bu_h^k((m-1)h) =\bu_k^{m-1},\quad t\in[(m-1)h,mh), \quad 1\leq m\leq n.
$$
Therefore, the approximated total derivatives read
\begin{align*}
	\partial^{\bullet,h}_t{\xhi}^h(t):=\frac {\xhi_{m}^{h}-\xhi_{m-1}^h\circ X_{-h}^{\widehat \bu_h^k(t)}}{h},\quad t\in[(m-1)h,mh), \quad 1\leq m\leq n.
\end{align*}
and
$$
\partial_t^{\bullet,h} \rho_h:=\frac{\rho_m^h-\rho_{m-1}^h\circ X_{-h}^{\widehat \bu_h^k(t)}}{h},\quad t\in[(m-1)h,mh), \quad 1\leq m\leq n.
$$ 
Then, we extend $\bv_h$ on $(0,nh)$, satisfying, for almost any $0<t\leq  nh$, $n\in \N$, 
\begin{align}
\nonumber	&\nonumber \langle \pt(\rho_h\bv_h),\bxi\rangle_{\Hds',\Hds}+\int_\Omega \beta\bA\bv_h\cdot \bA\bxi\dx+\frac12\langle\left(\partial_t^{\bullet,h}\rho_h-\pt \rho_h\right),\bv_h\cdot\bxi\rangle_{H^1(\om)',H^1(\om)}\\&\nonumber-\onehalf\int_\Omega \rho_h(\bxi\cdot \nabla)\bv_{h}\cdot\bv_h\dx-\onehalf\int_\Omega \rho_h(\bv_h\cdot \nabla)\bxi\cdot\bv_h\dx\\&\nonumber+\int_\Omega{\rho}_h(\bv_h\cdot \nabla )\bv_{h}\cdot \bxi\dx+\int_\Omega(\bv_h\otimes ({\rho_1-\rho_2})\nabla u_h):\nabla \bxi\dx  +\int_\Omega\nu(\xhi^h)D \bv_h:D \bxi\dx 
	\nonumber\\&\nonumber=-\int_\Omega\xhi^h\nabla w_h\cdot \bxi\dx,\quad \forall \bxi\in\Hds,\\&
	(\rho_h\bv_h)(0)=\rho_0\bv_0, \quad \text{ in }\Hds'.\label{velocity1}
\end{align}
Note that estimates \eqref{regv1} still hold on $(0,nh)$, so that the problem is well posed and guarantees the existence of a unique $\bv_h$, which, restricted to $(0,(n-1)h)$, coincides with the one computed on $(0,(n-1)h)$ at step $n-1$. Notice that the uniqueness of the solution is essential to guarantee that the extension of $\bv_h$ from $(0,(n-1)h)$ to $(0,nh)$ is uniquely defined. Clearly, also \eqref{energy} still holds, for almost any $0\leq s\leq t<nh$, with $s=0$ included. We can now update $\bu^k_h$ on $(0,nh)$ by solving \eqref{bvk} on $(0,nh)$ (again this solution is the extension of $\bu^k_h$ of step $n-1$, by uniqueness), and set $\xhi^h_{n}:=\xhi_{nh}^{T_h(nh)}$.

\subsection{Global energy estimates}
We now iterate the procedure of the previous section for any $n\in\mathbb N$.

In the following, we will see how to show that the constant in the right-hand side of the estimates in the above sections can be shown to be independent of $h$, so that we will be allowed to pass to the limit as $h\to 0$, recalling that $k=h^{-\frac1{8}},\beta=\tfrac1k=h^\frac1{8}$. This will give us the existence of the desired weak solution to the Navier-Stokes-Mullins-Sekerka system.

Recall that on any interval of the form $((n-1)h,nh)$, $n\in\N$, we have that \eqref{derivative3b} and \eqref{energy} hold, namely we have 
\begin{align}
	&	 E[\xhi_{n}^{h}]+\frac1{2h}\normh{\xhi_{n}^{h}-\xhi_{n-1}^h\circ X_{-h}^{ \bu_{n-1}^k}}^2\nonumber\\&\label{energy1a}+\int_{(n-1)h}^{nh}\frac{1}{2{T_h(\tau)}^2}\normh{\chi_\tau^{T_h(\tau)}-\xhi_{n-1}^h\circ X_{-T_h(\tau)}^{ \bu_{n-1}^k}}^2\d \tau\\& \nonumber\leq E[\xhi_{n-1}^h]-\int_{(n-1)h}^{nh}\frac{1}{T_h(\tau)}\left( \xhi_{n-1}^{h} \bu_{n-1}^k,\nabla \left((-\Delta_N)^{-1}\left(\chi_\tau^\Th{\tau}-\xhi_{n-1}^h\circ X_{-T_h(\tau)}^{ \bu_{n-1}^k}\right)\circ X_{T_h(\tau)}^{\bu^k_{n-1}}\right)\right)_\om\d \tau,
\end{align}
and, for almost any $0\leq s<t<nh$, and any $n\in \mathbb N$,
\begin{align}
	&\nonumber\onehalf\int_\Omega {\rho}_h(t) \normmm{\bv_h(t)}^2\dx+\int_s^t\beta\norm{\bA\bv_h}^2\d\tau+\int_s^t\norm{\sqrt{\nu(\xhi_{n}^{h})}D\bv_h}^2\d\tau\nonumber
	\\&
	\leq \onehalf\int_\Omega {\rho}_h(s) \normmm{\bv_h(s)}^2\dx-\int_s^t\int_\Omega \nabla w_h(\tau)\xhi_{n-1}^h\cdot \bv_h(\tau) dx\d \tau.\label{energyv}
\end{align}
Notice that the case $s=0$ is included (see \eqref{energy}).
We then introduce the function
\begin{align*}
	g_h(t):=t-\floor{\frac{t}h}h,
\end{align*}
so that, from \eqref{energy1a}, which is valid for any $n\in\mathbb N$, we immediately infer, by a telescoping argument, that, for any $n,m\in \mathbb N$, $m<n$, 
\begin{align}
	\nonumber	&E[\xhi_{n}^{h}]+\frac1{2}\int_{mh}^{nh}\normh{	\partial^{\bullet,h}_t{\xhi}^h(t)}^2\dt+\int_{mh}^{nh}\frac{1}{2{g_h(\tau)}^2}\normh{\chi_\tau^{g_h(\tau)}-\xhi^h(\tau)\circ X_{-g_h(\tau)}^{\widehat{\bu}_h^k(\tau)}}^2\d \tau\\& \leq E[\xhi_{m}^h]\nonumber\\&\quad -\int_{mh}^{nh}\frac{1}{g_h(\tau)}\left( \xhi^h{(\tau)} \widehat{\bu}_h^k(\tau),\nabla \left((-\Delta_N)^{-1}\left(\chi_\tau^{g_h(\tau)}-\xhi^h(\tau)\circ X_{-g_h(\tau)}^{\widehat{\bu}_h^k(\tau)}\right)\circ X_{g_h(\tau)}^{\widehat{\bu}_h^k(\tau)}\right)\right)_\om\d \tau.
	\label{energy1ab}
\end{align}
Note now that, by the minimizing properties of $\xhi_t^{T_h(t)}$, we have
\begin{align}
	\nonumber	E[\xhi_t^{T_h(t)}]&\leq \nonumber E[\xhi_t^{T_h(t)}]+\frac1{2T_h(t)}{\normh{\xhi_t^{T_h(t)}-\xhi_n^h\circ X_{-T_h(t)}^{\bu^k_n}}^2}\\&\leq E[\xhi_{n}^h]+\frac{1}{2T_h(t)}\normh{\xhi_n^h-\xhi_n^h \circ X_{-T_h(t)}^{\bu^k_n} }^2,\label{minimizer2}
\end{align}
for any $nh< t <(n+1)h$. Observe also that, due to \eqref{H1c}, recalling \eqref{basic}, we infer
\begin{align*}
	\frac{1}{2T_h(t)}\normh{\xhi_n^{h}-\xhi_n^h\circ X_{-T_h(t)}^{\bu^k_n}}^2\leq C(h)T_h(t),
\end{align*}
so that, since $T_h(t)\to 0$ for $t\to (nh)^+$ ($h>0$ fixed), by a slight abuse of notation, we can set the quantity $	\frac{1}{2T_h(t)}\normh{\xhi_n^{h}-\xhi_n^h\circ X_{-T_h(t)}^{\bu^k_n}}^2$ to be equal to zero when $t=nh$. This means that \eqref{minimizer2} actually holds for any $t\in [nh,(n+1)h)$. Let us now fix $T^*>0$, and $0<s<\kappa<\tau<T<T^*$. We can assume $h<\min\{T-\tau,\kappa-s\}$, so that we can set $n_0:=\floor{\frac{T}h}$ and $m_0:=\ceil{\frac sh}$, entailing
$$
0<s\leq m_0h<\kappa<\tau<n_0h\leq T,\quad T\in[n_0h,(n_0+1)h),\quad s\in((m_0-1)h,m_0h].
$$
Then, from what discussed above, we deduce from \eqref{energy1ab}-\eqref{minimizer2}, with $n=n_0,\ m=0$, 
\begin{align}
	\nonumber	&E[\xhi_{T}^{T_h(T)}]+\frac1{2}\int_{0}^{\tau}\normh{	\partial^{\bullet,h}_t{\xhi}^h(t)}^2\dt+\int_{0}^{\tau}\frac{1}{2{g_h(t)}^2}\normh{\chi_t^{g_h(t)}-\xhi^h(t)\circ X_{-g_h(t)}^{\widehat{\bu}_h^k(t)}}^2\d t\\& \leq E[\xhi_{0}]+\frac{1}{2T_h(T)}\normh{\xhi^h(T)-\xhi^h(T)\circ  X_{-T_h(T)}^{\widehat{\bu}_h^k(T)} }^2\nonumber\\&\quad -\int_{0} ^{n_0h}\frac{1}{g_h(t)}\left( \xhi^h{(t)} \widehat{\bu}_h^k(t),\nabla \left((-\Delta_N)^{-1}(\chi_t^{g_h(t)}-\xhi_{n-1}^h\circ X_{-g_h(t)}^{\widehat{\bu}_h^k(t)})\circ X_{g_h(t)}^{\widehat{\bu}_h^k(t)}\right)\right)_\om\d t.
	\label{energy1abcd}
\end{align}
Now we point out that, by considering the interval $((m_0-1)h,m_0h]$, we can argue as in \eqref{derivative} to deduce that 
\begin{align*}
	E[\xhi_{m_0}^h]&\leq E[\xhi_{m_0}^h]+\frac{1}{2h}\normh{\xhi_{m_0}^{h}-\xhi_{m_0-1}^h\circ X_{-h}^{\bu^{k}_{m_0-1}}}^2\\&\quad +\int_s^{m_0h}\frac{1}{2 g_h(\tau)^2}\normh{\xhi_\tau^{g_h(\tau)}-\xhi_{m_0-1}^h\circ X_{-h}^{\bu^{k}_{m_0-1}}}^2\d\tau\\&=E[\xhi_{s}^{T_h(s)}]+\frac{1}{2T_h(s)}\normh{\xhi_{s}^{T_h(s)}-\xhi_{m_0-1}^h\circ X_{-T_h(s)}^{\bu^{k}_{m_0-1}}}^2\\&\quad -\int_s^{m_0h}\frac1{g_h(\tau)}\left(\xhi^h(\tau)\widehat{\bu}_h^k(\tau)\cdot \nabla\left((-\Delta_N)^{-1}\left(\xhi_\tau^{g_h(\tau)}-\xhi^{h}(\tau)\circ X_{-g_h(\tau)}^{\widehat{\bu}_h^k(\tau)}\right)\right)\right)_\om\d\tau.
\end{align*}

Therefore, recalling \eqref{minimizer2}, and since the integrands in the left-hand side are all positive, we infer from \eqref{energy1ab}, with $n=n_0,\ m=m_0$, 
\begin{align}
	\nonumber	&E[\xhi_{T}^{T_h(T)}]+\frac1{2}\int_{\kappa}^{\tau}\normh{	\partial^{\bullet,h}_t{\xhi}^h(t)}^2\dt+\int_{\kappa}^{\tau}\frac{1}{2{g_h(t)}^2}\normh{\chi_t^{g_h(t)}-\xhi^h(t)\circ X_{-g_h(t)}^{\widehat{\bu}_h^k(t)}}^2\d t\\& \leq E[\xhi_{s}^{T_h(s)}]+\frac1{2T_h(s)}{\normh{\xhi_s^{T_h(s)}-\xhi^h(s)\circ X_{-T_h(s)}^{\widehat{\bu}_h^k(s)}}^2}+\frac{1}{2T_h(T)}\normh{\xhi^h(T)-\xhi^h(T)\circ  X_{-T_h(T)}^{\widehat{\bu}_h^k(T)} }^2\nonumber\\&\quad -\int_{s} ^{\floor{\frac Th}h}\frac{1}{g_h(t)}\left( \xhi^h{(t)} \widehat{\bu}_h^k(t),\nabla \left((-\Delta_N)^{-1}(\chi_t^{g_h(t)}-\xhi_{n-1}^h\circ X_{-g_h(t)}^{\widehat{\bu}_h^k(t)})\circ X_{g_h(t)}^{\widehat{\bu}_h^k(t)}\right)\right)_\om\d t.
	\label{energy1abc}
\end{align}
Following \cite{SH}, we introduce the varifolds associated to the varifold lift of $\xhi^{\Th{t}}_t$ in the interior and
on the boundary by
\begin{align*}
	\mu_t^{h,\Omega}:=\normmm{\nabla\xhi^{T_h(t)}_t}\otimes \left(\delta_{\frac{\nabla \xhi^{\Th{t}}_t}{\normmm{\nabla\xhi^{\Th{t}}_t}}(x)}\right)_{x\in \Omega}\in \mathcal M({\Omega}\times \mathbb S^{d-1}),
\end{align*}
and 
\begin{align*}
	\mu_t^{h,\partial\Omega}:=\xhi_t^{\Th{t}}\llcorner \partial\Omega \otimes \left(\delta_{\bn_{\partial\Omega}(x)}\right)_{x\in \partial\Omega}\in \mathcal M(\partial\Omega\times \mathbb S^{d-1}),
\end{align*}
where $\bn_{\partial\Omega}$ denotes the inner normal on $\partial\Omega$ and where with an abuse of
notation we do not distinguish between a function and its trace along $\partial\Omega$. In conclusion we introduce the total approximated varifold by 
\begin{align}
	\mu_t^h:=c_0\mu_t^{h,\Omega}+c_0\cos\gamma\mu_t^{h,\partial\Omega}\in \mathcal M(\overline{\Omega}\times \mathbb S^{d-1}),\quad t\in (0,T^*),
	\label{totvar}
\end{align} 
so that for any $t\in(0,T^*)$ it holds $E[\xhi_t^{T_h(t)}]=\normmm{\mu_t^h}_{\mathbb S^{d-1}}(\overline\Omega)$. Therefore, from \eqref{energy1abcd}, we have 
\begin{align}
	\nonumber	&\normmm{\mu_T}_{\mathbb S^{d-1}}(\overline\Omega)+\frac1{2}\int_{0}^{\tau}\normh{	\partial^{\bullet,h}_t{\xhi}^h(t)}^2\dt+\frac{1}{2}\int_{0}^{\tau}\norm{\nabla w_h(t)}^2\d t\\& \leq E[\xhi_0]+\frac{1}{2T_h(T)}\normh{\xhi^h(T)-\xhi^h(T)\circ  X_{-T_h(T)}^{\widehat{\bu}_h^k(T)} }^2\nonumber\\&\quad +\int_{0}^{\floor{\frac Th}h}\left( \xhi^h{(t)} \widehat{\bu}_h^k(t),\nabla \left(w_h(t)\circ X_{g_h(t)}^{\widehat{\bu}_h^k(t)}\right)\right)_\om\d t,
	\label{energy1abc2b}
\end{align}
as well as, from \eqref{energy1abc},

\begin{align}
	\nonumber	&\normmm{\mu_T^h}_{\mathbb S^{d-1}}(\overline\Omega)+\frac1{2}\int_{\kappa}^{\tau}\normh{	\partial^{\bullet,h}_t{\xhi}^h(t)}^2\dt+\frac{1}{2}\int_{\kappa}^{\tau}\norm{\nabla w_h(t)}^2\d t\\& \leq \normmm{\mu_s^h}_{\mathbb S^{d-1}}(\overline\Omega)+\frac{T_h(s)}{2}\norm{\nabla w_h(s)}^2+\frac{1}{2T_h(T)}\normh{\xhi^h(T)-\xhi^h(T)\circ  X_{-T_h(T)}^{\widehat{\bu}_h^k(T)} }^2\nonumber\\&\quad +\int_{s} ^{\floor{\frac Th}h}\left( \xhi^h{(t)} \widehat{\bu}_h^k(t),\nabla \left(w_h(t)\circ X_{g_h(t)}^{\widehat{\bu}_h^k(t)}\right)\right)_\om\d t.
	\label{energy1abc3}
\end{align}
Note also that the map $X_{s}^{\widehat{\bu}_h^k(t)}$, $s\in \mathbb R$, by the definitions \eqref{flowmap} and of the piecewise constant $\widehat{\bu}_h^k$, is such that, for any $s\in \mathbb R$ and $t\in[0,T]$,
 \begin{align}
	\begin{cases}\frac{\d}{\d s} X_{s}^{\widehat{\bu}_h^k(t)}=\widehat{\bu}_h^k(X_{s}^{\widehat{\bu}_h^k(t)}, t),\\
		X_0=Id.
	\end{cases}
	\label{map2}
\end{align}
By the minimizing property \eqref{M2} of $\xhi_{t}^{T_h(t)}$, Allard’s first
variation formula \cite{Allard}, and \cite[Lemma 10]{SH}, it follows that $\xhi_t^{T_h(t)}$ also satisfies the approximate Gibbs-Thomson relation:
\begin{align}
	\int_{\overline{\Omega}\times \mathbb S^{d-1}}(Id-s\otimes s):\nabla B(x)\d \mu_t^h(x,s)=\int_\Omega \xhi_t^{T_h(t)}\text{div}(w_h(\cdot,t)B)\d x,\label{allard}
\end{align}
for all $t\in(0,T^*)$ and all $B\in \mathcal S_{\xhi_t^{T_h(t)}}$. 


\subsection{Uniform estimates in $h$}
\subsubsection{Uniform bounds for $\widehat\bu_h^k$}
By the construction of $\bu_h^k$ in \eqref{bvk}, after multiplying the equation by $\bu_h^k$ and integrating over $\om$, we get, by Cauchy-Schwarz and Young's inequalities,
\begin{align*}
	\frac1{2k}\ddt \norm{\bA\bu_h^k}^2+\frac1k\norm{\bA \bu_h^k}^2+\norm{\bu_h^k}^2=(M_k(\bv_h),\bu_h^k)\leq \onehalf\norm{M_k(\bv_h)}^2+\onehalf\norm{\bu_h^k}^2,
\end{align*}
recalling $(\bP_\sigma M_k(\bv_h),\bu_h^k)=(M_k(\bv_h),\bu_h^k)$. Therefore, 
\begin{align*}
	\frac1{2k}\ddt \norm{\bA\bu_h^k}^2+\frac1k\norm{\bA\bu_h^k}^2+\onehalf\norm{\bu^k_h}^2\leq \onehalf\norm{M_k(\bv_h)}^2,
\end{align*}
and thus, integrating over $(0,t)$, $t\in(0,\Ts)$, $\Ts>0$, and recalling that $\bu_h^k(0)=\mathbf 0$, 
\begin{align}
	\frac1{\sqrt k}\norm{\bu_h^k}_{L^\infty(0,t;\Hds)}+\norm{\bu_h^k}_{L^2(0,t;\bL^2_\sigma(\om))}\leq \norm{M_k(\bv_h)}_{L^2(0,t;\bL^2(\om))}\leq   \norm{\bv_h}_{L^2(0,t;\bL^2_\sigma(\om))},
	\label{controls}
\end{align}
where we have used the fact that, from \eqref{cutoff0}, $\normmm{M_k(\bv_h)}\leq \normmm{\bv_h}$. It then holds, for $n=\floor{\tfrac sh}$, recalling the definition of the piecewise constant interpolant $\widehat \bu_h^k(s)=\bu_h^k(\floor{\tfrac s{h}}h)$,
\begin{align}
	\nonumber\norm{\widehat{\bu}_h^k(s)}\nonumber
	&
\leq \norm{\bu_h^k(\floor{\frac s{h}}h)- \bu^k_h(s)}+\norm{\bu_h^k(s)}\\&
	\leq \norm{\bu_h^k(\floor{\frac s{h}}h)- \bu^k_h(s)}+\norm{\bv_h(s)},\label{convvelox1a}
\end{align}
and thus
\begin{align}
	\nonumber&\int_0^t\norm{\widehat{\bu}_h^k(s)}^2\d s\\&
	\leq \int_0^t\norm{\bu_h^k(\floor{\tfrac s{h}}h)- \bu^k_h(s)}^2\d s+2\int_0^t\norm{\bu_h^k(\floor{\tfrac s{h}}h)- \bu^k_h(s)}\norm{\bv_h(s)}\d s+\int_0^t\norm{\bv_h(s)}^2\d s.\label{convvelox1}
\end{align}
Moreover, recalling \eqref{basic}, as we have set $n=\floor{\tfrac sh}$, and thus
$$
\normmm{s-nh}\leq h,
$$
we infer that, for any $s\in[0,\Ts]$, and any $\Ts>0$,
\begin{align}
	&	\norm{\bu_h^k(\floor{\tfrac sh}h)- \bu^k_h(s)}
	\leq \norm{\bu^k_h}_{C^{\frac12}([0,\Ts];\mathbf H^4_\sigma(\om))}	\normmm{s-nh}^\frac12\leq C(\Ts)h^{\frac14}.\label{c1a}
\end{align}
From \eqref{convvelox1} we can then conclude that 
\begin{align}
\int_0^t\norm{\widehat{\bu}_h^k(s)}^2\d s\leq C(\Ts)h^{\frac12}+2C(\Ts)h^{\frac14}\int_0^t\norm{\bv_h(s)}\d s+\int_0^t\norm{\bv_h(s)}^2\d s.
	\label{fund1}
\end{align}

\subsubsection{Uniform bounds for $\mu_h^t$ and $w_h$.}
Recall now that, from the fundamental estimate \eqref{basic}, which is valid for any $\Ts>0$ (setting $b=\Ts$), we infer
\begin{align}
	\norm{\widehat\bu^{k}_h}_{C([0,\Ts];\mathbf C^1(\om))}\leq 	\norm{\bu_h^k}_{C([0,\Ts];\mathbf C^1(\om))}\leq C(\Ts)h^{-\frac14}.
	\label{basica}
\end{align}
By repeating the same arguments as for \eqref{dr}, recalling \eqref{map2} and the crucial estimate \eqref{basica},

\begin{align}
	&\nonumber	\int_0^T\normm{\nabla (w_h\circ X_{g_h(t)}^{\widehat \bu_h^k(t)})}^2\dt \leq {d}(1+h)^2e^{2h\norm{\widehat \bu_h^k}_{C([0,T];\mathbf C^1(\overline{\Omega}))}}\int_0^T\norm{\nabla w_h}^2\dt\nonumber\\&
\leq 	C(\Ts)e^{2C(\Ts)h^{\frac34}}\int_0^T\norm{\nabla w_h}^2\dt\leq C_1(\Ts)\int_0^T\norm{\nabla w_h}^2\dt,
	\label{dr2}
\end{align}
where $C_1(\Ts)>0$ depends on $\Ts$ but not on $h$, since $h^\frac34\to0$ as $h\to0$ and thus we can without loss of generality choose $h\in(0,1)$. 
We can now show that we have uniform estimates in $h$, namely, from \eqref{energy1abc2b}, \eqref{fund1}, \eqref{dr2}, and since $\xhi^h\in\{0,1\}$, by H\"{o}lder's and Young's inequality, 
\begin{align}
\nonumber	&\normmm{\int_{0}^{T}\left( \xhi^h{(t)} \widehat{\bu}_h^k(t),\nabla \left(w_h(t)\circ X_{g_h(t)}^{\widehat{\bu}_h^k(t)}\right)\right)_\om\d t}\nonumber
	\leq \int_0^T\norm{\widehat{\bu}_h^k(t)}\norm{\nabla \left(w_h(t)\circ X_{g_h(t)}^{\widehat{\bu}_h^k(t)}\right)}\dt \\&\nonumber
	\leq \frac{1}{4C_1(\Ts)}\int_0^T\norm{\nabla \left(w_h(t)\circ X_{g_h(t)}^{\widehat{\bu}_h^k(t)}\right)}^2\dt+C(\Ts)\int_0^T\norm{\widehat{\bu}_h^k(t)}^2\dt\\&
	\leq 	\frac14\int_0^T\norm{\nabla w_h}^2\dt+C(\Ts)\int_0^T\norm{\bv_h(t)}^2\dt+C(\Ts)h^\frac14\nonumber\\&
	\leq 	\frac14\int_0^T\norm{\nabla w_h}^2\dt+C(\Ts)\int_0^T\norm{\bv_h(t)}^2\dt+C(\Ts)
	.\label{integrability}
\end{align}
Plugging this in \eqref{energy1abc2b}, we infer 
\begin{align}
	\nonumber	&\normmm{\mu_T^h}_{\mathbb S^{d-1}}(\overline\Omega)+\frac1{2}\int_{0}^{T}\normh{	\partial^{\bullet,h}_t{\xhi}^h(t)}^2\dt+\frac{1}{4}\int_{0}^{T}\norm{\nabla w_h(t)}^2\d t\\& \leq E[\xhi_0]+\frac{1}{2T_h(T)}\normh{\xhi^h(T)-\xhi^h(T)\circ  X_{-T_h(T)}^{\widehat{\bu}_h^k(T)} }^2+C(\Ts)\int_0^T\norm{\bv_h(t)}^2\dt+C(\Ts)\nonumber\\&
	\leq E[\xhi_0]+C(\Ts)\int_0^T\norm{\bv_h(t)}^2\dt+C(\Ts),\quad \forall T<T^*,
	\label{energyfinal0}
\end{align}
for any $T^*>0$. Here we used the fact that, from \eqref{H1c} and \eqref{basic}, recalling $\widehat{\bu}_h^k(T)=\bu^k_{n}$ for $n=\floor{\tfrac Th}$,
\begin{align}
	\nonumber&\frac{1}{2T_h(T)}\normh{\xhi^h(T)-\xhi^h(T)\circ  X_{-T_h(T)}^{\widehat{\bu}_h^k(T)} }^2
	=	\frac{1}{2T_h(T)}\normh{\xhi^h(T)-\xhi^h(T)\circ  X_{-T_h(T)}^{\bu_n^k} }^2
	\nonumber\\&\leq \frac{T_h(T)}{2}\norm{\bu_h^k(nh)}_{\mathbf C(\overline{\Omega})}^2\leq \frac{h}{2}C(\Ts)h^{-\frac12}=C(\Ts)h^\frac12,\label{convh1b}
\end{align}
 entailing that this quantity is bounded uniformly for any $h\in(0,1)$.
 
 \subsubsection{Uniform estimates for the velocity $\bv_h$}
 In order to close the estimate \eqref{energyfinal0} we need an estimate for $\bv_h$. Observe that we infer from \eqref{energyv} that, for almost any $0\leq s<T$, with $s=0$ included,
 \begin{align}
 	&\onehalf\int_\Omega {\rho}_h(T)  \normmm{\bv_h(T)}^2\dx+\beta\int_s^T\norm{\bA\bv_h}^2\dt+\int_s^T\norm{\sqrt{\nu(\xhi^{h}(t))}D \bv_h}^2\dt\nonumber
 	\\&\leq \onehalf\int_\Omega {\rho}_h(s)  \normmm{\bv_h(s)}^2\dx
 	-\int_s^T\int_\Omega \nabla w_h(t)\xhi^h(t)\cdot \bv_h(t) \dx\dt
 	.\label{energyv2}
 \end{align}
 By H\"{o}lder's and Young's inequalities we get
 \begin{align*}
 	&\normmm{\int_0^T\int_\Omega \nabla w_h(t)\xhi^h(t)\cdot \bv_h(t) \dx\dt}
 	\leq C\int_0^T\norm{\bv_h}\norm{\nabla w_h}\d t \\&
 	\leq \frac12\int_0^T\norm{\sqrt{\nu(\chi^h(t))}\bv_h}^2\dt+C_2\int_0^T\norm{\nabla w_h}^2\d t\quad \forall T<T^*,
 \end{align*}
 for some $C_2>0$ and for any $T^*>0$. 
 
 Therefore, from \eqref{energyv2} we deduce, after integrating over $(0,T)$, \begin{align}
 	&\onehalf\int_\Omega {\rho}_h(T)  \normmm{\bv_h(T)}^2\dx+\beta\int_0^T\norm{\bA\bv_h}^2\dt+\frac12\int_0^T\norm{\sqrt{\nu(\xhi^{h}(t))}D \bv_h}^2\dt\nonumber
 	\\&\leq C\norm{\bv_0}^2+C_2\int_0^T\norm{\nabla w_h}^2.
 	\label{energyv2b}
 \end{align}
 Summing up this inequality, multiplied by $\omega=\frac1{8C_2}$, with \eqref{energyfinal0}, we end up with 
 
 \begin{align}
 	\nonumber	&\normmm{\mu_T^h}_{\mathbb S^{d-1}}(\overline\Omega)+\frac\omega2 \int_\Omega {\rho}_h(T)  \normmm{\bv_h(T)}^2\dx
 +\frac1{2}\int_{0}^{T}\normh{	\partial^{\bullet,h}_t{\xhi}^h(t)}^2\dt\\&+\frac{1}{8}\int_{0}^{T}\norm{\nabla w_h(t)}^2\d t+\beta\int_0^T\norm{\bA\bv_h}^2\dt+\onehalf\int_0^T\norm{\sqrt{\nu(\xhi^{h}(t))}D \bv_h}^2\dt \nonumber\\&
 	\leq E[\xhi_0]+\omega C\norm{\bv_0}^2+C(\Ts)\int_0^T\norm{\bv_h(t)}^2\dt+C(\Ts),\quad \forall T<T^*,
 	\label{energyfinal}
 \end{align}
 entailing by the integral Gronwall Lemma, recalling $\rho_h\geq \rho_*>0$, 
 \begin{align}
 	\esssup_{T\leq \Ts}E[\xhi_T^{T_h(T)}]+\esssup_{T\leq \Ts}\frac\omega2 \int_\Omega {\rho}_h(T)  \normmm{\bv_h(T)}^2\dx\leq C(\Ts,\norm{\bv_0},E[\xhi_0]),
 	\label{dissiprel}
 \end{align}
 and thus also
 \begin{align}
 	&	\nonumber\norm{w_h}_{L^2(0,T;\Hz)}+\norm{\partial^{\bullet,h}_t{\xhi}^h(t)}_{L^2(0,T;\Hmz)}\\&+\norm{\bv_h}_{ L^2(0,T;\Hus)}+\sqrt\beta\norm{\bv}_{L^2(0,T;\Hds)}\leq C(\Ts,\norm{\bv_0},E[\xhi_0]),\label{regvbis}
 \end{align}
 for any $T<T^*$ and any $T^*>0$.

Furthermore, an application of Lemma \ref{lemma9}, exploiting the uniform bound on the energy $E$ from \eqref{dissiprel}, gives the existence of $\lambda_h\in L^2(0,T)$ and a constant $C(\Ts,\Omega,d,m_0,\xhi_0,\bv_0)>0$ such that, for all $t\in(0,T^*)$, 
\begin{align}
&	\int_{\overline\Omega\times \mathbb S^{d-1}}(Id-s\otimes s):\nabla B(x)\d \mu_t^h(x,s)=
\int_\Omega \chi_t^{T_h(t)}\text{div}((w_h(\cdot,t)+\lambda_h(t))B)\d x,\label{variation}
\end{align} 
for all $B\in C^1(\overline\Omega;\mathbb R^d)$ with $(B\cdot \bn_{\partial\Omega})|_{\partial\Omega}\equiv 0$, and 
\begin{align}
	\norm{w_h(\cdot,t)+\lambda_h(t)}_{H^1(\Omega)}\leq C(1+\norm{\nabla w_h(\cdot,t)}),
	\label{complete}
\end{align}
for almost any $t\in(0,T^*)$.

Considering the uniform estimates \eqref{regvbis} and \eqref{complete}, we have then shown that it holds
\begin{align}
	\norm{w_h+\lambda_h}_{L^2(0,T^*;H^1(\om))}+\norm{u_h}_{L^2(0,\Ts;\Hz)}\leq C(T^*,\norm{\bv_0},E[\xhi_0]),\label{bounds}
\end{align}
for any $T^*>0$, uniformly in $h$. From now on we omit the dependence on $\norm{\bv_0}$ and $E[\xhi_0]$ of the constants.

In conclusion, note that, by repeating the same arguments as to obtain \eqref{meanv}, we infer 
	\begin{align}
	&\nonumber	\int_0^\Ts\int_\om \frac{\xhi^h(t)-\xhi^h(t)\circ X_{-h}^{\widehat{\bu}_h^k(t)}}{h}\varphi\dx\dt
	\\&=-\int_0^\Ts\int_\om \nabla \varphi\cdot \widehat{\bu}_h^k(t)\int_0^{h} \frac{\xhi^h(t)\circ X_{-\tau}^{\widehat{\bu}_h^k(t)}}{h}\d\tau\dx\dt,\label{meanv1}
\end{align}	
with 
\begin{align}
	\norm{\int_{0}^{h}	\frac{ \xhi^h(t)\circ X_{-\tau}^{\widehat{\bu}_h^k(t)}}{h}\d\tau}_{L^\infty(\om\times(0,\Ts))}\leq 1,
	\label{strconvA}
\end{align} 
and thus 
\begin{align}
	\norm{\frac{\xhi^h(t)-\xhi^h(t)\circ X_{-h}^{\widehat{\bu}_h^k(t)}}{h}}_{L^2(0,\Ts;\Hmz)}\le \norm{\widehat{\bu}_h^k}_{L^2(0,\Ts;\Ls)}\leq C(\Ts),
	\label{bounddiff}
\end{align}
so that 
\begin{align}
	\normh{\pt \wchi^h(t)}\leq \normh{\partial_t^{\bullet,h}\xhi^h(t)}+	\frac{1}{h}\normh{\xhi^h(t)-\xhi^h(t)\circ X_{-h}^{\widehat{\bu}_h^k(t)}},
	\label{controlh}
	\end{align}
entailing from \eqref{regvbis} that 
\begin{align}
		\norm{\pt \wchi^h(t)}_{L^2(0,T;\Hmz)}\leq C(\Ts),
	\label{ptchi}
\end{align}
for any $T<\Ts$ and any $\Ts>0$. 
To conclude this section, we introduce the following function $F_h:(0,\infty)\to \R$:
\begin{align}
\nonumber	F_h(t):=&E[\xhi_t^{Th(t)}]+\frac{1}{2T_h(t)}\normh{\xhi_t^{T_h(t)}-\xhi^h(t)\circ X_{-T_h(t)}^{\widehat{\bu}_h^k(t)}}^2\\&- \int_{0}^{t}\left( \xhi^h{(\tau)} \widehat{\bu}_h^k(\tau),\nabla \left(w_h(\tau)\circ X_{g_h(\tau)}^{\widehat{\bu}_h^k(\tau)}\right)\right)_\om\d \tau,\label{F}
\end{align}
which is the discretized version of the nonincreasing function $F$ defined in the Introduction in \eqref{discretization}.
\subsubsection{Bounds for the velocity time derivative.}

We can now repeat the very same estimates as in \eqref{dtest}, to infer, for almost any $t\in(0,T)$, exploiting \eqref{regvbis}, 

\begin{align}
	&\nonumber\langle \bP_\sigma\pt(\rho_h\bv_h), \bxi\rangle_{\Hds',\Hds} \\&
	\leq \nonumber\beta\norm{\bA\bv_h}\norm{\bxi}_{\bH^2(\om)}\\&\nonumber\quad +C\left(\norm{\pt\rho_h}_{\Hmz}+\norm{\partial_t^{\bullet,h}\rho_h}_{\Hmz}\right)(\norm{\bv_h}_{\Hus(\om)}\norm{\bxi}_{\bL^\infty(\om)}+\norm{\bv_h}_{\bL^3(\om)}\norm{\bxi}_{\bW^{1,6}(\om)})\\&\nonumber
	\quad +C\norm{\bxi}_{\bL^\infty(\om)}\norm{\bv_h}_{\Hus}\norm{\bv_h}+C\norm{\bv_h}_{\bL^3(\om)}\norm{\bxi}_{\bW^{1,6}(\om)}\norm{\bv_h}+C\norm{\bv_h}\norm{\nabla\bv_h}\norm{\bxi}_{\bL^\infty(\om)}\\&\nonumber
	\quad+C\norm{\bv_h}_{\bL^3(\om)}\norm{\nabla u_h}\norm{\bxi}_{\bW^{1,6}(\om)}+C\norm{\bv_h}_{\Hus}\norm{\bxi}_{\Hus}+\norm{\xhi^h}_{L^\infty(\om)}\norm{\nabla w_h}\norm{\bxi}\\&\nonumber
	\leq C\left(\beta\norm{\bA\bv_h}+\left(\norm{\pt\rho_h}_{\Hmz}+\norm{\partial_t^{\bullet,h}\rho_h}_{\Hmz}\right)\norm{\bv_h}_{\Hus}+\norm{\bv_h}_{\Hus}\right.\\&\quad \left.+\norm{\bv_h}_{\Hus}^\onehalf\norm{\nabla u_h}+\norm{\nabla w_h}\right)\norm{\bxi}_{\Hds},\label{dtest2}
\end{align}
entailing, by means of \eqref{regvbis}, \eqref{bounds}, and \eqref{ptchi}, that 

\begin{align}
	&	\norm{\bP_\sigma (\rho_h\bv_h)}_{L^2(0,T;\bL^2_\sigma(\om))}+	\norm{\pt \bP_\sigma(\rho_h\bv_h)}_{L^\frac43(0,T;\Hds')}\leq C(\Ts),\label{regv1b}
\end{align}
uniformly in $h$, for any $T<\Ts$ and any $\Ts>0$.

\subsection{Monotonicity and boundedness of $F_h$.} We aim at showing that $F_h$, analogously to its limit counterpart $F$ defined in \eqref{discretization}, is nonincreasing in time. First note that, on each interval $((n-1)h,nh]$, it holds from \eqref{derivative} (after integration over $(s,t)$, and recalling the definition of $w_h$)
\begin{align*}
	F_h(t)\leq F_h(s),\quad \forall \ (s,t):\ (n-1)h<s\leq t\leq nh,\quad \forall n\in \mathbb N.
\end{align*}
Also, from the continuity on $((n-1)h,nh]$ of the first two summands of $F_h$, together with the integrability in time of the integrand in the third summand, given by \eqref{integrability}, we deduce that $F_h$ is continuous on any $((n-1)h,nh]$ interval, and thus the limit of $F_h(t)$ as $t\to (nh)^+$ exists for any $n\in \mathbb N$. We denote this limit by $F_h((nh)^+)$, whereas the corresponding limit from the left is simply denoted by $F_h(nh)$, since the function is left continuous in $nh$ for any $n\in \mathbb N$. In order to show that $F_h$ is monotone nonincreasing, it is then enough to ensure 
\begin{align}
 F_h((nh)^+)\leq F_h(nh),\quad \forall n\in \mathbb N.
	\label{monotone}
\end{align} 
Now, by the minimizing properties \eqref{M2} of $\xhi_t^{T_h(t)}$, it holds, for any fixed $n\in \mathbb N$ 
\begin{align*}
	&\nonumber E[\xhi_t^{\Th{t}}]+\frac{1}{2T_h(t)}\normh{\xhi_t^{T_h(t)}-\xhi^h(t)\circ X_{-T_h(t)}^{\widehat{\bu}_h^k(t)}}^2\\&\leq E[\xhi_n^h]+\frac{1}{2T_h(t)}\normh{\xhi^h_n-\xhi^h_n\circ X_{-T_h(t)}^{\bu^k_{n}}}^2, \quad \forall t\in (nh,(n+1)h],
\end{align*}
and thus, summing in both sides $- \int_{0}^{t}\left( \xhi^h{(\tau)} \widehat{\bu}_h^k(\tau),\nabla \left(w_h(\tau)\circ X_{g_h(\tau)}^{\widehat{\bu}_h^k(\tau)}\right)\right)_\om\d \tau$, we infer 
\begin{align}
		\nonumber &F_h(t)\leq F_h(t)+\tfrac{1}{2h}\normh{\xhi_n^h-\xhi_{n-1}^h\circ X_{-h}^{\bu_{n}^k}}^2\\&\nonumber\leq E[\xhi_n^h]+\tfrac{1}{2h}\normh{\xhi_n^h-\xhi_{n-1}^h\circ X^{\bu_{n}^k}_{-h}}^2\nonumber\\&\nonumber\quad - \int_{0}^{t}\left( \xhi^h{(\tau)} \widehat{\bu}_h^k(\tau),\nabla \left(w_h(\tau)\circ X_{g_h(\tau)}^{\widehat{\bu}_h^k(\tau)}\right)\right)_\om\d \tau+\frac{1}{2T_h(t)}\normh{\xhi^h_n-\xhi^h_n\circ X_{-T_h(t)}^{\bu^k_{n}}}^2\\&
		\leq F_h(nh)+\int_t^{nh}\left( \xhi^h{(\tau)} \widehat{\bu}_h^k(\tau),\nabla \left(w_h(\tau)\circ X_{g_h(\tau)}^{\widehat{\bu}_h^k(\tau)}\right)\right)_\om\d \tau\nonumber\\&\quad +\frac{1}{2T_h(t)}\normh{\xhi^h_n-\xhi^h_n\circ X_{-T_h(t)}^{\bu^k_{n}}}^2,
	\label{tolim}
\end{align}
for any $ t\in (nh,(n+1)h]$. Also, recalling \eqref{H1c}, \eqref{convh1b}, and the fact that $T_h(t)\to0$ as $t\to (nh)^+$, it holds (with $h>0$  fixed)
\begin{align*}
\frac{1}{2T_h(t)}\normh{\xhi^h_n-\xhi^h_n\circ X_{T_h(t)}^{\bu^k_{n}}}^2\leq T_h(t)	\norm{\bu_n^k}^2_{\mathbf C(\overline\Omega)}\leq C(\Ts)T_h(t)h^{-\frac12}\to 0,\quad \text{ as }t\to (nh)^+,
\end{align*}
as well as, since the integrand is integrable in time, 
$$
\int_t^{nh}\left( \xhi^h{(\tau)} \widehat{\bu}_h^k(\tau),\nabla \left(w_h(\tau)\circ X_{g_h(\tau)}^{\widehat{\bu}_h^k(\tau)}\right)\right)_\om\d \tau\to 0,\quad \text{as }t\to (nh)^+.
$$
Therefore, passing to the limit as $t\to (nh)^+$ in 
\eqref{tolim}, we infer \eqref{monotone}, from which we conclude that 
\begin{align}
	F_h(t)\leq F_h(s),\quad \forall 0<s\leq t.
	\label{monoton}
\end{align}
We can also prove that  $F_h(t)$ is uniformly bounded on $(0,+\infty)$. Indeed, by the monotonicity \eqref{monoton} we immediately see that 
\begin{align}
F_h(t)\leq F_h(h)\leq E[\xhi_0],\quad \forall t\geq h,
\label{fq}
\end{align}
where we used \eqref{derivative3b}. Furthermore, from \eqref{derivative3}, by letting $s\to0^+$, recalling \eqref{limitA}, we infer 
\begin{align*}
F_h(t)\leq E[\xhi_0],\quad \forall t\in(0,h],	
\end{align*}
which entails that we can also extend $F_h$ at $t=0$ by continuity. Therefore, this bound together with \eqref{fq}, allows us to conclude
that 
\begin{align}
	\sup_{t\geq 0}F_h(t)\leq E[\xhi_0].
	\label{unifcontr}
\end{align}
We also need a lower bound on $F_h$, which is easily obtained from \eqref{integrability} and \eqref{regvbis}. Indeed, we infer that, for any $T^*>0$,
\begin{align}
	\nonumber F_h(t)&\geq -\normmm{\int_0^t\left( \xhi^h{(\tau)} \widehat{\bu}_h^k(\tau),\nabla \left(w_h(\tau)\circ X_{g_h(\tau)}^{\widehat{\bu}_h^k(\tau)}\right)\right)_\om\d \tau}\\&\geq -C(T^*,E[\xhi_0],\norm{\bv_0}),\quad \forall t\in[0,T^*],
	\label{lowerbound}
\end{align}
so that, plugging this together with \eqref{unifcontr}, we deduce
\begin{align}
	\normmm{F_h(t)}\leq C(T^*,E[\xhi_0],\norm{\bv_0}),\quad \forall t\in[0,T^*],\label{uniformbound}
\end{align}
and for any $T^*>0$.

\subsection{Limit as $h\to0$}
\subsubsection{Compactness as $h\to0$}\label{compact}
From \eqref{regvbis} and \eqref{bounds} we immediately deduce that there exist $ \tw\in L^2_{loc}([0,\infty);\Hz)$, $ \widetilde{\lambda} \in L^2_{loc}([0,\infty))$, $u\in L^2_{loc}([0,\infty);\Hz)$, and $\bv\in L^2_{loc}([0,\infty);\Hus)\cap L^\infty_{loc}([0,\infty);\bL_\sigma^2(\om))$ (thus globally defined on $(0,\infty)$), such that, up to subsequences, as $h\to0$, 
\begin{align}
\label{whconv}	&w_h\rightharpoonup  \tw \quad\text{weakly in }L^2(0,\Ts;\Hz),\\&
\nonumber	\lambda_h\rightharpoonup  \widetilde\lambda \quad\text{weakly in }L^2(0,\Ts),\\&
	u_h\rightharpoonup u \quad\text{weakly in }L^2(0,\Ts;\Hz),\label{uhA}\\&
		\bv_h\rightharpoonup \bv\quad \text{weakly in } L^2(0,	\Ts;\Hus,\label{Hmz}
	\\&
	\bv_h\rightharpoonup \bv\quad \text{weakly* in } L^\infty(0,\Ts;\bL^2_\sigma(\om)),\label{velox}
\end{align}
for any $\Ts>0$, so that, clearly,
 \begin{align}
 	w_h+\lambda_h\rightharpoonup \tw+\tlambda,\quad \text{weakly in }L^2(0,T^*;H^1(\Omega)),\label{sum1}
 \end{align}
 for all $T^*>0$.
 
 Proceeding in the compactness argument, since $F_h$ defined in \eqref{F} is monotone nonincreasing by \eqref{monoton} and uniformly bounded by \eqref{uniformbound}, by Helly's selection theorem there exists a nonincreasing $\tF$ such that, up to subsequences,
 	\begin{align}
 		F_h(t)\to \tF(t),\quad \text{as }h\to 0,\quad \forall t\geq 0.\label{Helly}
 	\end{align}
This will be needed to identify the limit of $\normmm{\mu_t^h}_{\mathbb S^{d-1}}(\overline{\Omega})$ as $h\to0$.  

Furthermore, by \eqref{energyfinal}, we have
\begin{align}
	\int_0^\Ts \normh{\tfrac{\xhi_t^{T_h(t)}-\xhi^h(t)\circ X_{-T_h(t)}^{\widehat{\bu}_h^k(t)}}{T_h(t)}}^2\dt\leq C(\Ts),\label{bounddt}
\end{align}
for any $\Ts>0$, and then 
\begin{align*}
	\int_0^\Ts \frac{1}{{T_h(t)}}\normh{{\xhi_t^{T_h(t)}-\xhi^h(t)\circ X_{-T_h(t)}^{\widehat{\bu}_h^k(t)}}}^2\dt\leq h\int_0^\Ts \normh{\tfrac{\xhi_t^{T_h}-\xhi^h(t)\circ X_{-T_h(t)}^{\widehat{\bu}_h^k(t)}}{T_h(t)}}^2\dt\leq C(\Ts)h\to 0,
\end{align*}
as $h\to 0$. Setting $q_h(t):=\frac{1}{{T_h(t)}}\normh{{\xhi_t^{T_h}-\xhi^h(t)\circ X_{-T_h(t)}^{\widehat{\bu}_h^k(t)}}}^2\geq 0$, this entails $q_h\to 0$ in $L^1(0,\Ts)$ for any $\Ts>0$, entailing, up to a subsequence,
\begin{align}
\frac{1}{{T_h(t)}}\normh{{\xhi_t^{T_h(t)}-\xhi^h(t)\circ X_{-T_h(t)}^{\widehat{\bu}_h^k(t)}}}^2\to 0,\quad \text{ for almost any }t\in (0,\infty).\label{ae}
\end{align}
In conclusion, we also observe that $w_h\rightharpoonup \tw$ weakly in $L^2(0,\Ts; \Hz)$ for any $\Ts>0$, and as noticed in \eqref{map2}, $X_s^{\widehat \bu_h^k(t)}$ satisfies the assumption \eqref{map1} of Lemma \ref{convh1}, with $\bv_n(t)=\widehat \bu_h^k(t)$ for any $t\in[0,\Ts]$. As by \eqref{basica} we have that 
\begin{align*}
	\norm{\widehat{\bu}_h^k}_{C([0,\Ts];\mathbf C^1(\overline\Omega))}\leq C(\Ts)h^{-\frac14},
\end{align*}
for any $\Ts>0$, we see that all the assumptions of the first part of Lemma \ref{convh1} are satisfied, with $\gamma=\frac14$ in \eqref{unifC1}, and thus we also infer from \eqref{w1} that 
\begin{align}
	\nabla(w_h\circ X_{h_h(\cdot)}^{\widehat{\bu}_h^k(\cdot)})\rightharpoonup \nabla \tw,\quad\text{ weakly in }L^2(0,\Ts;\mathbf L^2(\Omega)),\quad\text{ as }h\to 0,
	\label{fondamental}
\end{align}
for any $\Ts>0$.

We now need to prove a compactness result for the phase indicator. Differently from \cite{SH}, this case is more involved due to the presence of the flow map $X_s^\bv$. Our aim is to use the well-known Aubin-Lions-Simon compactness theorem. Let us first observe that, for any $\delta\in(0,\Ts)$, 
\begin{align}
	\int_0^{\Ts-\delta}\normh{\wchi^h(t+\delta)-\wchi^h(t)}^2\dt \leq  C\delta,\quad \forall h\in(0,1),
	\label{control1b}
\end{align}
where $C>0$ does not depend on $h$. Indeed, we can write, by the Fundamental Theorem of Calculus and Cauchy-Schwarz's inequality,
\begin{align}
\nonumber\int_0^{\Ts-\delta}\normh{\wchi^h(t+\delta)-\wchi^h(t)}^2\dt &\leq 
	 \int_0^{\Ts-\delta}\int_{t}^{t+\delta}\delta\normh{\pt \wchi^h(\tau)}^2\d\tau\dt\\&
	 \leq \delta(\Ts-\delta)
	\int_{0}^{	\Ts}\normh{\pt \wchi^h(\tau)}^2\d\tau
	\leq C_A(\Ts)\delta,\label{ineq1}
\end{align}
recalling \eqref{controlh}, for some $C_A(\Ts)$ independent of $h$.

Also, notice that, for almost any $t\in(0,\Ts)$, for $n=\floor{\tfrac th}$,
\begin{align}
&	\nonumber\normh{\chi^h(t)-{\wchi}^h(t)}=\normh{\chi_n^h-\wchi^h(t)}=\normh{\wchi^h(nh)-\wchi^h(t)}\\&
	\leq \int_{nh}^{t}\normh{\pt\wchi^h(\tau)}\d\tau\leq \sqrt{t-nh}\left(\int_0^\Ts \normh{\pt\wchi^h(\tau)}^2\right)^\onehalf\leq C(\Ts)\sqrt h,
	\label{relation}
\end{align}
entailing
\begin{align}
	\int_0^{\Ts}\normh{\chi^h(t)-{\wchi}^h(t)}^2\dt\leq C_B(\Ts)h,
	\label{controla}
\end{align}
for some $C_B(\Ts)$ independent of $h$.
Let us now fix $\epsilon>0$ and choose $\wtilde h_1>0$ such that $C_A(\Ts)\wtilde h_1\leq \epsilon$. We consider a generic sequence $h_n\to 0$. Then clearly there exists $M\subset \mathbb N$ of finitely many indices such that $h_n>\widetilde h_1$ for any $n\in M$ (whereas $h_n\leq\wtilde  h_1$ for any $n\in \N\setminus M$)). Since the right translation operator is continuous in $L^2(0,T^*-h;\Hmz)$ (this can be easily shown by approximation with continuous functions), there exists $\delta_1$ such that 
\begin{align}
	\int_0^{\Ts-\delta}\normh{\chi^{h_n}(t+\delta)-\chi^{h_n}(t)}^2\dt\leq \epsilon,\quad \forall n\in M,\quad \forall 0<\delta<\delta_1.\label{subseq}
\end{align}
Up to reducing the size of $\delta_1$, namely setting it to be bounded above by $\tfrac \epsilon{C_A}$, we have, by applying triangle inequality and by \eqref{ineq1},
\begin{align*}
	&\int_0^{\Ts-\delta}\normh{\chi^{h}(t+\delta)-\chi^{h}(t)}^2\dt\\&
	\leq 2\int_0^{\Ts-\delta}\normh{\chi^{h}(t+\delta)-\wchi^{h}(t+\delta)}^2\dt+2\int_0^{\Ts-\delta}\normh{\chi^{h}(t)-\wchi^{h}(t)}^2\dt\\&
	+2\int_0^{\Ts-\delta}\normh{\wchi^{h}(t+\delta)-\wchi^{h}(t)}^2\dt\\&
	\leq 6\epsilon,\quad \forall h\leq \wtilde h_1,\quad \forall 0<\delta<\delta_1.
\end{align*}
This, together with \eqref{subseq}, entails
\begin{align}
	\sup_{n\in \mathbb N}\int_0^{\Ts-\delta}\normh{\chi^{h_n}(t+\delta)-\chi^{h_n}(t)}^2\dt\to 0,\quad \text{ as }\delta\to0.
	\label{controlB}
\end{align}
Concerning $\chi_t^{T_h(t)}$, observe that it holds, recalling \eqref{H1c}, \eqref{basica}, and \eqref{bounddt},
\begin{align}
\nonumber	&\int_0^\Ts\normh{\xhi_t^{\Th{t}}-\xhi^h(t)}^2\dt \\&\leq 2\nonumber	\int_0^\Ts\Th{t}^2\normh{\frac{1}{T_h(t)}\left(\xhi_t^{\Th{t}}-\xhi^h(t)\circ X_{-T_h(t)}^{\widehat{\bu}_h^k(t)}\right)}^2\dt\\&\quad \nonumber+2\int_0^\Ts\normh{\xhi^h(t)-\xhi^h(t)\circ X_{-T_h(t)}^{\widehat{\bu}_h^k(t)}}^2\dt\\&
	\leq 2C(\Ts)h^2+2\int_0^\Ts T_h(t)^2\norm{\widehat{\bu}_h^k}_{C([0,T];\mathbf C(\overline{\Omega}))}^2\leq C_C(\Ts)(h^2+h^\frac32),\label{CC}
\end{align}
 where $C_C(\Ts)$ does not depend on $h$.
Again let us fix $\epsilon>0$ and choose $\wtilde h_2>0$ such that $C_C(\Ts)(\wtilde h_2^2+\wtilde h_2^\frac32)\leq \epsilon$. Let us consider a generic sequence $h_n\to 0$. Then clearly there exists $M_2\subset \mathbb N$ of finitely many indices such that $h_n>\wtilde h_2$ for any $n\in M_2$ (and $h_n\leq \wtilde h_2$ for any $n\in \N\setminus M_2$). By continuity of the right translation operator in $L^2(0,T^*-h;\Hmz)$, there exists $\delta_2$ such that 
\begin{align}
	\int_0^{\Ts-\delta}\normh{\chi^{T_{h_n}(t+\delta)}_{t+\delta}-\chi^{T_{h_n}(t)}_{t}}^2\dt\leq \epsilon,\quad \forall n\in M_2,\quad \forall 0<\delta<\delta_2.\label{subseqA}
\end{align}
Then, choosing the same subsequence as the one in \eqref{controlB}, by triangle inequality, exploiting \eqref{controlB} and \eqref{CC} we infer that  
\begin{align*}
	\sup_{h_n\leq \wtilde h_2}\int_0^{\Ts-\delta}\normh{\chi^{T_{h_n}(t+\delta)}_{t+\delta}-\chi^{T_{h_n}(t)}_{t}}^2\dt\leq 6\epsilon,
\end{align*}
for any $\delta>0$ sufficiently small.
Together with \eqref{subseqA} this implies,
\begin{align}
	\sup_{n\in\mathbb N}\int_0^{\Ts-\delta}\normh{\chi^{T_{h_n}(t+\delta)}_{t+\delta}-\chi^{T_{h_n}(t)}_{t}}^2\dt\to 0,\quad \text{ as } \delta\to0.\label{subseq1}
\end{align}
Note that, by \eqref{energyfinal}, we also have 
\begin{align*}
	\wchi^h(\cdot),\xhi^h(\cdot),\xhi_{\cdot}^{\Th{\cdot}}	\in L^\infty_{w*}(0,\Ts;BV(\om;\{0,1\})),
\end{align*}
with norms uniformly bounded in $h$, and since $BV(\om;\{0,1\})\hookrightarrow\hookrightarrow L^{p}(\Omega)\hookrightarrow H^1(\Omega)'$, where $p\in(\tfrac{2d}{d+2},\tfrac{d}{d-1})$ (since $d=2,3$), recalling \eqref{control1b}, \eqref{controlB}, \eqref{subseq1}, we can apply Aubin-Lions-Simon Theorem, together with a diagonal argument, to infer the existence of $\xhi\in L^2_{loc}([0,\infty);L^1(\Omega))$, such that, up to a subsequence
\begin{align*}
	\wchi^h(\cdot),\xhi^h(\cdot),\chi_\cdot^{\Th{\cdot}}\to \xhi(\cdot), \text{ strongly in }L^2(0,\Ts;L^1(\Omega)),
\end{align*}
for any $\Ts>0$, where the identification of the three limits has been possible thanks to \eqref{controla} and \eqref{CC}. This also entails, up to the same subsequence, 
\begin{align*}
	\wchi^h(\cdot),\xhi^h(\cdot),\chi_\cdot^{\Th{\cdot}}\to \xhi(\cdot), \text{ for almost any   }(x,t)\in \Omega\times(0,T^*),
\end{align*}
which entails that $\xhi\in \{0,1\}$ almost everywhere in $\Omega\times(0,\infty)$.
Moreover, since clearly all the functions are bounded in $[0,1]$, we infer 
\begin{align}
		\wchi^h(\cdot),\xhi^h(\cdot),\chi_\cdot^{\Th{\cdot}}\to \xhi(\cdot), \text{ strongly in }L^p(\Omega\times(0,\Ts)),\quad \forall p\in[1,\infty).
	\label{finalconv}
\end{align}
As observed in \cite{AbelsRoger}, since $\xhi^h$ is uniformly bounded in $L^\infty_{w*}(0,\Ts;BV(\om;\{0,1\}))$, this implies that $\xhi_h\to \xhi$ weakly* in  $L^\infty_{w*}(0,\Ts;BV(\om;\{0,1\}))$, and thus, $\xhi\in L^\infty_{w*,loc}([0,\infty);\mathcal{M}_{m_0})$.

Additionally, since $\pt \wchi^h$ is uniformly bounded in $L^2(0,\Ts;\Hmz)$ for any $\Ts>0$ by \eqref{ptchi}, we deduce $\xhi\in H^1_{loc}([0,\infty);\Hmz)$, and
\begin{align}
\pt \wchi^h\rightharpoonup \pt \xhi\quad \text{ weakly in }L^2(0,\Ts;\Hmz),
\label{dtchi1}
\end{align}
for any $\Ts>0$, entailing 
\begin{align}
	\xhi\in H^1(0,\Ts;\Hmz)\hookrightarrow C([0,\Ts];\Hmz),
	\label{regchi}
\end{align}
for any $\Ts>0$. Also, by means of Lemma \ref{embeddingimp} we see that 
\begin{align}
	\xhi\in C([0,\Ts];L^p(\om)),\quad \forall \Ts>0,\quad 	\forall p\in [1,\infty).
	\label{regchi1a}
\end{align}

 It is now immediate to observe that the condition for the initial values easily passes to the limit (see \cite[eq. (104)]{SH}), and allows us to infer that 
$$
\text{Tr}_{|t=0}\ \xhi=\xhi_0\quad \text{ in }\Hmz.
$$

Now we need to show the strong convergence of $\bv_h$ in $L^2(0,\Ts;\bL^2_\sigma(\om))$. This can be obtained by following  the arguments from \cite[Lemma 3.8]{AbelsLengeler}. Namely, from the uniform bounds \eqref{regv1b}, since $\bL^2_\sigma(\om)\hookrightarrow\hookrightarrow \Hus'\hookrightarrow \Hds'$, we can apply Aubin-Lions Theorem to infer 
\begin{align}
	\bP_\sigma(\rho_h\bv_h)\to \boldsymbol\zeta,\quad \text{strongly  in }L^2(0,\Ts;\Hus'),
	\label{zeta}
\end{align} 
for some $\boldsymbol\zeta\in L^2(0,\Ts;\Hus')$. Now from \eqref{finalconv} we infer that, also,
	\begin{align}
		\rho_h\to \rho,\quad \text{ strongly in }L^p(\Omega\times(0,\Ts)),
		\label{rhoconvergence}
	\end{align}
where $\rho:=\rho_1\xhi+\rho_2(1-\chi)$. Therefore, recalling \eqref{velox} and that $\bP_\sigma:L^2(0,\Ts;\bL^2(\om))\to L^2(0,T;\bL^2_\sigma(\om))$ is weakly continuous, we deduce $\boldsymbol\zeta=\bP_\sigma(\rho\bv)$. By means of \eqref{velox} and \eqref{zeta}, we also infer that
\begin{align}
	&\nonumber\int_0^\Ts\int_\Omega \rho_h\normmm{\bv_h}^2\dx\dt=	\int_0^\Ts\int_\Omega \bP_\sigma(\rho_h\bv_h)\cdot\bv_h\dx\dt\\&
	\to \int_0^\Ts\int_\Omega \bP_\sigma(\rho\bv)\cdot\bv\dx\dt=	\int_0^\Ts\int_\Omega \rho\normmm{\bv}^2\dx\dt,\quad \text{as }h\to0.
	\label{strongconv0}
\end{align} 
Additionally, since $\rho_h\geq \rho_*>0$, we infer from \eqref{rhoconvergence} and \eqref{strongconv0} that
\begin{align*}
	\rho_h^\frac12\bv_h\to \rho^\frac12\bv,\ \text{strongly in } L^2(0,\Ts;\bL^2(\om)),\text{ and } 	\rho_h^{-\frac12}\to \rho^{-\onehalf},\ \text{almost everywhere in }\om\times(0,\Ts).
\end{align*}
Therefore,
$$
\bv_h\to \bv, \quad \text{almost everywhere in }\om\times(0,\Ts),
$$
and thus, by Lebesgue's Dominated Convergence Theorem, we infer 
\begin{align}
	\bv_h\to \bv,\quad \text{ strongly in }L^2(\om\times(0,\Ts)),
	\label{velconv}
\end{align}
as $h\to0$, for any $\Ts>0$.

Additionally, recalling \eqref{regchi1a}, we can prove that $\bv\in C_w([0,\Ts];\bL_\sigma^2(\om))$ and $\bv(0)=\bv_0$ almost everywhere, by following the arguments in \cite[Section 5.2]{ADG}.


Now we need to address the convergence of $\widehat{\bu}_h^k$ to $\bv$. From \eqref{c1a} we recall that 
\begin{align}
	&\esssup_{t\in(0,\Ts)}\norm{\widehat{\bu}_h^k(t)-\bu_h^k(t)}^2\leq C(\Ts)h^\frac12,\label{ctrl1a}
\end{align}
and thus, from \eqref{fund1}, we infer that (along the same subsequence for which \eqref{velconv} holds) 
\begin{align}
\limsup_{h\to0}	\int_0^\Ts\norm{\widehat{\bu}_h^k(t)}^2\dt\leq \limsup_{h\to0}\int_0^\Ts\norm{\bv_h(t)}^2\dt= \int_0^\Ts\norm{\bv(t)}^2\dt. \label{limsupcontrol}
\end{align}
This also entails that $\widehat \bu_h^k\rightharpoonup \boldsymbol\zeta_1$ weakly in $L^2(0,\Ts;\Ls)$. If we identify $\boldsymbol\zeta_1=\bv$ almost everywhere in $\om\times(0,\infty)$, then, since $L^2(0,\Ts;\Ls)$ is a Hilbert space, from \eqref{limsupcontrol} we conclude that 
\begin{align}
\widehat \bu_h^k\to \bv,\quad \text{ strongly in }L^2(0,\Ts;\Ls),\label{strvkh}
\end{align}
for any $\Ts>0$, as $h\to0$.
To identify $\boldsymbol\zeta_1$, we come back to \eqref{bvk}. First we multiply it by $\mathbf w\in C^\infty_c((0,\Ts);\mathbf C^\infty_{c,\sigma}(\om))$ and integrate over $\om\times(0,\Ts)$, to obtain, after integration by parts where necessary, and recalling the estimates \eqref{controls},
\begin{align}
	\nonumber&\normmm{\int_0^\Ts\int_\om(\bu^k_h-M_k(\bv_h))\cdot\bw\dx\dt}\leq	\frac1k\normmm{\int_0^\Ts\int_\Omega \bA\bu^k_h\cdot \pt \bA\bw\dx\dt}+\frac1k\normmm{\int_0^\Ts\int_\om \bA\bu_h^k\cdot \bA\bw\dx\dt}\\&\nonumber
	\leq \frac1k\norm{\bu_h^k}_{L^2(0,\Ts;\Hds)}\left(\norm{\pt\bA\bw}_{L^2(0,\Ts;\bL^2(\om))}+\norm{\bA\bw}_{L^2(0,\Ts;\bL^2(\om))}\right)\\&
	\leq \frac{C(\Ts)}{\sqrt k}\left(\norm{\pt\bA\bw}_{L^2(0,\Ts;\bL^2(\om))}+\norm{\bA\bw}_{L^2(0,\Ts;\bL^2(\om))}\right)\to 0,\label{bw1}
\end{align}
as $h\to0$, where we recall that we set $k=h^{-\frac18}$.
Now, it holds, by \eqref{dissiprel}, for any component $i=1,\ldots,d$,
\begin{align}
	k\esssup_{T\leq \Ts}\mathcal L^d(\{x\in\om:\ \normmm{\bv_{h,i}(x,t)}>k\})\leq \esssup_{T\leq \Ts}\int_\om \normmm{\bv_h(t)}\dx\leq C(\Ts),
	\label{vanishing}
\end{align}
so that, by \eqref{dissiprel} and \eqref{regvbis}, since $\normmm{M_k(\bv_{h})_i}\leq \normmm{\bv_{h,i}}$ and $\norm{\bv_h}_{L^2(0,\Ts;\bL^4(\om))}\leq C(\Ts)$, 
\begin{align*}
	&\norm{M_k(\bv_{h})_i-\bv_{h,i}}_{L^2(0,\Ts;L^2(\om))}^2= \int_0^\Ts\int_{\{x\in\om:\ \normmm{\bv_{h,i}(x,t)}>k\}}\normmm{M_k(\bv_{h})_i-\bv_{h,i}}^2\dx\dt\\&
	\leq C(\Ts)\esssup_{T\leq \Ts}\mathcal{L}^d(\{x\in\om:\ \normmm{\bv_{h,i}(x,t)}>k\})^\frac12\norm{\bv_h}_{L^2(0,\Ts;\bL^4(\om))}^2\\&
	\leq \frac{C(\Ts)}{\sqrt k}\to 0,\quad \text{as }h\to0,
\end{align*}
for any $i=1,\ldots,d$, entailing in particular that $(\bP_\sigma M_k(\bv_h)-\bv_{h},\bw)\to 0$ as $h\to0$, for any $\bw\in C^\infty_c((0,\Ts);\mathbf C^\infty_{c,\sigma}(\om))$. This result, together with \eqref{velconv}, allows to deduce 
\begin{align*}
	(\bP_\sigma M_k(\bv_h)-\bv,\bw)\to 0,
\end{align*}
for any $\bw\in C^\infty_c((0,\Ts);\mathbf C^\infty_{c,\sigma}(\om))$, which, in combination with \eqref{ctrl1a} and \eqref{bw1}, gives 
\begin{align*}
	(\widehat{\bu}_h^k,\bw)\to (\bv,\bw),
\end{align*}
for any $\bw\in C^\infty_c((0,\Ts);\mathbf C^\infty_{c,\sigma}(\om))$, as $k\to\infty$.
This allows us to conclude that $\boldsymbol{\zeta}_1=\bv$, thus entailing \eqref{strvkh}.

To end this compactness section, notice that, recalling \eqref{fondamental}, \eqref{finalconv}, and \eqref{strvkh}, it is now immediate to infer that 
\begin{align}
	&\int_0^T\left( \xhi^h{(\tau)} \widehat{\bu}_h^k(\tau),\nabla \left(w_h(\tau)\circ X_{g_h(\tau)}^{\widehat{\bu}_h^k(\tau)}\right)\right)_\om\d \tau\to \int_{0}^T \left( \xhi{(\tau)} \bv(\tau),\nabla \tw(\tau)\right)_\om\d \tau,\label{reusable1}
\end{align}
as $h\to 0$, for any $T>0$.

Additionally, observe that, since \eqref{map2} holds, we can resort to the last statement of Lemma \ref{basic1}. Indeed, setting $\xhi_n=\xhi^h$, $\bv_n=\widehat{\bu}_h^k$, along a suitable subsequence for which \eqref{finalconv} and \eqref{strvkh} hold, assumptions \eqref{b1} and \eqref{b2} are satisfied thanks to \eqref{dissiprel} and \eqref{regvbis}. As a consequence, we obtain from \eqref{ww1} that 
\begin{align}
\frac {\xhi^{h}-\xhi^h\circ X_{-h}^{ \widehat{\bu}_h^k}}{h}\to \bv\cdot \nabla \xhi,\ \text{strongly in }L^2(0,\Ts; \Hmz),\label{partialchi}
\end{align}
as $h\to 0$, where we recall that $\langle\bv\cdot \nabla \xhi, \varphi\rangle_{\Hmz,\Hz}=-(\xhi\bv,\nabla\varphi)_\om$ for any $\varphi\in \Hz$. Recalling also \eqref{dtchi1}, this also entails 
\begin{align}
	\partial^{\bullet,h}_t\xhi_h=\partial_t\widehat\xhi^h+\frac {\xhi^{h}-\xhi^h\circ X_{-h}^{ \widehat{\bu}_h^k}}{h}\rightharpoonup \partial_t\xhi+\bv\cdot \nabla \xhi,\ \text{weakly in }L^2(0,\Ts; \Hmz).
	\label{pchia}
\end{align}

\subsubsection{Limit varifold I. Disintegration of the measure} 

Thanks to the uniform bound of the energy \eqref{dissiprel}, by the weak*-compactness properties of finite Radon measures, we infer that there exists $\widetilde{\mu}:{\mathcal B((0,\infty)	\times\overline\Omega\times \mathbb S^{d-1})}\to \mathbb R^+$, such that, up to subsequences, 
\begin{align}
	\mathcal L^1\llcorner{(0,\Ts)}\otimes (\mu_h^t)_{t\in(0,\Ts)}\rightharpoonup \widetilde{\mu},\quad \text{ weakly* in } \mathcal M((0,\Ts)\times\overline\Omega\times \mathbb S^{d-1}),\label{convmeas}
\end{align}
for any $\Ts>0$, where $\mathcal B((0,\infty)	\times\overline\Omega\times \mathbb S^{d-1})$ is the Borel $\sigma$-algebra on $(0,\infty)	\times\overline\Omega\times \mathbb S^{d-1}$. We aim at showing that the
limit measure $\widetilde\mu$ can be sliced in time.
Let us now select, once and for all, a nonrelabeled subsequence $h$ such that all the convergences above take place as $h\to0$.
To show the disintegration of the measure $\widetilde\mu$, we perform an argument in the spirit of \cite[Lemma 2]{LH}, namely we want to use a cylindrical test function. Let $\eta \in C_c((0,\infty))$ with $\eta\geq 0$ and $\zeta \in C_c(\overline{\Omega}\times \mathbb S^{d-1})$ with $\zeta \in [0, 1]$. Using the definitions of $\mu_t^h$ and $F_h$ (see \eqref{F}), and the fact that $\zeta\leq 1$, 
\begin{align}	
&\nonumber\int_{(0,\infty)\times \overline{\Omega}\times \mathbb S^{d-1}}\zeta(x,s)\eta(t)\d 	\left(\mathcal L^1\llcorner{(0,\infty)}\otimes (\mu_h^t)_{t\in(0,\infty)}\right)\\&\nonumber\leq \int_0^\infty \int_{\overline{\Omega}\times \mathbb S^{d-1}}\eta(t)\d \mu_h^t\d t\\&\nonumber
=\int_0^\infty \eta(t)\normmm{\mu_t^h}_{\mathbb S^{d-1}}(\overline{\Omega})\d t\\&\label{fl}
\leq \int_{\text{supp}(\eta)} \eta(t)F_h(t)\dt+\int_{\text{supp}(\eta)} \eta(t)\left(\int_0^t\left( \xhi^h{(\tau)} \widehat{\bu}_h^k(\tau),\nabla \left(w_h(\tau)\circ X_{g_h(\tau)}^{\widehat{\bu}_h^k(\tau)}\right)\right)\d \tau\right)\dt.
\end{align}
We need now to pass to the limit as $h\to0$. Since $\eta(t)$ is of compact support in $(0,\infty)$, there exists $\Ts>0$ such that $\text{supp}(\eta)\subset (0,\Ts)$. Since $\eta$ is also bounded, recalling \eqref{Helly}, this entails by dominated convergence that 
$$
\int_{\text{supp}(\eta)} \eta(t)F_h(t)\dt\to \int_{\text{supp}(\eta)} \eta(t)\widetilde F(t)\dt,\quad \text{ as }h\to 0.
$$
Also, note that, recalling \eqref{fondamental}, \eqref{finalconv}, and \eqref{strvkh}, as in \eqref{reusable1}, it is immediate to infer that 
\begin{align}
	\nonumber&\int_{\text{supp}(\eta)} \eta(t)\left(\int_0^t\left( \xhi^h{(\tau)} \widehat{\bu}_h^k(\tau),\nabla \left(w_h(\tau)\circ X_{g_h(\tau)}^{\widehat{\bu}_h^k(\tau)}\right)\right)\d \tau\right)\dt\\&\to \int_{\text{supp}(\eta)} \eta(t)\left(\int_0^t\left( \xhi{(\tau)} \bv(\tau),\nabla \tw(\tau)\right)\d \tau\right)\dt,\quad \text{ as }h\to 0\label{reusable}.
\end{align}
Therefore, we can pass to the limit in \eqref{fl}, and obtain, recalling the convergence \eqref{convmeas}, 
\begin{align}	
	&\nonumber\int_{(0,\infty)\times \overline{\Omega}\times \mathbb S^{d-1}}\zeta(x,s)\eta(t)\d \widetilde\mu\\&
	\leq \int_{\text{supp}(\eta)} \eta(t)\widetilde F(t)\dt+\int_{\text{supp}(\eta)} \eta(t)\left(\int_0^t\left( \xhi{(\tau)} \bv(\tau),\nabla \tw(\tau)\right)\d \tau\right)\dt.\label{bound1}
\end{align}
Since the right-hand side is finite for any nonnegative $\eta \in C_c((0, \infty))$, this implies that the projection of $\widetilde\mu$ onto the time variable is a Radon measure. In particular, we can disintegrate $\widetilde\mu$, meaning that there
exists a Radon measure $\sigma$ on $(0, \infty)$ and a weakly-* $\sigma$-measurable family of Radon probability measures $(\tau_t)_t$ on $\overline\Omega\times\mathbb S^{d-1}$, such that $\widetilde\mu = \sigma\otimes(\tau_t)_t$
. Therefore, \eqref{bound1} entails that $\sigma\leq   \left(\widetilde{F}(t)+\int_0^t\left( \xhi{(\tau)} \bv(\tau),\nabla \tw(\tau)\right)\d \tau\right)\mathcal L^1\llcorner{(0,\infty)}$, so
that, by the Radon-Nikodym Theorem, there exists an $\mathcal L^1$-measurable function $\widetilde\sigma$ on $(0, \infty)$ with $\widetilde\sigma(t) \leq \left(\widetilde{F}(t)+\int_0^t\left( \xhi{(\tau)} \bv(\tau),\nabla \tw(\tau)\right)\d \tau\right)$ for almost any $t \in (0, \infty)$ such that $\sigma = \widetilde{\sigma}\mathcal L^
1\llcorner(0, \infty)$. Thus we have obtained the desired representation
\begin{align}
	\widetilde{\mu}=\mathcal L^1\llcorner(0,\infty)\otimes (\widetilde\mu_t)_{t\in(0,\infty)},
	\label{representation}
\end{align}
with $\widetilde\mu_t=\widetilde{\sigma}(t)\tau_t$.

\subsubsection{Limit varifold II. Identification of the limit.} 
Thanks to the validity of \eqref{dissiprel} and \eqref{variation}, we can now follow the proof of \cite[Theorem 1, Step 6]{SH} (note that we can still apply \cite[Proposition 11]{SH} also in this setting, with, in the notation of the proposition, $\xhi_0=\xhi_{n-1}^h\circ X_{-\Th{t}}^{\bu_{n-1}^k}\in \mathcal M_{m_0}$, for any $n\in \mathbb N$ and any $t\in(0,\infty)$) to infer that (see \cite[eq. (110)]{SH})
\begin{align*}
	&\sup_{B\in C^1(\overline\om;\R^d), \norm{B}_{L^\infty(\Omega)}\leq1}\normmm{\delta \mu_t^{h,\Omega}}(B)+	\sup_{B\in C^1(\overline\om;\R^d), \norm{B}_{L^\infty(\Omega)}\leq1}\normmm{\delta \mu_t^{h,\partial\Omega}}(B)\\&\leq C(1+\normmm{\mu_t^h}_{\mathbb S^{d-1}}(\overline\Omega)^\frac32+\norm{\nabla w_h(t)}^\frac{3d}{2}).
\end{align*}
Note that we extended $\mu_{h}^{t,\Omega}\in \mathcal M(\Omega\times\mathbb S^{d-1})$ as to be $\mu_{h}^{t,\Omega}\in\mathcal M(\overline\Omega\times \mathbb S^{d-1})$ by defining $\mu_t^{h,\Omega}(A):=\mu_t^{h,\Omega}(A\cap (\Omega\times \mathbb S^{d-1}))$, for any Borel set $A\subset \overline\Omega\times \mathbb S^{d-1}$.
By means of \eqref{regvbis} and Fatou's lemma, we know that 
$$
t\mapsto \liminf_{h\to0} \norm{\nabla w_h(\cdot,t)}^2
$$
is an integrable function over $(0,\Ts)$, for any $\Ts>0$. Therefore, for almost any $t>0$, we can choose a suitable subsequence $\{h_i(t)\}_{i\in \N}$, with $h_i(t)\to0$ as $i\to\infty$, such that 
\begin{align}
\limsup_{i\in \N}\norm{\nabla w_{h_i(t)}(\cdot,t)}^2\leq C<+\infty.
\label{wi}
\end{align}
Along this subsequence we then apply Allard’s result on compactness for integer varifolds \cite{Allard}, deducing that, for almost any $t\in(0,\infty)$, there exists a (nonrelabeled) subsequence $h_i(t)\to0$ as $i\to\infty$ such that
\begin{align}
	\mu_t^{h_i(t),\Omega}\rightharpoonup \mu_t^{\Omega},\quad\text{weakly* in }\mathcal M(\overline\Omega\times \mathbb S^{d-1})\label{subseq1a}
\end{align}
with $\mu^{\Omega}_t$ a $(d-1)$-integer-rectifiable varifold. Analogously, we can also prove that
\begin{align}
	\mu_t^{h_i(t),\partial \Omega}\rightharpoonup \mu_t^{\partial\Omega},\quad\text{weakly* in }\mathcal M(\partial \Omega\times \mathbb S^{d-1}),\label{subseq1b}
\end{align}
where $\mu_t^{\partial\Omega}=g_t\mathcal H^{d-1}\llcorner \partial\Omega\otimes (\delta_{\mathbf n_{\partial\Omega}(x)})_{x\in \pom}$, with $g_t\in BV(\partial\Omega;\{0,1\})$. Without relabeling, we also extend $\mu_t^{\partial\Omega}$ to $\mathcal M(\overline\Omega\times \mathbb S^{d-1})$ by setting $\mu_t^{\partial\Omega}(A):=\mu_t^{\partial\Omega}(A\cap (\partial\Omega\times \mathbb S^{d-1}))$, for any Borel set $A\subset \overline\Omega\times \mathbb S^{d-1}$. 
We can now introduce 
\begin{align}
	\mu_t:=c_0\mu_t^{\Omega}+c_0\cos\gamma\mu_t^{\partial\Omega}\in \mathcal M(\overline{\Omega}\times \mathbb S^{d-1}),
	\label{mut}
\end{align}
entailing from \eqref{subseq1a}-\eqref{subseq1b} that 
\begin{align*}
	\mu_t^{h_i(t)}\rightharpoonup \mu_t,\quad\text{weakly* in }\mathcal M(\overline\Omega\times \mathbb S^{d-1}).
\end{align*}
Therefore, since $\overline{\Omega}\times \mathbb S^{d-1}$ is compact, it also holds, for almost any $t\in(0,\infty)$,
\begin{align}
	\normmm{\mu_t^{h_i(t)}}_{\mathbb S^{d-1}}(\overline\Omega)\to \normmm{\mu_t}_{\mathbb S^{d-1}}(\overline\Omega),\quad\text{ as }i\to\infty.
	\label{convtotvar}
\end{align}
Analogously, recalling \eqref{subseq1a}-\eqref{subseq1b},
\begin{align}
	\normmm{\mu_t^{h_i(t),\om}}_{\mathbb S^{d-1}}(\overline\Omega)\to \normmm{\mu_t^\om}_{\mathbb S^{d-1}}(\overline\Omega),\quad \normmm{\mu_t^{h_i(t),\pom}}_{\mathbb S^{d-1}}(\pom)\to \normmm{\mu_t^\pom}_{\mathbb S^{d-1}}(\pom), \quad\text{ as }i\to\infty.
	\label{convtotvar1}
\end{align}
Since, by construction, $\normmm{\mu_t^{h,\om}}_{\mathbb S^{d-1}}(\overline{\Omega})\leq \normmm{\mu_t^h}_{\mathbb S^{d-1}}(\overline{\Omega})$ and $\normmm{\mu_t^{h,\pom}}_{\mathbb S^{d-1}}(\pom)\leq \normmm{\mu_t^h}_{\mathbb S^{d-1}}(\pom)$ for any $h>0$, then, exploiting the pointwise convergences \eqref{convtotvar}-\eqref{convtotvar1}, we infer
\begin{align}
	\normmm{\mu_t^\om}_{\mathbb S^{d-1}}({\overline\Omega})\leq \normmm{\mu_t}_{\mathbb S^{d-1}}({\overline\Omega}),\quad \normmm{\mu_t^\pom}_{\mathbb S^{d-1}}(\pom)\leq \normmm{\mu_t}_{\mathbb S^{d-1}}(\pom)
	\label{inequality}
\end{align} 
for almost any $t>0$. Note also that, since $\normmm{\mu_t^{h}}_{\mathbb S^{d-1}}(\overline{\Omega})=c_0\normmm{\mu_t^{h,\om}}_{\mathbb S^{d-1}}(\overline{\Omega})+c_0\cos\gamma\normmm{\mu_t^{h,\pom}}_{\mathbb S^{d-1}}(\partial{\Omega})$ for any $h>0$, from \eqref{convtotvar} and \eqref{convtotvar1} we also infer $E[\mu_t]:=\normmm{\mu_t}_{\mathbb S^{d-1}}(\overline{\Omega})=c_0\normmm{\mu_t^{\om}}_{\mathbb S^{d-1}}(\overline{\Omega})+c_0\cos\gamma\normmm{\mu_t^{\pom}}_{\mathbb S^{d-1}}(\partial{\Omega})$, for almost any $t>0$.

We now identify the limit $\widetilde{F}$ in \eqref{Helly}, for almost any $t>0$. Indeed, recalling \eqref{Helly}, \eqref{ae} and \eqref{reusable1}, we easily infer from \eqref{convtotvar} that, for almost any $t\in(0,\infty)$,
\begin{align}
F_{h_i(t)}(t)\to \widetilde{F}(t)=\normmm{\mu_t}_{\mathbb S^{d-1}}(\overline\Omega)-\int_{0}^t \left( \xhi{(\tau)} \bv(\tau),\nabla \tw(\tau)\right)_\om\d \tau,\text{ as }i\to\infty.\label{Fh1}
\end{align}
Since $\wtilde F$ is nonincreasing, this means that
\begin{align}
	t\mapsto\normmm{\mu_t}_{\mathbb S^{d-1}}(\overline\Omega)-\int_{0}^t \left( \xhi{(\tau)} \bv(\tau),\nabla \tw(\tau)\right)_\om\d \tau\quad\text{ is nonincreasing in }(0,\infty),
	\label{nonincreasing}
\end{align}
and, from \eqref{Helly},
\begin{align}
F_{h}(t)\to \widetilde{F}(t)=\normmm{\mu_t}_{\mathbb S^{d-1}}(\overline\Omega)-\int_{0}^t \left( \xhi{(\tau)} \bv(\tau),\nabla \tw(\tau)\right)_\om\d \tau,\text{ as }h\to 0,\label{Fh1b}
\end{align}
for any $t\geq 0$.
\subsubsection{Gibbs-Thomson law as $h\to0$ and generalized mean curvature.}
\label{gibbs}
For any $\Ts>0$, let us consider a vector field $ B:\overline\om\times(0,\Ts)\to \R^d$ continuous and compactly supported in $(0,\Ts)$, such that, for any $t\in(0,\Ts)$, $B(\cdot,t)\in C^1(\overline{\Omega},\mathbb R^d)$ with $(B(\cdot,t)\cdot \mathbf n_{\partial \Omega})|_{\partial\Omega}\equiv0$. Then by \eqref{convmeas} and \eqref{representation} it holds 
\begin{align}
	\int_0^\Ts \delta\mu_t^h(B)\dt\to \int_0^\Ts \delta\widetilde{\mu}_t(B)\dt. \label{v1}
\end{align}
As an immediate consequence, we can multiply the Gibbs-Thomson relation \eqref{variation} by a smooth and
compactly supported test function on $(0, \Ts)$, integrate in time, and pass to the limit as $h\to0$, by means of \eqref{sum1} and \eqref{finalconv}, and then, after again a localization in time, deduce that, for almost any  $t>0$, it holds 
\begin{align}
	&	\delta\wtilde\mu_t(B)=\int_{\overline\Omega\times \mathbb S^{d-1}}(Id-s\otimes s):\nabla B(x)\d \widetilde{\mu}_t(x,s)\nonumber\\&=
	\int_\Omega \chi(\cdot,t)\text{div}((\tw(\cdot,t)+\tlambda(t))B)\d x,\label{variation2}
\end{align} 
for any $B\in C^1(\overline{\Omega};\R^d)$ with $B\cdot\mathbf n_{|\partial\Omega}\equiv 0$ (see also \cite[Eq. (118)]{SH}). By Allard's first variation formula the left-hand side of the identity above is $\delta\widetilde{\mu}(B)$. Additionally, by Proposition \ref{prop5}, we can also write \eqref{variation2} as 
\begin{align}
	\delta\widetilde{\mu}_t(B)=-\int_\Omega c_0\frac{\tw(\cdot,t)+\tlambda(t)}{c_0}\frac{\nabla\xhi(\cdot,t)}{\normmm{\nabla\chi(\cdot,t)}}\cdot B\d \normmm{\nabla\xhi(\cdot,t)}.\label{v2}
\end{align}
Then an application of \cite[Lemma 4.2]{Roger} gives that, for any $B\in L^2(0,\Ts;\mathbf C^1_c(\Omega))$,
\begin{align}
	\int_0^\Ts \delta\mu_t^h(B) \to -\int_0^\Ts\int_\Omega c_0H_{\xhi(\cdot,t)}\frac{\nabla\xhi(\cdot,t)}{\normmm{\nabla\chi(\cdot,t)}}\cdot B\d \normmm{\nabla\xhi(\cdot,t)}\dt,\label{cona}
\end{align}
for any $\Ts>0$, where $H_\xhi$ is the generalized mean curvature in the sense of \cite{Roger} (see also Definition \ref{admissibility}), which is intrinsic to the surface $\text{supp}\normmm{\nabla\xhi(\cdot,t)}$ for almost every $t>0$. In particular, we also deduce that 
\begin{align}
	H_{\xhi(x,t)}=Tr_{\text{supp}\normmm{\nabla \xhi(,t)}}\frac{ \tw(\cdot,t)+\tlambda(t)}{c_0}(x),\label{Hrog}
\end{align}
for $\normmm{\nabla \xhi(\cdot,t)}$-almost any $x\in \om$ and almost any $t>0$.

Recalling \eqref{v1} and \eqref{v2}, by localizing \eqref{cona} in time, we additionally infer that, since $\Ts>0$ is arbitrary, for almost any $t>0$ it holds
\begin{align}
\delta\widetilde{\mu}_t(B)= 	-\int_\Omega c_0H_{\xhi(\cdot,t)}\frac{\nabla\xhi(\cdot,t)}{\normmm{\nabla\chi(\cdot,t)}}\cdot B\d \normmm{\nabla\xhi(\cdot,t)},
	\label{equivA}
\end{align}  
$B\in C^1_c({\Omega};\R^d)$.
Thanks to Proposition \ref{prop5} (or by means of \cite{Schatzle}), we also infer that, for almost any $t>0$,
\begin{align}
H_{\xhi(\cdot,t)}\in L^s(\Omega; d\normmm{\nabla\xhi}),
\label{chi1A}
\end{align}
where $s\in[2,4]$ if $d=3$, and $s\in [2,\infty)$ if $d=2$. 
Additionally, from \eqref{variation2} and \eqref{Hrog} we see that for almost any  $t>0$, it holds 
\begin{align}
	&	\delta\wtilde{\mu}_t(B)=\int_{\overline\Omega\times \mathbb S^{d-1}}(Id-s\otimes s):\nabla B(x)\d \widetilde{\mu}_t(x,s)\nonumber\\&\nonumber=
	-\int_\Omega (\tw(\cdot,t)+\tlambda(t))B\cdot \frac{\nabla \xhi(\cdot,t)}{\normmm{\nabla\xhi(\cdot,t)}}\d \normmm{\nabla\xhi(\cdot,t)}\\&
    =-\int_\Omega c_0H_{\xhi(\cdot,t)}\frac{\nabla\xhi(\cdot,t)}{\normmm{\nabla\chi(\cdot,t)}}\cdot B\d \normmm{\nabla\xhi(\cdot,t)}
    ,\label{variation3}
\end{align} 
for any $B\in C^1(\overline{\Omega};\R^d)$ with $B\cdot\mathbf n_{|\partial\Omega}\equiv 0$.

\subsubsection{Identification of $\delta\mu_t$.} In the previous argument we have seen how to characterize $\delta\widetilde\mu_t$, but we still need to show that this variation actually coincides with $\delta\mu_t$, where $\mu_t$ is defined in \eqref{mut}. 
Recalling \eqref{complete}, and considering, for almost any $t>0$, the same subsequence $\{h_i(t)\}_i$ for which \eqref{wi} holds, we can find, for almost any $t>0$, $\widehat w(t)\in H^1(\Omega)$ such that, up to a further subsequence, $w_{h_i(t)}+\lambda_{h_i(t)}\rightharpoonup \widehat w(t)$ weakly in $H^1(\Omega)$ as $i\to\infty$. Also, by \eqref{finalconv}, for almost any $t$ we have $\xhi_t^{T_h(t)}(\cdot,t)\to \xhi(\cdot,t)$ strongly in $L^1(\Omega)$ as $h\to0$. Furthermore, for any $B\in C^1_c(\Omega;\R^d)$, it holds $\delta\mu^h_t(B)=\delta\mu_t^{h,\Omega}(B)$, for almost any $t>0$. Then, also exploiting the convergences detailed in the previous sections, we can apply the results of \cite[Lemma 4.1]{Roger} and
\cite[Theorems 1.1 and 1.2]{Schatzle} to the sequence $\{(\xhi^{T_{h_i(t)}(t)}_t,\vert{\mu_t^{h_i(t),\Omega}\vert}_{\mathbb S^{d-1}}\llcorner \Omega, w_{h_i(t)}+\lambda_{h_i(t)})\}_i$, to obtain, additionally to \eqref{subseq1a}-\eqref{Fh1}, that 
\begin{itemize}
	\item The $(d-1)$-integer-rectifiable measure $\tmu\llcorner\Omega\in \mathcal{M}(\Omega)$ is such that the compatibility condition \eqref{PC} holds true.
		\item The trace of the function $\frac{ \widehat w(\cdot,t)}{c_0}$ coincides  $\normmm{\mu_t^\Omega}_{\mathbb S^{d-1}}\llcorner \Omega$-almost everywhere in $\Omega$ with the generalized mean curvature vector $\bH_{\tmu\llcorner \Omega}:\text{supp}(\tmu\llcorner\Omega)\to \R^d$
		(which exists well defined), i.e.,
	\begin{align}
		\label{Hmu}\bH_{\tmu\llcorner\Omega}=\begin{cases}
			\frac{ \widehat w(\cdot,t)}{c_0}\frac{\nabla\xhi(\cdot,t)}{\normmm{\nabla\xhi(\cdot,t)}},\quad\text{if }x\in \text{supp}\normmm{\nabla\xhi(\cdot,t)},\\
			0,\quad \text{otherwise}.
		\end{cases}
	\end{align}
	\item Since the mean curvature $H_\xhi(t)$ defined above is intrinsic to $\text{supp}\normmm{\nabla \xhi(\cdot,t)}$, and both the the $(d-1)$-integer rectifiable varifolds associated to $\normmm{\wtilde\mu_t}_{\mathbb S^{d-1}}\llcorner \om$ and $\normmm{\mu_t^{\om}}_{\mathbb S^{d-1}}\llcorner \om$ satisfy the assumptions of \cite[Proposition 1.1]{Roger}, the mean curvature vector $\bH_{\tmu\llcorner \Omega}$ coincides with the curvature vector associated to $H_\xhi$ in the sense of \eqref{PJ}-\eqref{PK}.
	Also, \eqref{PH} holds for all compactly supported variations $B\in C^1_c(\Omega;\R^d)$. In particular, property \eqref{PG} holds true due to \eqref{chi1A}.

\end{itemize} 

It is now clear that the couple $(\xhi,\mu)$ is the right candidate to be the one in the notion of admissible couple of Definition \ref{admissibility}. 
To extend the validity of \eqref{PH} for all tangential variations, we  pass to the limit in \eqref{variation}, by means of \eqref{subseq1a}-\eqref{subseq1b} and the convergences above, namely that $w_{h_i(t)}+\lambda_{h_i(t)}\rightharpoonup \widehat w(t)$ weakly in $H^1(\Omega)$ for almost any $t>0$ as $i\to\infty$. We then get, also recalling \eqref{Hmu} and the property that $H_\xhi(t)$ is intrinsic to $\text{supp}\normmm{\nabla \xhi(\cdot,t)}$,
\begin{align}
	&\nonumber	\delta\mu_t(B)=\int_{\overline\Omega\times \mathbb S^{d-1}}(Id-s\otimes s):\nabla B(x)\d \mu_t(x,s)\\&=\nonumber
	\int_\Omega \chi(\cdot,t)\text{div}(\widehat w(\cdot,t)B)\d x=-\int_\Omega \widehat w(\cdot,t)B\cdot \frac{\nabla \xhi(\cdot,t)}{\normmm{\nabla\xhi(\cdot,t)}}\d \normmm{\nabla\xhi(\cdot,t)}\\&
    =-\int_\Omega c_0H_{\xhi(\cdot,t)}\frac{\nabla\xhi(\cdot,t)}{\normmm{\nabla\chi(\cdot,t)}}\cdot B\d \normmm{\nabla\xhi(\cdot,t)},\label{variationF}
\end{align} 
for all $B\in C^1(\overline\Omega;\mathbb R^d)$ with $(B\cdot \bn_{\partial\Omega})|_{\partial\Omega}\equiv 0$, for almost any  $t>0$. This gives \eqref{PH} for all $B\in C^1(\overline\Omega;\R^d)$. 

With these properties at hand it is now straightforward to see, by means of \eqref{variation3}, that  
 \begin{align}
 	\delta{\mu}_t(B)= 	-\int_\Omega c_0H_{\xhi(\cdot,t)}\frac{\nabla\xhi(\cdot,t)}{\normmm{\nabla\chi(\cdot,t)}}\cdot B\d \normmm{\nabla\xhi(\cdot,t)}=\delta\wtilde\mu_t(B),
 	\label{equivB1}
 \end{align}
for all $B\in C^1_c(\Omega;\mathbb R^d)$, and thus $\delta\mu_t(B)=\delta\widetilde{\mu}_t(B)$, for almost any  $t>0$. It then holds, from \eqref{variation2}, that
\begin{align}
	&	\delta\mu_t(B)=\delta\wtilde\mu_t(B)=
	\int_\Omega \chi(\cdot,t)\text{div}((\tw(\cdot,t)+\tlambda(t))B)\d x,\label{variation4}
\end{align} 
for any $B\in C^1(\overline{\Omega};\R^d)$ with $B\cdot\mathbf n_{|\partial\Omega}\equiv 0$, for almost any $t>0$.

To show the validity of property \eqref{PI}, it is enough to observe that, thanks to the results above, we can apply Proposition \ref{prop5} to the varifold $\mu_t$ (recalling its decomposition \eqref{mut} and the properties above), and infer property \eqref{PI} from \eqref{totalvar}.  

\subsubsection{Admissibility of the couple $(\xhi,\mu)$.}
We now explicitely verify that the couple $(\xhi,\mu)$ is admissible in the sense of Definition \ref{admissibility}. In particular, Property (1) is immediate after the definition \eqref{mut} of $\mu_t$, and the integer rectifiability of the corresponding measures for almost any $t\in(0,\Ts)$, and any $\Ts>0$, is given after \eqref{subseq1a}-\eqref{subseq1b}.  Concerning Property (2), the relation \eqref{PC} has been already shown in the previous section. Concerning the relations \eqref{PD}-\eqref{PF} they can be easily shown by means of the corresponding definition of $\mu_t^{h,\Omega}$ and $\mu_t^{h,\partial\Omega}$, the Gau{\ss} Theorem for BV functions and passing to the limit as $h\to 0$, see \cite[Step 12, Proof of item (ii) of Definition 3]{SH}. Then properties (3)-(4) of Definition \ref{admissibility}, i.e., \eqref{PG}-\eqref{PK}, they have just been shown in the previous section. In conclusion, Property (5), i.e., the measurability in time of $\mu_t$ is a direct consequence of \eqref{nonincreasing}, since the map $t\mapsto \int_{0}^t \left( \xhi{(\tau)} \bv(\tau),\nabla \tw(\tau)\right)\d \tau$ is also measurable.   This allows to conclude that the couple $(\xhi,\mu)$ is admissible in the sense of Definition \ref{admissibility}.  

\subsubsection{Validity of the (optimal) energy dissipation inequality.} 
By the compactness results of Section \ref{compact}, together with \eqref{Fh1b}, the definition of $u_h$ in \eqref{uh1}, and lower semicontinuity of the norms involved, we can now pass to the limit in \eqref{energy1abc3}, first taking $h\to0^+$, then $\tau\to T$ and $\kappa\to s$,  to infer  
\begin{align}
	\nonumber	&E[{\mu_T}]+\frac1{2}\int_{s}^{T}\norm{\nabla u(t)	}^2\dt+\frac{1}{2}\int_{s}^{T}\norm{\nabla \tw(t)}^2\d t\\& \leq E[{\mu_s}]+\int_{s} ^{T} \int_\om \xhi{(t)}  \nabla \tw(t)	\cdot \bv(t)\d t,
	\label{energy1abc3b}
\end{align}
for almost any $T\in (0,\Ts]$, almost any $s\in (0,T)$, and any $\Ts>0$. Note that we have used the fact that $\floor{\frac Th}h\to T$ as $h\to0$, as well as 
\begin{align*}
	\frac{T_h(s)}{2}\norm{\nabla w_h(s)}^2+\frac{1}{2T_h(T)}\normh{\xhi^h(T)-\xhi^h(T)\circ  X_{-T_h(T)}^{\widehat{\bu}_h^k(T)} }^2\to 0,
\end{align*}
for almost any $s\in(0,T)$ and almost any $T\in(0,\Ts]$, recalling the definition of  $w_h$ in \eqref{wh}, together with \eqref{convh1b} and \eqref{ae}. The validity of \eqref{energy1abc3b} for $s=0$ can be inferred analogously from \eqref{energy1abc2b}, where we recall that we have set $E[\mu_t]:=E[\xhi_0]$ when $t=0$.

Concerning the velocity part, we have from the convergences \eqref{Hmz} and \eqref{finalconv}, that, by weak lower semicontinuity of the norms, for any $\Ts>0$, 
\begin{align*}
&\int_s^T\norm{\sqrt{\nu(\xhi)}D \bv}^2\dt\leq \liminf_{h\to0}\int_s^T\norm{\sqrt{\nu(\xhi^h)}D \bv_h}^2\dt.
\end{align*} 
Also, recalling \eqref{finalconv} and \eqref{strongconv0}, for almost any $t>0$, it holds  
\begin{align}
	\frac12\int_\Omega \rho_h(t)\normmm{\bv_h(t)}^2\dx \to 	\frac12\int_\Omega \rho(t)\normmm{\bv(t)}^2\dx,\quad\text{ as }h\to 0.\label{converg}
\end{align}
Morever it is easy to infer, from \eqref{whconv}, \eqref{finalconv}, and \eqref{velconv} that, for any $0<s<T\leq\Ts$,
\begin{align*}
	-\int_s^T\int_\Omega \nabla w_h(t)\xhi^h(t)\cdot \bv_h(t) \dx\dt \to  -\int_s^T\int_\Omega \nabla \tw(t)\xhi(t)\cdot \bv(t) \dx\dt,\quad \text{as }h\to0.
\end{align*}

Therefore, recalling that $\beta\int_s^T\norm{\bA\bv_h}^2\dt\geq0$, we can pass to the limit as $h\to0$ in \eqref{energyv2}, to obtain, for almost any $0<s<T\leq \Ts$,

\begin{align}
	&\onehalf\int_\Omega {\rho}(T)  \normmm{\bv(T)}^2\dx+\int_s^T\norm{\sqrt{\nu(\xhi(t))}\nabla \bv}^2\dt\nonumber
	\\&\leq \onehalf\int_\Omega {\rho}(s)  \normmm{\bv(s)}^2\dx
	-\int_s^T\int_\Omega \nabla \tw(t)\xhi(t)\cdot \bv(t) \dx\dt.\label{energyv2b1}
\end{align}

In conclusion, notice that, since by \eqref{regchi1a} and recalling the regularity of $\bv$ we have that the function $t\mapsto \int_\Omega {\rho}(t)  \normmm{\bv(t)}^2\dx$ is lower semicontinuous, we can deduce that \eqref{energyv2b1} also holds for any $T>0$.

The case when $s=0$ is analogous, and can be easily obtained setting first $s=0$ in \eqref{energyv2} and then passing to the limit as $h\to0$ in an analogous way. 

\subsubsection{Limit systems as $h\to0$.} 
First, we notice that, from the definition of \eqref{uh1}, it holds, for any $\zeta\in C_c^\infty(\overline{\Omega}\times[0,\Ts))$,
\begin{align*}
\int_0^{\Ts}\int_\Omega\partial_t^{\bullet,h}\xhi^h(t)\zeta(t)\dx\dt=-\int_0^{\Ts}\int_\Omega \nabla u_h(t)\cdot \nabla\zeta(t)\dx\dt,
\end{align*}
entailing from \eqref{uhA} and \eqref{pchia}, by recalling  the regularity $\zeta\in C_c^\infty(\overline{\Omega}\times[0,\Ts])$,  that 
\begin{align}
	&\int_0^{\Ts}\langle\partial_t\xhi(t),\zeta(t)\rangle_{H^1(\om)',H^1(\om)}\dt-	\int_0^{\Ts}\int_\Omega\xhi(t)\bv(t)\cdot\nabla\zeta(t)\dx\dt\nonumber\\&=-\int_0^{\Ts}\int_\Omega \nabla u(t)\cdot \nabla\zeta(t)\dx\dt,
	\label{finalA}
\end{align}
giving \eqref{kinetic}, simply after an integration by parts in time.
Concerning the velocity equation, we notice that we can write
\begin{align}
	\partial^{\bullet,h}_t\rho_h=\partial_t\rho_h+\frac {\rho^{h}-\rho^h\circ X_{-h}^{ \widehat{\bu}_h^k}}{h},
	\end{align}
where we denoted by $\rho^h$ the piecewise constant interpolant corresponding to the sequence $\{\rho^{h}_n\}_n$. Therefore, we can rewrite some terms in the weak formulation \eqref{velocity1} as 
\begin{align}
\nonumber	\frac12\langle\left(\partial_t^{\bullet,h}\rho_h-\pt \rho_h\right),\bv_h\cdot\bxi\rangle_{H^1(\om)',H^1(\om)}&\nonumber=\frac12\langle\tfrac {\rho^{h}-\rho^h\circ X_{-h}^{ \widehat{\bu}_h^k}}{h},\bv_h\cdot\bxi\rangle_{H^1(\om)',H^1(\om)}\\&=\frac12\int_\om \frac {\rho^{h}-\rho^h\circ X_{-h}^{ \widehat{\bu}_h^k}}{h}\bv_h\cdot\bxi\dx,\label{reformulation}
\end{align} 
for any $\bxi\in \Hds$.
Observe now that, since $\bxi\in C_c^\infty((0,\infty);\Hds)$, by means of \eqref{Hmz} we infer 
\begin{align}
	\bv_h\cdot\bxi\rightharpoonup \bv\cdot\bxi,\quad\text{ weakly in }L^2(0,\Ts;H^1(\om)),\ \text{ as }h\to 0,\label{weakvxi}
\end{align}
for any $\Ts>0$. Recalling \eqref{partialchi} and the definition of $\rho_h$, it holds
\begin{align}
 \frac {\rho^{h}-\rho^h\circ X_{-h}^{ \widehat{\bu}_h^k}}{h}\to \bv\cdot\nabla \rho,\quad \text{ strongly in }L^2(0,\Ts;\Hmz),
	\label{rhoconv}
\end{align}
for any $\Ts>0$.
To pass to the limit in \eqref{reformulation} (integrated in time over $(0,t)$, $t>0$) we can then write 
\begin{align*}
	&	\int_0^t\int_\om \left(\tfrac {\rho^{h}-\rho^h\circ X_{-h}^{ \widehat{\bu}_h^k}}{h}\bv_h\cdot\bxi+\rho\bv\cdot  \nabla(\bxi\cdot\bv)\right)\dx\d s\\&=-\int_0^t	\int_\om \rho\bv\cdot \nabla(\bxi\cdot(\bv_h-\bv))\dx\d s+	\int_0^t\langle\tfrac {\rho^{h}-\rho^h\circ X_{-h}^{ \widehat{\bu}_h^k}}{h}-\nabla \rho\cdot \bv,\ \bxi\cdot\bv_h-\tfrac{\int_\om\bv_h\cdot\bxi\dx}{\mathcal L^d{(\Omega)}}\rangle_{\Hmz,\Hz},
\end{align*}
and the first term in the right-hand side converges to 0, due to \eqref{weakvxi} and the fact that $\rho\bv\in L^\infty(0,\Ts;\Ls)$, for any $\Ts>0$, by \eqref{dissiprel}. On the other hand, for the second term we  have, for any $t\in(0,\Ts)$ and any $\Ts>0$, by Poincaré's inequality,
\begin{align*}
	&\normmm{\int_0^t\langle\tfrac {\rho^{h}-\rho^h\circ X_{-h}^{ \widehat{\bu}_h^k}}{h}-\nabla \rho\cdot \bv,\ \bxi\cdot\bv_h-\tfrac{\int_\om\bv_h\cdot\bxi\dx}{\mathcal L^d{(\Omega)}}\rangle_{\Hmz,\Hz}}\\&
	\leq C \norm{\frac {\rho^{h}-\rho^h\circ X_{-h}^{ \widehat{\bu}_h^k}}{h}-\nabla \rho\cdot \bv}_{L^2(0,\Ts;\Hmz)}\norm{\nabla (\bxi\cdot\bv_h)}_{L^2(0,\Ts;\bL^2(\om))} \\&
	\leq C(\Ts)\norm{\frac {\rho^{h}-\rho^h\circ X_{-h}^{ \widehat{\bu}_h^k}}{h}-\nabla \rho\cdot \bv}_{L^2(0,\Ts;\Hmz)}\to 0,\quad\text{ as }h\to0,
\end{align*}
where we used \eqref{rhoconv} as well as the uniform bounds on $\nabla\bv_h$ from \eqref{regvbis}. 

Note now that, for any $\bxi\in C_c^\infty((0,\infty);\Hds)$, by \eqref{regvbis},
\begin{align*}
	\normmm{\int_0^t\int_\Omega \beta\bA\bv_h\cdot \bA\bxi\dx\d s}&\leq \beta\norm{\bv_h}_{L^2(0,\Ts;\Hds)}\norm{\bxi}_{L^2(0,\Ts;\Hds)}\\&\leq C(\Ts)\sqrt\beta \norm{\bxi}_{L^2(0,\Ts;\Hds)}\to 0,
\end{align*}
as $h\to0$, for any $t\leq \Ts$ and any $\Ts>0$, where we recall that we set $\beta=h^\frac18$.

Additionally, we have, for any $\bxi\in C_c^\infty((0,\infty);\Hds)$, 
\begin{align*}
	&	\normmm{\int_0^t\int_\om \rho_h(\bxi\cdot \nabla)\bv_h\cdot \bv_h\dx\d s}\\&\leq 	\normmm{\int_0^t\int_\om (\rho_h-\rho)(\bxi\cdot \nabla)\bv_h\cdot \bv\dx\d s}+	\normmm{\int_0^t\int_\om \rho_h(\bxi\cdot \nabla)\bv_h\cdot (\bv_h-\bv)\dx\d s}\\&\quad +	\normmm{\int_0^t\int_\om \rho(\bxi\cdot \nabla)(\bv_h-\bv)\cdot \bv\dx\d s}\\&
	\leq 
	C\norm{\rho_h-\rho}_{L^6(\om\times(0,\Ts))}\norm{\bxi}_{L^\infty(0,\Ts;\bL^\infty(\om))}\norm{\bv_h}_{L^2(0,\Ts;\Hus)}\norm{\bv}_{L^3(\om\times(0,\Ts))}\\&
	\quad +C\norm{\rho_h}_{L^\infty(\om\times(0,\Ts))}\norm{\bxi}_{L^\infty(0,\Ts;\bL^\infty(\om))}\norm{\bv_h}_{L^2(0,\Ts;\Hus)}\norm{\bv_h-\bv}_{L^2(0,\Ts;\bL^2_\sigma(\om))}\\&
	\quad +\normmm{\int_0^t\int_\om (\nabla(\bv_h-\bv)\cdot \bxi)\cdot(\bv\rho)\dx\d s}\to 0,\quad \text{as }h\to0,
\end{align*}
where we made use of \eqref{regvbis}, \eqref{Hmz}, \eqref{rhoconvergence}, \eqref{velconv}, as well as the embedding $ L^\infty(0,\Ts;\bL^2_\sigma(\om))\cap L^2(0,\Ts;\Hus)\hookrightarrow L^3(\om\times(0,\Ts))$ and the fact that $\rho\bv\in L^2(0,\Ts;\bL^2(\om))$, which is used in the convergence of the last summand (recalling \eqref{Hmz}). All the convergences in the advective terms in the weak formulation \eqref{velocity1} can be treated analogously.

Therefore, we can pass to the limit as $h\to0$ in \eqref{velocity1}, after integrating over $(0,t)$, $t>0$, to deduce
\begin{align}
	\nonumber	&\nonumber -\int_0^t\rho\bv\cdot\pt\bxi\dx\d s \nonumber-\int_0^t\int_\Omega \rho(\bxi\cdot \nabla)\bv\cdot\bv\dx\d s-\int_0^t\int_\Omega \rho(\bv\cdot \nabla)\bxi\cdot\bv\dx\d s\\&\nonumber+\int_0^t\int_\Omega{\rho}(\bv\cdot \nabla )\bv\cdot \bxi\dx\d s+\int_0^t\int_\Omega(\bv\otimes ({\rho_1-\rho_2})\nabla u):\nabla \bxi\dx \d s +\int_0^t\int_\Omega\nu(\xhi)D \bv:D \bxi\dx\d s 
	\nonumber\\&=\int_\om \rho_0\bv_0\cdot \bxi\dx -\int_0^t\int_\Omega\xhi\nabla \wtilde{w}\cdot \bxi\dx\d s,\quad \forall \bxi\in\mathbf C^\infty_c([0,\infty);\Hds).\label{velocity2k}
\end{align}

Note now that

\begin{align}
	\nonumber	&-\int_0^t\int_\Omega \rho(\bxi\cdot \nabla)\bv\cdot\bv\dx\d s-\int_0^t\int_\Omega \rho(\bv\cdot \nabla)\bxi\cdot\bv\dx\d s+\int_0^t\int_\Omega{\rho}(\bv\cdot \nabla )\bv\cdot \bxi\dx\d s\\&=-\int_0^t\int_\om \rho\bv\otimes \bv:\nabla \bxi\dx\d s.\label{identityk}
\end{align}
Indeed, this can be obtained by approximating $\rho,\bv$ with smooth functions, then observing that, after an integration by parts, recalling the boundary conditions,
\begin{align*}
	&-\int_0^t\int_\Omega \rho(\bxi\cdot \nabla)\bv\cdot\bv\dx\d s-\int_0^t\int_\Omega \rho(\bv\cdot \nabla)\bxi\cdot\bv\dx\d s+\int_0^t\int_\Omega{\rho}(\bv\cdot \nabla )\bv\cdot \bxi\dx\d s\\&=\int_0^t\int_\om (\nabla\rho\cdot\bv)\bv\cdot\bxi\dx\d s+\int_0^t\int_\Omega{\rho}(\bv\cdot \nabla )\bv\cdot \bxi\dx\d s\\&
	=\int_0^t\int_\om \Div(\rho\bv\otimes\bv)\cdot\bxi\dx \d s
	=-\int_0^t\int_\om \rho\bv\otimes \bv:\nabla \bxi\dx\d s,
\end{align*}
and then concluding by a density argument. 

Finally, notice that we can repeat the same estimate as in \eqref{dtest2}, this time in the limit equation as $h\to0$, and see from \eqref{velocity2k} and \eqref{identityk} that we can choose $\bxi\in \Hus\cap \bW^{1,\infty}(\om)$, and obtain, since $\bv\in L^\infty(0,\Ts;\bL^2_\sigma(\om))$ for any $\Ts>0$,
\begin{align}
	&\nonumber\normmm{\langle \pt\bP_\sigma(\rho\bv), \bxi\rangle_{(\Hus\cap \bW^{1,\infty(\om)})',\Hus\cap \bW^{1,\infty}(\om)}} \\&\nonumber
	\leq\nonumber
	C\norm{\bv}^2\norm{\bxi}_{\bW^{1,\infty}(\om)}
	+C\norm{\bv}\norm{\nabla u}\norm{\bxi}_{\bW^{1,\infty}(\om)}+C\norm{\bv}_{\Hus}\norm{\bxi}_{\Hus}+\norm{\xhi}_{L^\infty(\om)}\norm{\nabla \tw}\norm{\bxi}\\&
	\leq C(\Ts)\left(\norm{\bv}_{\Hus}+\norm{\nabla u}+\norm{\nabla \wtilde{w}}\right)\norm{\bxi}_{\Hus\cap \bW^{1,\infty}(\om)},\label{dtest2bk}
\end{align}
for almost any $t>0$.  Since $\wtilde{w},u\in L^2(0,\Ts;\Hz)$ and $\bv\in L^2(0,\Ts;\Hus)$, this entails

\begin{align}
	&\norm{\pt\bP_\sigma (\rho\bv)}_{L^2(0,\Ts;(\Hus\cap \bW^{1,\infty}(\om))')}\leq C(\Ts).\label{regv1b1kk}
\end{align}
\subsection{$(\xhi,\mu,\bv)$ is a varifold solution to Navier-Stokes-Mullins-Sekerka system}\label{ssf}
Checking that, for any $\Ts>0$, the triple $(\xhi,\bv,\mu)$ is a varifold solution according to Definition \ref{weaksol} is now immediate. Indeed, properties (1)-(3) have been assessed in the previous sections. Concerning Property (4), the potential $\wtilde{w}\in L^2(0,T;\Hz)$ is already a good candidate to be the curvature potential for the mean curvature, thanks to \eqref{variation4}. Nevertheless, we also need that the curvature potential belongs to $\mathcal G_{\xhi(\cdot,t)}$ for almost any $t>0$, where this space is defined in \eqref{Gchi}. Now, from \eqref{variation4} we deduce that 
\begin{align}
	&	\delta\mu_t(B)=
	\int_\Omega \chi(\cdot,t)\text{div}(\wtilde{w}(\cdot,t)B)\d x=-\langle B\cdot \nabla\xhi(\cdot,t),\wtilde{w}(t)\rangle_{\Hmz,\Hz},\label{variationFkk1}
\end{align} 
for all $B\in \mathcal S_{\xhi(\cdot,t)}$ and almost any $t>0$. Here the map $B\cdot \nabla \xhi(\cdot,t)\mapsto \delta\mu_t(B)$ is linear and continuous from $\mathcal W_{\xhi(\cdot,t)}$ (see definition \eqref{Wchi}) to $\R$, so that it can be uniquely extended by density as a continuous linear functional $L_t\in\mathcal V_{\xhi(\cdot,t)}'\subset \Hz$ with, from \eqref{variationFkk1}, 
\begin{align}
	\norm{L_t}_{\mathcal V_{\xhi(\cdot,t)}'}\leq \norm{\nabla\wtilde{w}(t)},
	\label{normcontrol1}
\end{align} 
as well as 
\begin{align}
	\langle L_t, B\cdot\nabla \xhi(\cdot,t)\rangle_{\mathcal V_{\xhi(\cdot,t)}', \mathcal V_{\xhi(\cdot,t)}}=\delta\mu_t(B),\quad \forall B\in \mathcal S_{\xhi(\cdot,t)},
	\label{extension}
\end{align}
for almost any $t>0$. Now, since $\mathcal V_{\xhi(\cdot,t)}$ is a Hilbert space, by Riesz Lemma there exists a unique $z(t)\in \mathcal V_{\xhi(\cdot,t)}$ such that  
\begin{align*}
	\langle L_t, f\rangle_{\mathcal V_{\xhi(\cdot,t)}', \mathcal V_{\xhi(\cdot,t)}}=(z(t), f)_{\mathcal V_{\xhi(\cdot,t)}}=(z(t), f)_{\Hmz},\quad \forall f\in \mathcal V_{\xhi(\cdot,t)}.
\end{align*}
Recalling that $\mathcal G_{\xhi(\cdot,t)}=-(-\Delta_N)^{-1}\mathcal V_{\xhi(\cdot,t)}$, together with the property that $-\Delta_N$ is the Riesz isomorphism between $\Hz$ and $\Hmz$, we deduce that $w_0(t):=-(-\Delta_N)^{-1}z(t)\in \mathcal G_{\xhi(\cdot,t)}	\subset\Hz$ is such that 
\begin{align*}
	\langle L_t, f\rangle_{\mathcal V_{\xhi(\cdot,t)}', \mathcal V_{\xhi(\cdot,t)}}=\langle (-\Delta_N)^{-1}z(t), f\rangle_{\Hz,\Hmz}=-\langle w_0(t), f\rangle_{\Hz,\Hmz},\quad \forall f\in \mathcal V_{\xhi(\cdot,t)}.
\end{align*}
Therefore, we can conclude from \eqref{extension} that 
\begin{align}
	\delta\mu_t(B)	=\langle L_t, B\cdot\nabla\xhi\rangle_{\mathcal V_{\xhi(\cdot,t)}', \mathcal V_{\xhi(\cdot,t)}}=-\langle w_0(t),B\cdot \nabla\chi(\cdot,t)\rangle_{\Hz,\Hmz},\quad \forall B\in \mathcal S_{\xhi(\cdot,t)},\label{deltamuf}
\end{align} 
and, from \eqref{normcontrol1} and the Riesz isometry properties, that 
\begin{align}
	\norm{\nabla w_0(t)}=\norm{w_0(t)}_{\Hz}= \norm{z(t)}_{\Hmz}=\norm{z(t)}_{\mathcal V_{\xhi(\cdot,t)}}=\norm{L_t}_{\mathcal V_{\xhi(\cdot,t)}'}\leq \norm{\nabla \wtilde{w}(t)},
	\label{encontrol}
\end{align}
for almost any $t>0$. Therefore, since $\wtilde w\in L^2(0,\Ts;\Hz)$, for any $\Ts>0$, also $w_0\in L^2(0,\Ts;\Hz)$, with $w_0(t)\in \mathcal G_{\xhi(\cdot,t)}$ for almost any $t>0$. Then, by means of \eqref{deltamuf}, we can see that 
\begin{align}
	&	\delta\mu_t(B)=
	\int_\Omega \chi(\cdot,t)\text{div}({w}_0(\cdot,t)B)\d x, \quad \forall B\in \mathcal S_{\xhi(\cdot,t)}.\label{variationFkk2}
\end{align} 
This entails that we can apply Lemma \ref{lemma9}, so that, since $w_0\in L^2(0,\Ts;\Hz)$, for any $\Ts>0$, there exists $\lambda \in L^2_{loc}([0,\infty))$ such that 
$$
w:=w_0+\lambda \in L^2(0,\Ts; H^1(\om)), \quad \forall \Ts>0,
$$
satisfies \eqref{GThom} and \eqref{C1}. This means that the entire Property (5) of Definition \ref{weaksol} (together with \eqref{C1}) is satisfied by $w$, which is then the chemical potential associated to the mean curvature. Moreover, Property (6) of the same definition is obtained from \eqref{velocity2k}, recalling \eqref{identityk}. Indeed, since $\bxi(t)\in \mathcal S_{\xhi(\cdot,t)}$ for any $\bxi\in C^\infty_c([0,\Ts);\mathbf C_{c,\sigma}^\infty(\om))$ and almost any $t\geq0$, it holds, from \eqref{variationFkk1} and \eqref{deltamuf},
\begin{align*}
	&\int_\om \xhi(\cdot,t)\nabla \wtilde{w}(t)\cdot \bxi(t)\dx=-\langle \bxi(t)\cdot \nabla\xhi(\cdot,t),\wtilde{w}(t)\rangle_{\Hmz,\Hz}\\&=-\langle w_0(t),\bxi(t)\cdot \nabla\chi(\cdot,t)\rangle_{\Hz,\Hmz}=\int_\om \xhi(\cdot,t)\nabla {w}(t)\cdot \bxi(t)\dx,
\end{align*}
so that we can substitute $\wtilde{w}$ with $w$ in \eqref{velocity2k} and infer the validity of \eqref{velocity}. This also entails, by density, that 
\begin{align*}
	&\int_\om \xhi(\cdot,t)\nabla \wtilde{w}(t)\cdot \bxi(t)\dx=\int_\om \xhi(\cdot,t)\nabla {w}(t)\cdot \bxi(t)\dx,
\end{align*}
for any $\bxi\in L^2(0,\Ts;\bL^2_\sigma(\om))$, almost any $t\in(0,\Ts)$, and any $\Ts>0$. Therefore, by plugging $\bxi=\bv$, we also infer 
\begin{align}
	&\int_\om \xhi(\cdot,t)\nabla \wtilde{w}(t)\cdot \bv(t)\dx=\int_\om \xhi(\cdot,t)\nabla {w}(t)\cdot \bv(t)\dx,\label{ff1}
\end{align}
for almost any $t>0$. 

As a consequence, due to \eqref{encontrol} and \eqref{ff1}, we can substitute $\wtilde w$ with $w$ in \eqref{energy1abc3b} and \eqref{energyv2b1}, so that the energy estimates of Property (7) in Definition \ref{weaksol} hold.

Therefore  $(\xhi,\bv,\mu)$, with corresponding potentials $w$ and $u$, is a varifold solution to the Navier-Stokes-Mullins-Sekerka system according to Definition \ref{weaksol}, for any $\Ts>0$, and this concludes the proof of Theorem \ref{thm1}.

\section{Proof of Theorem \ref{consistency}}
\label{proofthm0}
\subsection{Step 1. Classical solutions are varifold solutions}
Let $(\cA,\bv)$ be a classical solution for the Navier-Stokes-Mullins-Sekerka system in the sense of \eqref{A1}-\eqref{A10}. First, we define the associated varifold as $\normmm{\mu_t^\om}_{\mathbb S^{d-1}}:=\normmm{\nabla\xhi(\cdot,t)}=\cH^{d-1}\llcorner(\pAs(t)\cap \om)$, $\normmm{\mu_t^{\pom}}_{\mathbb S^{d-1}}:=\xhi(\cdot,t)\cH^{d-1}\llcorner\pom=\cH^{d-1}\llcorner (\pAs(t)\cap \pom)$, with $\bn_{\pcA(t)}=\frac{\nabla\xhi(\cdot,t)}{\normmm{\nabla\xhi(\cdot,t)}}$. Therefore, we simply have that $E[\mu_t]=E[\xhi(\cdot,t)]$ as defined in \eqref{E1}. Then, admissibility of the couple $(\xhi,\mu)$ according to Definition \ref{admissibility} immediately comes from the regularity of the classical solution. Now, as in \cite[Eq.(155)]{SH}, we can deduce that properties \eqref{PH}-\eqref{PK} hold with $H_{\xhi(\cdot,t)}=\bH_{\pcA(t)}\cdot\bn_{\pcA(t)}$, where the stated integrability is a consequence of the boundary condition \eqref{A3} and a standard trace estimate for $\ovu$. Then, properties (2)-(3) are immediately satisfied from the regularity of the classical solution.  Analogously, $u=\ovu$ is the kinetic potential due to \eqref{A2}, verifying Property (5). Also, \eqref{A6}-\eqref{A10} allow to easily verify the validity of \eqref{velocity}, by recalling \eqref{byparts1} and the trace identity \eqref{A3}. Concerning Property (7), the energy estimate \eqref{NS} is a direct consequence of \eqref{enid1A}. Then also also \eqref{advectiveMullins} holds. Indeed, recalling \eqref{enid2A}, and since $\bv\equiv \mathbf 0$ along $\pom$ and $\Div\bv=0$ in $\om$, we have
\begin{align}
	&\nonumber\ddt E[	\xhi(\cdot,t)]=-\int_\om\normmm{\nabla\ovu(t)}^2\dx-\int_\om \bv(t)\cdot\bn_{\pcA(t)}\ovu(\cdot,t)\d\cH^{d-1} \\&
	=-\int_\om\normmm{\nabla\ovu(t)}^2\dx+\int_\om \xhi(\cdot,t)\bv(t)\cdot\nabla \ovu(\cdot,t)\dx.
	\label{energycl}
\end{align} 
Since then $u=\ovu$ and $\nabla w=\nabla \ovu$, we obtain \eqref{advectiveMullins} from \eqref{energycl}. In conclusion, note that, again from \cite[Eq.(155)]{SH}, we see that, after an integration by parts, 
\begin{align*}
	&\delta\mu_t(B)=\delta E[\xhi(\cdot,t)](B)=-c_0\int_{\pAs(t)\cap \om} \bH_{\pcA(t)}\cdot B\d\cH^{d-1}\\&=-\int_{\pAs(t)\cap \om} \ovu(\cdot,t)\bn_{\pcA(t)}\cdot B\d\cH^{d-1}
	=\int_\om\chi(\cdot,t)\Div(\ovu(\cdot,t)B)\dx,
\end{align*}    
for any $B\in C^1(\overline{\om};\R^d)$, with $B\cdot\bn_\pom\equiv 0$ along $\pom$. Therefore, we can set $\tilde{w}_0(t)=\ovu-\tfrac{\int_\om\ovu(t)\dx}{\mathcal L^d(\om)}$ and, as in the final part of the proof of Theorem \ref{thm1} (see Section \ref{ssf}), obtain the existence of $w_0(t)\in \mathcal G_{\xhi(\cdot,t)}$ such that 
\begin{align*}
	\delta\mu_t(B)=\int_\om\chi(\cdot,t)\Div(w_0(\cdot,t) B)\dx
\end{align*}
for $B\in \mathcal S_{\xhi(\cdot,t)}$. Then we apply Lemma \ref{lemma9} to infer that there exists $\lambda\in L^2(0,\Ts)$, such that $w=w_0+\lambda$ satisfies the conditions  \eqref{GThom} as well as \eqref{C1}. Since, as in \eqref{encontrol}, $\norm{\nabla w_0(t)}\leq \norm{\nabla \tilde w_0}=\norm{\nabla \ovu(t)}$, we can substitute $w$ in place of $\ovu$ in properties (6)-(7), so that $w$ is actually the curvature potential satisfying Property (4). Therefore, $(\cA,\bv)$, with potentials $u$ and $w$, is also a varifold solution on $(0,\Ts)$ according to Definition \ref{weaksol}.
\subsection{Step 2. Smooth varifold solutions are classical solutions}
Let $(\xhi,\mu,\bv)$ be a varifold solution with smooth interface, i.e., satisfying Definition \ref{admissibility} and such that $\chi$ is represented as $\chi(x,t) = \xhi_{\cA(t)}(x)$, where $\cA = (\cA (t))_t\in[0,\Ts)$
is a time-dependent family of smoothly evolving sets $\cA (t) \subset \om$, $t \in [0, \Ts)$, exactly as in the previous step. Also, assume that the velocity $\bv$ is smooth. Fix $t\in(0,\Ts)$ such that \eqref{wkk1} holds. Then, following word by word \cite[Step 2, proof of Lemma 4]{SH}, we can deduce that $\bH_{\pcA(t)}\cdot\bn_{\pcA(t)}=H_{\xhi(\cdot,t)},\quad \cH^{d-1}\text{-a.e. on }\pcA(t)\cap \om$, so that, from \eqref{wkk1}, we deduce 
\begin{align}
	w(\cdot,t)\bn_{\pcA(t)}=c_0\bH_{\pcA(t)},\quad \cH^{d-1}\text{-a.e. on }\pcA(t)\cap \om.\label{curvatureovh}
\end{align}
Also, recalling the characterization \eqref{pchi}, since $w(t)\subset \mathcal G_{\chi(\cdot,t)}$, we infer that $w$ is harmonic in $\cA$ and in the interior of $\om\setminus \cA$. By elliptic regularity theory (as in \cite{SH}), we can further obtain a continuous representative for $w$ and conclude that \eqref{curvatureovh} holds everywhere on $\pcA(t)\cap \om$, so that \eqref{A1} and \eqref{A3} hold.
Concerning the contact angle condition \eqref{A5}, this can be retrieved as in \cite[Step 2, proof of Lemma 4]{SH}. Here we only point out that the assumption that $\delta\mu_t(B)$ can be substituted by $\delta E[\xhi(\cdot,t)](B)$ in \eqref{PH} is essential. 

 Now, by \eqref{kinetic}, it holds that, for any $t\in(0,\Ts)$, the kinetic potential $u$ satisfies 
$$
\pt\xhi(\cdot,t)+(\bv\cdot\bn_{\pcA(t)})\bn_{\pcA(t)}\cdot\frac{\nabla\xhi(\cdot,t)}{\normmm{\nabla\xhi(\cdot,t)}}=(-V_{\pcA(t)}+(\bv\cdot\bn_{\pcA(t)})\bn_{\pcA(t)})\cdot\frac{\nabla\xhi(\cdot,t)}{\normmm{\nabla\xhi(\cdot,t)}}=-A_N u.
$$
Namely, thanks to the regularity of $\text{supp}\normmm{\nabla\xhi{(\cdot,t)}}=\pcA(t)\cap \om$, 
	\begin{align}
		\Delta u(\cdot,t)=0,\quad\text{ in }\om\setminus\pcA(t).	\label{potu}
	\end{align}
Again, by a very similar argument to the one for obtaining \cite[Eq.(165)]{SH} (i.e., elliptic regularity theory for the Neumann problem, and integration by parts), we infer 
\begin{align}
\label{CA}	&V_{\pcA(t)}=(\bv\cdot\bn_{\pcA(t)})\bn_{\pcA(t)}-(\bn_{\pcA}\cdot[[\nabla u(\cdot,t)]])\bn_{\pcA(t)},\quad\text{ on }\pcA(t)\cap \om,\\&(\bn_{\pcA(t)}\cdot\nabla )u(\cdot,t)=0,\quad \text{ on }\pom\setminus \overline{\pcA(t)\cap \om}.
\label{CB}
\end{align}
At this point, since, following \cite{SH}, the higher-order Sobolev regularity of $w,u$ can be obtained, then, after an integration by parts, also exploiting the regularity of $\bv$, we can easily show the validity of \eqref{A6}-\eqref{A10}, this time with $\widetilde{\mathbf J}=(\rho_1-\rho_2)\nabla u$. Notice that the validity of \eqref{A9} can be shown thanks to \eqref{curvatureovh}. To complete the proof of the validity of \eqref{A1}-\eqref{A10}, we thus only need to prove that $\nabla u=\nabla w$ in $\overline\om$. We show that 
\begin{align}
	\norm{\nabla u(t)-\nabla w(t)}=0,\quad \forall t\in[0,\Ts),\label{id1}
\end{align}
which entails, by the regularity of $u,w$, that $\nabla u=\nabla w$ in $\overline\om$, allowing to conclude the proof of \eqref{A1}-\eqref{A10}. To this aim, note that we can obtain, arguing as in \eqref{enid2A}, and recalling  \eqref{curvatureovh}, \eqref{CA}, and \eqref{CB},
\begin{align*}
	\ddt E[\chi(t)]&=-\int_{\pcA(t)\cap \om}V_{\pcA(t)\cap \om}\cdot c_0\bH_{\pcA(t)}\d\cH^{d-1}\\&=-\int_{\pcA(t)\cap \om}(\bv(t)\cdot\bn_{\pcA(t)})w(\cdot,t)\d\cH^{d-1}-\int_\om \nabla u(\cdot,t)\cdot \nabla w(\cdot,t)\dx\\&
	=\int_\om \xhi(\cdot,t)\bv(t)\cdot\nabla w(t)\dx -\int_\om \nabla u(\cdot,t)\cdot \nabla w(\cdot,t)\dx.
\end{align*}
Subtracting this identity (in integrated form) from \eqref{advectiveMullins}, and recalling that $E[\xhi(\cdot,T)]\leq E[\mu_T]$ (note that $(\xhi,\mu)$ is admissible according to Definition \ref{admissibility}), we immediately infer
\begin{align}
0\leq 	\frac12 \int_0^T\norm{\nabla w(t)-\nabla u(t)}^2\d t \leq E[\xhi(\cdot,T)]-E[\mu_T]\leq 0,\label{muid1}
\end{align}
for almost any $T\in(0,\Ts)$, which gives \eqref{id1}.  From \eqref{muid1} we also deduce that $E[\chi(·, T )] = E[\mu_T]$ for almost any $T\in(0, \Ts)$. This identity entails, exploiting \eqref{PC} (in the form $\normmm{\nabla \xhi(\cdot,t)}\leq\tmu\llcorner \om$) and \eqref{PD}, that $\normmm{\nabla \xhi(\cdot,t)}=\tmu\llcorner \om$, which gives the first identity in \eqref{H1}, and $(\cos\gamma)\xhi(\cdot,t)\cH^{d-1}\llcorner\pom=(\cos\gamma)(\tmu\llcorner \pom+\tmup)$. Since $\mu_t\in \mathcal M(\overline{\om}\times \mathbb S^{d-1})$ is an integer
rectifiable oriented varifold, by considering the Radon-Nikodym derivative of both sides of the latter identity with respect to $\cH^{d-1}\llcorner\pom$, we find that
$$
(\cos\gamma )\chi(\cdot,t)=(\cos\gamma)\frac{\tmu\llcorner \pom+\tmup}{\cH^{d-1}\llcorner\pom}(\cdot,t)=(\cos \gamma)\frac{\tmup}{\cH^{d-1}\llcorner\pom}(\cdot,t)+m(\cdot,t),
$$
with $m(\cdot,t):\pom \to \mathbb N\cup \{0\}$, so that, since $\xhi(\cdot,t)\in\{0,1\}$ and $\frac{\tmup}{\cH^{d-1}\llcorner\pom}(\cdot,t)\geq0$, necessarily $m\equiv0$. This entails the second identity of \eqref{H1}, as well as \eqref{H2}. The proof is thus concluded.

\bigskip\bigskip

\textbf{Acknowledgments.} The authors want to warmly thank Tim Laux for helpful and inspiring discussions, especially concerning the results in Appendix \ref{measurability}. Part of this contribution was completed while AP was visiting HA at the Faculty of Mathematics of the University of Regensburg, whose hospitality is kindly acknowledged. AP also gratefully aknowledges support from the Alexander von Humboldt Foundation for the stay in Regensburg. This research was funded in part by the Austrian Science Fund (FWF) \href{https://doi.org/10.55776/ESP552}{10.55776/ESP552}.
AP is a member of Gruppo Nazionale per l’Analisi Matematica, la Probabilità e le loro Applicazioni (GNAMPA) of
Istituto Nazionale per l’Alta Matematica (INdAM).  
For open access purposes, the authors have applied a CC BY public copyright license to
any author accepted manuscript version arising from this submission. 
\\

\appendix
\section{Energy identity for the Navier-Stokes-Mullins-Sekerka system}
\label{App:A}
In this appendix we give more details for the computations of some identities for the Navier-Stokes-Mullins-Sekerka system \eqref{A1}-\eqref{A10}, which have been shown without justification in the introduction. First, note that the following integration by parts formula holds for any sufficiently regular functions $\mathbf f:\overline{\Omega}\to \R^d$, $\eta:\overline{\om}\to \R$:
	\begin{align}
		-\int_{\om\setminus\pcA(t)}\eta\Div \mathbf f\dx=\int_\om \nabla\eta\cdot\mathbf f\dx +\int_{\pcA(t)\cap\om}\eta(\bn_{\pcA(t)}\cdot[[\mathbf f]])\d \mathcal H^{d-1}+\int_\pom \eta(\bn_{\pom}\cdot\mathbf f)\d \mathcal H^{d-1}.
			\label{byparts1}
	\end{align}

Then, in the same setting and with the same notation as in the introduction, it is immediate to notice that the total mass is conserved in the flow, since, recalling \eqref{A7},
	\begin{align*}
	&	\ddt \mathcal L^d(\cA(t))=\ddt \mathcal L^d(\cA(t)\cap 	\om)=-\int_{\pcA(t)\cap \om}V_{\pcA(t)}\cdot\bn_{\pcA(t)}\d \mathcal H^{d-1}\\&=-\int_{\pcA(t)\cap \om}(\bv\cdot\bn_{\pcA})\d \mathcal H^{d-1}+\int_{\pcA(t)\cap \om}(\bn_{\pcA(t)}\cdot[[\nabla\ovu]])\d \mathcal H^{d-1}=0,
	\end{align*}
which justifies \eqref{massflow}. To obtain  \eqref{enid2A0}, let us denote by $\tau_{\pcA(t)\cap \om}$ a vector field on $\partial(\overline{\pcA(t)\cap \pom})\subset \pom$, tangent to $\pcA(t)\cap \om$, normal to $\partial(\overline{\pcA(t)\cap \pom})$, and	pointing outward with respect to $\pcA(t)\cap \om$ (i.e., it is the conormal vector to $\pcA(t)\cap \om$). Also, denote by $\tau_{\pcA(t)\cap \pom}$ a vector field again over  $\partial(\overline{\pcA(t)\cap \pom})\subset \pom$, but tangent to $\pcA(t)\cap \pom$, normal to $\partial(\overline{\pcA(t)\cap \pom})$ and pointing towards $\pcA(t)\cap \pom$ (i.e., it is the inner conormal vector to $\pcA(t)\cap \pom$). Notice that, by \eqref{A5}, $\tau_{\pcA(t)\cap \om}\cdot \tau_{\pcA(t)\cap \pom}=\bn_{\pcA(t)}\cdot \bn_{\pom}=\cos \gamma$.
	Then, since $V_{\pcA(t)}\cdot \bn_{\pom}=0$ and $V_{\pcA(t)}=(V_{\pcA(t)}\cdot \tau_{\pcA(t)\cap \pom})\tau_{\pcA(t)\cap \pom}$ (indeed $V_{\pcA(t)}$ is normal to $\partial(\overline{\pcA(t)\cap \pom})$), recalling \eqref{byparts1}, and the fact that the contact angle $\gamma$ does not depend on time, it holds by standard computations and integrations by parts,
	\begin{align}
		&\nonumber\ddt E[\xhi(t)]
		=\nonumber c_0\int_{\pcA(t)\cap \om}\Div_{\pcA(t)\cap \Omega}(V_{\pcA(t)})\d\cH^{d-1}
		+c_0\cos\gamma\int_{\pcA(t)\cap \pom}\Div_{\pcA(t)\cap \pom}(V_{\pcA(t)})\d\cH^{d-1}
		\\&\nonumber
		=-\int_{\pcA(t)\cap \om}V_{\pcA(t)}\cdot c_0\bH_{\pcA(t)}\d \mathcal H^{d-1}+c_0\int_{\partial(\pcA(t)\cap \om)}V_{\pcA}\cdot (\tau_{\pcA(t)\cap \om}-\cos\gamma\tau_{\pcA(t)\cap \pom})\d\cH^{d-2}
		\\&\nonumber=-\int_{\pcA(t)\cap \om}V_{\pcA(t)}\cdot c_0\bH_{\pcA(t)}\d \mathcal H^{d-1}\nonumber
		=-\int_{\pcA(t)\cap \om }(\bv(t)\cdot\bn_{\pcA(t)})\bn_{\pcA(t)}\cdot c_0\bH_{\pcA(t)}\d \mathcal H^{d-1}\\&\nonumber\quad+\int_{\pcA(t)\cap \om }(\bn_{\pcA(t)}\cdot[[\nabla\ovu(\cdot,t)]])\bn_{\pcA(t)}\cdot c_0\bH_{\pcA(t)}\d \mathcal H^{d-1}\\&
		=-\int_{\pcA(t)\cap \om }(\bv(t)\cdot\bn_{\pcA(t)})\ovu(\cdot,t)\d \mathcal H^{d-1}-\int_{\om}\normmm{\nabla\ovu(\cdot,t)}^2\dx,\label{enid2A}
	\end{align}
    which gives \eqref{enid2A0}. In conclusion, identity \eqref{enid1A0} can be computed from the fluid velocity equations \eqref{A6}-\eqref{A10}, together with \eqref{A3}:  
	\begin{align}
		&\nonumber\frac12\ddt \int_\om \rho(\xhi(t))\normmm{\bv(t)}^2\dx +\int_\om \nu(\xhi(t))\normmm{D\bv(t)}^2\dx\\&=\int_{\pcA(t)\cap \om }c_0\bH_{\pcA(t)}\cdot \bv\d \mathcal H^{d-1}=\int_{\pcA(t)\cap \om }(\bv(t)\cdot\bn_{\pcA(t)})\ovu(\cdot,t)\d \mathcal H^{d-1},\label{enid1A}
	\end{align}
	where we have also used \eqref{densityA}.

\section{Bochner measurability of $\xhi_t^{T_h(t)}$ with values in $L^p(\Omega)$, $p\geq1$}
\label{measurability}
\subsection{Measurability}
Concerning the proof of Theorem \ref{thm1}, referring to Section \ref{firstminmovement}, in order to prove that we can choose the approximating function $\xhi_t^{T_h(t)}$ to be measurable in time in a suitable sense, we resort to the Kuratowski-Ryll-Nardzewski Measurable Selection Theorem. For the sake of simplicity, we concentrate on the interval $(0,h]$, $h>0$. We set $\mathcal F:=\mathcal B((0,h])$, i.e., the Borel $\sigma$-algebra on $(0,h]$. 
Now we introduce the (possibly) set-valued function
\begin{align*}
	\psi(\omega):=\arg\min_{\xhi\in \mathcal M_{m_0}} \left(E[\chi]+\frac 1{2\omega}\normh{\chi- \xhi_0\circ X_{-\omega}^{\bu_0^k}}^2\right),
		\quad \text{ for }\omega\in (0,h],
\end{align*}
for $k\in \mathbb N$ fixed. This set is always nonempty by means of an application of the direct method of Calculus of Variations. Let us now show that $\psi(\omega)$ is closed with respect to the strong $L^1(\Omega)$ topology. Indeed, assume that $\{\xhi_n\}_n\subset \psi(\omega)$ is such that $\xhi_n\to \xhi$ strongly in $L^1(\Omega)$.  Moreover, by the optimality of any element $x\in \psi(\omega)$, it holds
\begin{align}
E[x]+\frac 1{2\omega}\normh{x- \xhi_0\circ X_{-\omega}^{\bu_0^k}}^2\leq E[\xhi_0]+\frac 1{2\omega}\normh{\xhi_0- \xhi_0\circ X_{-\omega}^{\bu_0^k}}^2 \leq C_0,\label{pp1}
\end{align}
since $\frac 1{2\omega}\normh{\xhi_0- \xhi_0\circ X_{-\omega}^{\bu_0^k}}^2\to 0$ as $\omega\to 0$ by Lemma \ref{basic1}. Therefore, since the sequence $\xhi_n$ has bounded total variation of its distributional gradients, by \cite[Proposition 3.13]{AmbrosioFuscoPallara} this also implies that $\xhi_n\to \xhi$ weakly* in $BV(\om;\{0,1\})$. Therefore, it holds, by lower semicontinuity of $E[\xhi]$ with respect to weak* convergence in $BV(\om;\{0,1\})$ (see for instance \cite[Proposition 1.2]{Modica}),
\begin{align}
	E[\xhi]\leq \liminf_{n\to \infty}E[\xhi_n],\label{liminf}
\end{align}
and also note that, for any $\phi\in \Hz$, since $\{\xhi_n\}_n\subset \mathcal M_{m_0}$,
\begin{align*}
&\normmm{({\xhi_n- \xhi_0\circ X_{-\omega}^{\bu_0^k}}-(\xhi-\xhi_0\circ X_{-\omega}^{\bu_0^k}),\phi)}\\&=(\xhi_n-\xhi,\phi)\leq \norm{\xhi_n-\xhi}\norm{\phi}\leq \sqrt 2\norm{\chi_n-\xhi}_{L^1(\Omega)}^\onehalf\norm{\phi}_{\Hz},
\end{align*}
entailing that 
\begin{align*}
	\frac 1{2\omega}\normh{\xhi_n- \xhi_0\circ X_{-\omega}^{\bu_0^k}}^2\to 	\frac 1{2\omega}\normh{\xhi-\xhi_0\circ X_{-\omega}^{\bu_0^k}}^2,\text{ as }n\to \infty.
\end{align*}
Therefore, this result together with \eqref{liminf} allows to conclude, passing to the limit as $n\to \infty$, that 
$$
E[\xhi]+	\frac 1{2\omega}\normh{\xhi- \xhi_0\circ X_{-\omega}^{\bu_0^k}}^2\leq \liminf_{n\to \infty}\left(E[\xhi_n]+	\frac 1{2\omega}\normh{\xhi_n- \xhi_0\circ X_{-\omega}^{\bu_0^k}}^2\right),
$$
thus entailing that $\xhi\in \psi(\omega)$ and then that $\psi(\omega)$ is closed in the strong $L^1(0,h)$ topology for any $\omega\in (0,h]$.
We now define $X:=L^1(0,h)$, and we endow it with its strong topology. Since $X$ is separable, $X$ is a Polish space. We also denote by $\mathcal{B}(X)$ the corresponding Borel $\sigma$-algebra.
Notice that, since $\Mm\subset X$,
$$
\psi(\omega)\subset X,\quad \forall \omega\in (0,h],
$$
and, as observed, $\psi(\omega)$ is a nonempty closed subset of $X$ for any $\omega\in (0,h]$.
The Kuratowski-Ryll-Nardzewski Measurable Selection Theorem then states that if, for any $U\in \mathcal B(X)$, it holds 
\begin{align}
	\left\{\omega\in (0,h]:\ \psi(\omega)\cap U\not =\emptyset\right\} \in\mathcal F,\label{condition}
\end{align}
then $\psi$ admits a selection $\xhi_\omega:(0,h]\to X$ which is $\mathcal F$-$\mathcal B(X)$ measurable (i.e., for any $B\in \mathcal B(X)$ it holds $(\xhi_\omega)^{-1}(B)\in \mathcal F$). This concludes the proof if we set $\xhi_t^{T_h(t)}:=\xhi_t$ for  any $t\in(0,h]$ (here we denoted $\omega$ by $t$). We are then left to show \eqref{condition}. First, we introduce the function
\begin{align}
	f(\omega):=E[x]+\frac 1{2\omega}\normh{x- \xhi_0\circ X_{-\omega}^{\bu_0^k}}^2,\quad x\in \psi(\omega),\quad \forall \omega\in (0,h],
\end{align}
which is well defined since $f(\omega)$ is exactly the evaluation of the functional which is minimized by any element $x\in \psi(\omega)$, and thus it \textit{does not} depend on $x\in \psi(\omega)$. Now, following the same argument as to obtain \eqref{derivative}, we infer that $f$ is locally Lipschitz continuous in $(0,h]$. In particular, it is continuous in $(0,h]$. Recalling \eqref{pp1}, we can now define the following functional $G:X\times (0,h]\to \mathbf R^+\cup\{+\infty\}$ as
\begin{align*}
	G(\chi,\omega):=\begin{cases} E[\chi]+\frac 1{2\omega}\normh{\chi- \xhi_0\circ X_{-\omega}^{\bu_0^k}}^2-f(\omega),\ \text{if }\xhi\in \Mm\cap\{x\in \Mm:\ \normmm{\nabla x}(\Omega)\leq C_0\},\\
		+\infty, \qquad\qquad \qquad\quad\qquad\qquad\qquad\text{if }\xhi\in X\setminus (\Mm\cap\{x\in \Mm:\ \normmm{\nabla x}(\Omega)\leq C_0\}).
		\end{cases}
\end{align*}
Since $\Mm\cap\{x\in \Mm:\ \normmm{\nabla x}(\Omega)\leq C_0\}$ is a closed subset of $X$ (which can be proven by the same argument as to show that $\psi(\omega)$ is closed in $X$), it is immediate to see that $G$ is $\mathcal B(X)$ measurable with respect to the first variable for any $\omega\in (0,h]$, and, from the discussion above on $f$, continuous in its second variable $\omega$ for any $\xhi\in X$, i.e., $G$ is a Charathéodory functional. Therefore, it is jointly $\mathcal B(X)\otimes \mathcal F$ measurable. To conclude the proof we are left to observe that, thanks to \eqref{pp1}, given $U\in \mathcal B(X)$,
\begin{align*}
		\left\{\omega\in(0,h]:\ \psi(\omega)\cap U\not =\emptyset\right\}=\Pi_{\omega}\left(G^{-1}(0)\cap (U\times (0,h])\right),
\end{align*}
where $\Pi_\omega$ is the projection of a subset of $X\times (0,h]$ over its second component, which is continuous with continuous inverse. Since $G$ is Charathéodory, this implies $\Pi_{\omega}\left(G^{-1}(0)\cap (U\times (0,h])\right)\in \mathcal F$, entailing \ref{condition} and concluding the proof.
\subsection{Bochner measurability}
In the previous subsection we have shown that  $\psi$ admits a selection $\xhi_\omega$ which is $\mathcal F$-$\mathcal B(X)$ measurable, and thus that we can find $\xhi_t^{Th(t)}$ which is $\mathcal F$-$\mathcal B(X)$ measurable. We now aim at showing that this selection is also (strongly) Bochner-measurable with values in $L^p(\Omega)$ for any $p\in[1,\infty)$. First, we immediately see that our selection $\xhi_t^{Th(t)}$ is weakly measurable. Indeed, for any element $g\in X'$ we have that $g(\xhi_\cdot^{T_h(\cdot)})$ is $\mathcal F$-measurable, since it is the composition of a continuous functional on $X$ with an $\mathcal F$-$\mathcal B(X)$ measurable function. Then, since $X$ is separable, by Pettis' Theorem this implies that $\xhi_t^{T_h(t)}$ is also strongly Bochner-measurable with values in $L^1(\Omega)$. In conclusion, since $\norm{\xhi_t^{T_h(t)}}_{L^\infty(\Omega)}\leq 1$ for almost any $t\in (0,h]$, we also have that $\xhi_t^{T_h(t)}$ is also (strongly) Bochner-measurable with values in $L^p(\Omega)$ for any $p\geq1$, concluding the argument.
\section{On an embedding in vector-valued continuous spaces}
\label{App:C}
Here we prove the following embedding, which has been used to derive some regularity properties of the phase variable $\xhi$. Note that this result is more refined than the calssical Aubin-Lions type argument, due to the weak* measurability with respect to the first space. In particular, we have the following 
\begin{lemma}
	\label{embeddingimp}
	Let $T>0$. Then it holds
	\begin{align}
		L^\infty_{w*}(0,T;BV(\Omega;\{0,1\})) \cap H^1(0,T;H^1(\om)')\hookrightarrow C([0,T];L^p(\Omega)),
		\label{embedding1}
	\end{align}
	for any $p\in [1,\infty)$.
\end{lemma}
\begin{proof}
Let us assume $	\xhi\in L^\infty_{w*}(0,T;BV(\Omega;\{0,1\}))\cap H^1(0,T;H^1(\om)')$. To prove the result, we first show that $\xhi(t)\in BV(\Omega;\{0,1\})$ for any $t\in[0,T]$, as well as 
	\begin{align}
		\sup_{t\in[0,T]}\norm{\xhi(t)}_{BV(\Omega;\{0,1\})}\leq  2\norm{\xhi}_{L^\infty(0,T;BV(\Omega;\{0,1\}))}.\label{B2}
		\end{align}
		Notice that this is not straightforward, since $BV(\om,\{0,1\})$ is not separable nor reflexive. Now, let us fix $t\in[0,T]$. Since $\xhi\in L^\infty_{w*}(0,T;BV(\Omega;\{0,1\}))$, we can find a sequence $\{t_n\}_{n\in \N}\subset(0,T)$ such that $t_n\to t$ as $n\to \infty$ and
		\begin{align}
			\sup_{n\in\N}\norm{\xhi(t_n)}_{BV(\Omega;\{0,1\})}\leq 2\norm{\xhi}_{L^\infty(0,T;BV(\Omega;\{0,1\}))}.
			\label{bound}
		\end{align}
		Then, recalling the definition of BV functions, we can write, for any $f\in BV(\om;\{0,1\})$, 
		$$
		\norm{\nabla f}_{\mathcal M(\Omega)}=\sup_{\varphi\in C^1_c(\om):\ \norm{\varphi}_{C(\om)}\leq 1}\normmm{\int_\om f\Div\varphi\dx }. 
		$$
		Now, recalling also that $\xhi\in H^1(0,T;H^1(\om)')\hookrightarrow C([0,T]; H^{1}(\Omega)')$, we have 
		\begin{align*}
		&	\normmm{\int_\om \xhi(t)\Div\varphi\dx }=\normmm{\langle \xhi(t),\Div\varphi \rangle_{H^1(\om)',H^1(\om)}}\\&=\lim_{n\to\infty}\normmm{\langle \xhi(t_n),\Div\varphi \rangle_{H^1(\om)',H^1(\om)}}=	\lim_{n\to\infty}\normmm{\int_\om \xhi(t_n)\Div\varphi\dx }\\&
		\leq 2\norm{\xhi}_{L^\infty(0,T;BV(\om;\{0,1\}))}, \quad \forall \varphi \in C_c^1(\om),\quad \norm{\varphi}_{C(\om)}\leq 1.
		\end{align*}
		where we have used \eqref{bound} in the estimate
		\begin{align*}
			\normmm{\int_\om \xhi(t_n)\Div\varphi\dx }\leq \norm{\nabla \xhi(t_n)}_{\mathcal M(\om)}, \quad \forall \varphi \in C_c^1(\om),\quad \norm{\varphi}_{C(\om)}\leq 1,\quad \forall n\in \N.
			\end{align*}
		This entails 
		$$
			\normmm{\int_\om \xhi(t)\Div\varphi\dx }\leq 2\norm{\xhi}_{L^\infty(0,T;BV(\om;\{0,1\}))}, \quad \forall \varphi \in C_c^1(\om),\quad \norm{\varphi}_{C(\om)}\leq 1,
		$$
		so that, taking the supremum over all $\varphi  \in C_c^1(\om)$ such that $\norm{\varphi}_{C(\om)}\leq 1$ we infer that 
		$$
		\norm{\xhi(t)}_{BV(\om;\{0,1\})}\leq 2\norm{\xhi}_{L^\infty(0,T;BV(\om;\{0,1\}))},
		$$
		and since this holds for any $t\in[0,T]$, we have proven \eqref{B2}.
		
		Now we argue as in the well-known Aubin-Lions Lemma, namely, let us fix $t\in[0,T]$ and consider a generic sequence $\{t_n\}_{n\in \N}\subset [0,T]$ such that $t_n\to t$ as $n\to\infty$. Then, recalling the uniform bound \eqref{B2}, together with the embedding $BV(\om)\hookrightarrow\hookrightarrow L^1(\om)$, we infer that there exist a subsequence $\{t_{n_k}\}_{k\in \N}$ and $\wtilde{\xhi}\in BV(\om;\{0,1\})$ such that, as $k\to\infty$, 
	    \begin{align}
	    \xhi(t_{n_k})\to \wtilde\xhi,\quad \text{ strongly in }L^1(\om),
	    		    \end{align} 
		and, by interpolation, since $\xhi\in L^\infty(\om\times(0,T))$, also 
		\begin{align}
		\xhi(t_{n_k})\to \wtilde\xhi,\quad \text{ strongly in }L^p(\om),
		\end{align}
		for any $p\in[1,\infty)$. Now, since it is well known (see, e.g., \cite{Strauss}) that $\xhi\in L^\infty(\om\times(0,T))\cap H^1(0,T;H^1(\om)')\hookrightarrow C_w([0,T];L^2(\om))$, we also infer that 
		$$
		\chi(t_n)\rightharpoonup \xhi(t),\quad \text{weakly in }L^2(\om),
		$$
		as $n\to \infty$, allowing the identification $\wtilde \xhi=\xhi(t)$. By the arbitrarieness of $t\in[0,T]$ and the uniqueness of the limit this then entails that 
		\begin{align*}
			\xhi(s)\to \xhi(t),\quad \text{strongly in }L^p(\om),
			\end{align*}
		as $s\to t$, for any $p\in[1,\infty)$, entailing that $\xhi\in C([0,T];L^p(\om))$ for any $p\in[1,\infty)$ and concluding the proof.
\end{proof}
	\bibliographystyle{siam}
	\bibliography{Bib}
	
\end{document}